\documentclass[10pt, oneside, a4paper]{article}
\usepackage{amsmath, amsthm, amssymb, amsfonts}
\usepackage{geometry}
\usepackage{graphicx}
\usepackage{natbib}
\usepackage{setspace}
\usepackage{algorithm}
\usepackage{algorithmicx}
\usepackage{algpseudocode}
\usepackage{longtable}
\usepackage{caption}

\newtheorem{theorem}{Theorem}[section]
\newtheorem{lemma}[theorem]{Lemma}

\newtheorem{corollary}[theorem]{Corollary}
\theoremstyle{definition}
\newtheorem{definition}{Definition}[section]

\newtheorem{assumption}{Assumption}
\newtheorem{remark}[theorem]{Remark}

\input{tcilatex}

\begin{document}

\title{Further Results on Size and Power of Heteroskedasticity and
Autocorrelation Robust Tests, with an Application to Trend Testing\thanks{%
We thank the referees for helpful comments on a previous version of the
paper. Financial support of the second author by the Danish National
Research Foundation (Grant DNRF 78, CREATES) and by the Program of Concerted
Research Actions (ARC) of the Universit\'{e} libre de Bruxelles is
gratefully acknowledged. Address correspondence to Benedikt P\"{o}tscher,
Department of Statistics, University of Vienna, A-1090 Oskar-Morgenstern
Platz 1. E-Mail: benedikt.poetscher@univie.ac.at.}}
\author{Benedikt M. P\"{o}tscher and David Preinerstorfer \\
Department of Statistics, University of Vienna\\
ECARES, Universit\'{e} libre de Bruxelles}
\date{First version: August 2017\\
Second version: January 2018\\
Third version: April 2019}
\maketitle

\begin{abstract}
We complement the theory developed in \cite{PP2016} with further finite
sample results on size and power of heteroskedasticity and autocorrelation
robust tests. These allow us, in particular, to show that the sufficient
conditions for the existence of size-controlling critical values recently
obtained in \cite{PP3} are often also necessary. We furthermore apply the
results obtained to tests for hypotheses on deterministic trends in
stationary time series regressions, and find that many tests currently used
are strongly size-distorted.
\end{abstract}

\section{Introduction}

Heteroskedasticity and autocorrelation robust tests in regression models
suggested in the literature (e.g., tests based on the covariance estimators
in \cite{NW87, NW94}, \cite{A91}, and \cite{A92}, or tests in \cite{KVB2000}%
, \cite{KiefVogl2002, KV2002, KV2005}) often suffer from substantial size
distortions or power deficiencies. This has been repeatedly documented in
simulation studies, and has been explained analytically by the theory
developed in \cite{PP2016} to a large extent. Given a test for an affine
restriction on the regression coefficient vector, the results in \cite%
{PP2016} provide several sufficient conditions that imply size equal to one,
or severe biasedness of the test (resulting in low power in certain regions
of the alternative). The central object in that theory is the set of
possible covariance matrices of the regression errors, i.e., the covariance
model, and, in particular, its set of concentration spaces. Concentration
spaces are defined as the column spaces of all singular matrices belonging
to the boundary of the covariance model (cf. Definition 2.1 in \cite{PP2016}%
). In \cite{PP2016} it was shown that the position of the concentration
spaces relative to the rejection region of the test often lets one deduce
whether size distortions or power problems occur. Loosely speaking, if a
concentration space lies in the \textquotedblleft
interior\textquotedblright\ of the rejection region, the test has size equal
to one, whereas if a concentration space lies in the \textquotedblleft
exterior\textquotedblright\ (the \textquotedblleft
interior\textquotedblright\ of the complement) of the rejection region, the
test is biased and has nuisance-minimal power equal to zero.\footnote{%
The situation is a bit more complex. For example, sometimes a modification
of the rejection region, which leaves the rejection probabilities unchanged,
is required in order to enforce the interiority (exteriority) condition; see
Theorem 5.7 in \cite{PP2016}.} These interiority (exteriority) conditions
can be formulated in terms of test statistics and critical values, can be
easily checked in practice, and have been made explicit in \cite{PP2016} at
different levels of generality concerning the test statistic and the
covariance model (cf. their Corollary 5.17, Theorem 3.3, Theorem 3.12,
Theorem 3.15, and Theorem 4.2 for more details).

Given a test statistic, the results of \cite{PP2016} just mentioned -- if
applicable -- all lead to implications of the following type: (i) size
equals one for any choice of critical value (e.g., testing a zero
restriction on the mean of a stationary AR(1) time series falls under this
case); or (ii) all critical values smaller than a certain real number
(depending on observable quantities only) lead to a test with size one.
While implication (i) certainly rules out the existence of a
size-controlling critical value, implication (ii) does not, because it only
makes a statement about a certain range of critical values. Hence, the
question when a size-controlling critical value actually exists has not
sufficiently been answered in \cite{PP2016}. Focusing exclusively on size
control, \cite{PP3} recently developed conditions under which size can be
controlled at any level.\footnote{%
We note that, apart from the results mentioned before, \cite{PP2016} also
contains results that ensure size control (and positive infimal power). The
scope of these results is, however, substantially more narrow than the scope
of the results in \cite{PP3}.} It turns out that these conditions can, in
general, not be formulated in terms of concentration spaces of the
covariance model alone. Rather, they are conditions involving a different,
but related, set $\mathbb{J}$, say, of linear spaces obtained from the
covariance model. This set $\mathbb{J}$ consists of nontrivial projections
of concentration spaces as well as of spaces which might be regarded as
\textquotedblleft higher-order\textquotedblright\ concentration spaces (cf.
Section 5 and Appendix B.1 of \cite{PP3} for a detailed discussion). Again,
the conditions in \cite{PP3} do not depend on unobservable quantities, and
hence can be checked by the practitioner. \cite{PP3} also provide algorithms
for the computation of size-controlling critical values, which are
implemented in the \textsf{R}-package \textbf{acrt} (\cite{csart}).

Summarizing we arrive at the following situation: \cite{PP2016} provide --
inter alia -- sufficient conditions for non-existence of size-controlling
critical values in terms of the set of concentration spaces of a covariance
model, whereas \cite{PP3} provide sufficient conditions for the existence of
size-controlling critical values formulated in terms of a different set of
linear spaces derived from the covariance model. Combining the results in 
\cite{PP2016} and \cite{PP3} does in general \textit{not} result in
necessary and sufficient conditions for the existence of size-controlling
critical values. [This is partly due to the fact that different sets of
linear spaces associated with the covariance model are used in these two
papers.] Rather, there remains a range\ of problems for which the existence
of size-controlling critical values can be neither disproved by the results
in \cite{PP2016} nor proved by the results in \cite{PP3}.

In the present paper we close the \textquotedblleft gap\textquotedblright\
between the negative results in \cite{PP2016} on the one hand, and the
positive results in \cite{PP3} on the other hand. We achieve this by
obtaining new negative results that are typically more general than the ones
in \cite{PP2016}. Instead of directly working with concentration spaces of a
given covariance model (as in \cite{PP2016}) our main strategy is
essentially as follows: We first show that size properties of (invariant)
tests are preserved when passing from the given covariance model to a
suitably constructed auxiliary covariance model which has the property that
the concentration spaces of this auxiliary covariance model coincide with
the set $\mathbb{J}$ of linear spaces derived from the initial covariance
model (as used in the results of \cite{PP3}). Then we apply results in \cite%
{PP2016} to the concentration spaces of the auxiliary covariance model to
obtain a necessary condition for the existence of size-controlling critical
values. [This result is first formulated for arbitrary covariance models,
and is then further specialized to the case of stationary autocorrelated
errors.] The so-obtained new result now allows us to prove that the
conditions developed in \cite{PP3} for the possibility of size control are
not only sufficient, but are -- under certain (weak) conditions on the test
statistic -- also necessary. Additionally, we also study power properties
and provide conditions under which a critical value leading to size control
will lead to low power in certain regions of the alternative; we also
discuss conditions under which this is not so.

Obtaining results for the class of problems inaccessible by the results of 
\cite{PP2016} and \cite{PP3} is not only theoretically satisfying. It is
also practically important as this class\ contains empirically relevant
testing problems: As a further contribution we thus apply our results to the
important problem of testing hypotheses on polynomial or cyclical trends in
stationary time series, the former being our main focus. Testing for trends
certainly is an important problem (not only) in economics, and has received
a great amount of attention in the literature. Using our new results we can
prove that many tests currently in use (e.g., conventional tests based on
long-run-variance estimators, or more specialized tests as suggested in \cite%
{vogelsang1998} and \cite{bunzelvogelsang2005}) suffer from severe size
problems whenever the covariance model is not extremely small (that is, is
large enough to contain all covariance matrices of stationary autoregressive
processes of order two or a slight enlargement of that set, a weak condition
that is satisfied by the covariance models used in \cite{vogelsang1998} or 
\cite{bunzelvogelsang2005}; cf.~also the last paragraph preceding Section
5.1.1). Furthermore, our results show that this problem can not be resolved
by increasing the critical values used (as it is established that no
size-controlling critical value exists).

The structure of the article is as follows: Section \ref{frame} introduces
the framework and some notation. In Section \ref{Sec_2} we present results
concerning size properties of nonsphericity-corrected F-type tests. This is
done on two levels of generality: In Subsection \ref{Sec_general_results} we
present results for general covariance models, whereas in Subsection \ref%
{sec_stationary-results} we present results for covariance models obtained
from stationary autocorrelated errors. In these two sections it is also
shown that the conditions for size control obtained in Theorems 3.2, 3.8,
6.5, 6.6 and in Corollary 5.6 of \cite{PP3} are not only sufficient but are
also necessary in important scenarios. In Section \ref{sec_power} we present
results concerning the power of tests based on size-controlling critical
values. Finally, in Section \ref{sec_trends} we discuss consequences of our
results for testing restrictions on coefficients of polynomial and cyclical
regressors. All proofs as well as some auxiliary results are given in the
appendices.

\section{Framework\label{frame}}

\subsection{The model and basic notation}

Consider the linear regression model 
\begin{equation}
\mathbf{Y}=X\beta +\mathbf{U},  \label{lm}
\end{equation}%
where $X$ is a (real) nonstochastic regressor (design) matrix of dimension $%
n\times k$ and where $\beta \in \mathbb{R}^{k}$ denotes the unknown
regression parameter vector. We always assume $\limfunc{rank}(X)=k$ and $%
1\leq k<n$. We furthermore assume that the $n\times 1$ disturbance vector $%
\mathbf{U}=(\mathbf{u}_{1},\ldots ,\mathbf{u}_{n})^{\prime }$ is normally
distributed with mean zero and unknown covariance matrix $\sigma ^{2}\Sigma $%
, where $\Sigma $ varies in a prescribed (nonempty) set $\mathfrak{C}$ of
symmetric and positive definite $n\times n$ matrices and where $0<\sigma
^{2}<\infty $ holds ($\sigma $ always denoting the positive square root).%
\footnote{%
Since we are concerned with finite-sample results only, the elements of $%
\mathbf{Y}$, $X$, and $\mathbf{U}$ (and even the probability space
supporting $\mathbf{Y}$ and $\mathbf{U}$) may depend on sample size $n$, but
this will not be expressed in the notation. Furthermore, the obvious
dependence of$\ \mathfrak{C}$ on $n$ will also not be shown in the notation.}
The set $\mathfrak{C}$ will be referred to as the covariance model. We shall
always assume that $\mathfrak{C}$ allows $\sigma ^{2}$ and $\Sigma $ to be
uniquely determined from $\sigma ^{2}\Sigma $.\footnote{%
That is, $\mathfrak{C}$ has the property that $\Sigma \in \mathfrak{C}$
implies $\delta \Sigma \notin \mathfrak{C}$ for every $\delta \neq 1$.}
[This entails virtually no loss of generality and can always be achieved,
e.g., by imposing some normalization assumption on the elements of $%
\mathfrak{C}$ such as normalizing the first diagonal element of $\Sigma $ or
the norm of $\Sigma $ to one, etc.] The leading case will concern the
situation where $\mathfrak{C}$ results from the assumption that the elements 
$\mathbf{u}_{1},\ldots ,\mathbf{u}_{n}$ of the $n\times 1$ disturbance
vector $\mathbf{U}$ are distributed like consecutive elements of a zero mean
weakly stationary Gaussian process with an unknown spectral density, but
allowing for more general covariance models is useful.

The linear model described in (\ref{lm}) together with the Gaussianity
assumption on $\mathbf{U}$ induces a collection of distributions on the
Borel-sets of $\mathbb{R}^{n}$, the sample space of $\mathbf{Y}$. Denoting a
Gaussian probability measure with mean $\mu \in \mathbb{R}^{n}$ and
(possibly singular) covariance matrix $A$ by $P_{\mu ,A}$, the induced
collection of distributions is then given by 
\begin{equation}
\left\{ P_{\mu ,\sigma ^{2}\Sigma }:\mu \in \mathrm{\limfunc{span}}%
(X),0<\sigma ^{2}<\infty ,\Sigma \in \mathfrak{C}\right\} .  \label{lm2}
\end{equation}%
Since every $\Sigma \in \mathfrak{C}$ is positive definite by assumption,
each element of the set in the previous display is absolutely continuous
with respect to (w.r.t.) Lebesgue measure on $\mathbb{R}^{n}$.

We shall consider the problem of testing a linear (better: affine)
hypothesis on the parameter vector $\beta \in \mathbb{R}^{k}$, i.e., the
problem of testing the null $R\beta =r$ against the alternative $R\beta \neq
r$, where $R$ is a $q\times k$ matrix always of rank $q\geq 1$ and $r\in 
\mathbb{R}^{q}$. Set $\mathfrak{M}=\limfunc{span}(X)$. Define the affine
space 
\begin{equation*}
\mathfrak{M}_{0}=\left\{ \mu \in \mathfrak{M}:\mu =X\beta \text{ and }R\beta
=r\right\}
\end{equation*}%
and let 
\begin{equation*}
\mathfrak{M}_{1}=\left\{ \mu \in \mathfrak{M}:\mu =X\beta \text{ and }R\beta
\neq r\right\} .
\end{equation*}%
Adopting these definitions, the above testing problem can then be written
more precisely as 
\begin{equation}
H_{0}:\mu \in \mathfrak{M}_{0},\ 0<\sigma ^{2}<\infty ,\ \Sigma \in 
\mathfrak{C}\quad \text{ vs. }\quad H_{1}:\mu \in \mathfrak{M}_{1},\
0<\sigma ^{2}<\infty ,\ \Sigma \in \mathfrak{C}.  \label{testing problem}
\end{equation}%
We also define $\mathfrak{M}_{0}^{lin}$ as the linear space parallel to $%
\mathfrak{M}_{0}$, i.e., $\mathfrak{M}_{0}^{lin}=\mathfrak{M}_{0}-\mu _{0}$
for some $\mu _{0}\in \mathfrak{M}_{0}$. Obviously, $\mathfrak{M}_{0}^{lin}$
does not depend on the choice of $\mu _{0}\in \mathfrak{M}_{0}$. The
previously introduced concepts and notation will be used throughout the
paper.

The assumption of Gaussianity is made mainly in order not to obscure the
structure of the problem by technicalities. Substantial generalizations away
from Gaussianity are possible exactly in the same way as the extensions
discussed in Section 5.5 of \cite{PP2016}; see also Appendix E of \cite{PP3}%
. The assumption of nonstochastic regressors can be relaxed somewhat: If $X$
is random and, e.g., independent of $\mathbf{U}$, the results of the paper
apply after one conditions on $X$. For arguments supporting conditional
inference see, e.g., \cite{RO1979}.

We next collect some further terminology and notation used throughout the
paper. A (nonrandomized) \textit{test} is the indicator function of a
Borel-set $W$ in $\mathbb{R}^{n}$, with $W$ called the corresponding \textit{%
rejection region}. The \textit{size} of such a test (rejection region) is
the supremum over all rejection probabilities under the null hypothesis $%
H_{0}$, i.e., 
\begin{equation*}
\sup_{\mu \in \mathfrak{M}_{0}}\sup_{0<\sigma ^{2}<\infty }\sup_{\Sigma \in 
\mathfrak{C}}P_{\mu ,\sigma ^{2}\Sigma }(W).
\end{equation*}%
Throughout the paper we let $\hat{\beta}_{X}(y)=\left( X^{\prime }X\right)
^{-1}X^{\prime }y$, where $X$ is the design matrix appearing in (\ref{lm})
and $y\in \mathbb{R}^{n}$. The corresponding ordinary least squares (OLS)
residual vector is denoted by $\hat{u}_{X}(y)=y-X\hat{\beta}_{X}(y)$. If it
is clear from the context which design matrix is being used, we shall drop
the subscript $X$ from $\hat{\beta}_{X}(y)$ and $\hat{u}_{X}(y)$ and shall
simply write $\hat{\beta}(y)$ and $\hat{u}(y)$. We use $\Pr $ as a generic
symbol for a probability measure. Lebesgue measure on the Borel-sets of $%
\mathbb{R}^{n}$ will be denoted by $\lambda _{\mathbb{R}^{n}}$, whereas
Lebesgue measure on an affine subspace $\mathcal{A}$ of $\mathbb{R}^{n}$
(but viewed as a measure on the Borel-sets of $\mathbb{R}^{n}$) will be
denoted by $\lambda _{\mathcal{A}}$, with zero-dimensional Lebesgue measure
being interpreted as point mass. The set of real matrices of dimension $%
l\times m$ is denoted by $\mathbb{R}^{l\times m}$ (all matrices in the paper
will be real matrices). Let $B^{\prime }$ denote the transpose of a matrix $%
B\in \mathbb{R}^{l\times m}$ and let $\mathrm{\limfunc{span}}(B)$ denote the
subspace in $\mathbb{R}^{l}$ spanned by its columns. For a symmetric and
nonnegative definite matrix $B$ we denote the unique symmetric and
nonnegative definite square root by $B^{1/2}$. For a linear subspace $%
\mathcal{L}$ of $\mathbb{R}^{n}$ we let $\mathcal{L}^{\bot }$ denote its
orthogonal complement and we let $\Pi _{\mathcal{L}}$ denote the orthogonal
projection onto $\mathcal{L}$. For an affine subspace $\mathcal{A}$ of $%
\mathbb{R}^{n}$ we denote by $G(\mathcal{A})$ the group of all affine
transformations on $\mathbb{R}^{n}$ of the form $y\mapsto \delta
(y-a)+a^{\ast }$ where $\delta \neq 0$ and $a$ as well as $a^{\ast }$ belong
to $\mathcal{A}$. [If $\mathcal{A}$ is a linear space, $G(\mathcal{A})$
consists precisely of all transformations of the form $y\mapsto \delta y+%
\bar{a}$ with $\delta \neq 0$ and $\bar{a}\in \mathcal{A}$.] The $j$-th
standard basis vector in $\mathbb{R}^{n}$ is written as $e_{j}(n)$.
Furthermore, we let $\mathbb{N}$ denote the set of all positive integers. A
sum (product, respectively) over an empty index set is to be interpreted as $%
0$ ($1$, respectively). Finally, for a subset $A$ of a topological space we
denote by $\limfunc{cl}(A)$ the closure of $A$ (w.r.t. the ambient space).

\subsection{Classes of test statistics \label{Classes}}

The rejection regions we consider will be of the form $W=\left\{ y\in 
\mathbb{R}^{n}:T(y)\geq C\right\} $, where the critical value $C$ satisfies $%
-\infty <C<\infty $ and the test statistic $T$ is a Borel-measurable
function from $\mathbb{R}^{n}$ to $\mathbb{R}$. With the exception of
Section \ref{sec_power}, the results in the present paper will concern the
class of nonsphericity-corrected F-type test statistics as defined in (28)
of Section 5.4 in \cite{PP2016} that satisfy Assumption 5 in that reference.
For the convenience of the reader we recall the definition of this class of
test statistics. We start with the following assumption, which is Assumption
5 in \cite{PP2016}:

\begin{assumption}
\label{Ass_5}(i) Suppose we have estimators $\check{\beta}:\mathbb{R}%
^{n}\backslash N\rightarrow \mathbb{R}^{k}$ and $\check{\Omega}:\mathbb{R}%
^{n}\backslash N\rightarrow \mathbb{R}^{q\times q}$ that are well-defined
and continuous on $\mathbb{R}^{n}\backslash N$, where $N$ is a closed $%
\lambda _{\mathbb{R}^{n}}$-null set. Furthermore, $\check{\Omega}(y)$ is
symmetric for every $y\in \mathbb{R}^{n}\backslash N$. (ii) The set $\mathbb{%
R}^{n}\backslash N$ is assumed to be invariant under the group $G(\mathfrak{M%
})$, i.e., $y\in \mathbb{R}^{n}\backslash N$ implies $\delta y+X\eta \in 
\mathbb{R}^{n}\backslash N$ for every $\delta \neq 0$ and every $\eta \in 
\mathbb{R}^{k}$. (iii) The estimators satisfy the equivariance properties $%
\check{\beta}(\delta y+X\eta )=\delta \check{\beta}(y)+\eta $ and $\check{%
\Omega}(\delta y+X\eta )=\delta ^{2}\check{\Omega}(y)$ for every $y\in 
\mathbb{R}^{n}\backslash N$, for every $\delta \neq 0$, and for every $\eta
\in \mathbb{R}^{k}$. (iv) $\check{\Omega}$ is $\lambda _{\mathbb{R}^{n}}$%
-almost everywhere nonsingular on $\mathbb{R}^{n}\backslash N$.
\end{assumption}

Nonsphericity-corrected F-type test statistics are now of the form%
\begin{equation}
T(y)=\left\{ 
\begin{array}{cc}
(R\check{\beta}(y)-r)^{\prime }\check{\Omega}^{-1}(y)(R\check{\beta}(y)-r),
& y\in \mathbb{R}^{n}\backslash N^{\ast }, \\ 
0, & y\in N^{\ast },%
\end{array}%
\right.  \label{F-type}
\end{equation}%
where $\check{\beta}$, $\check{\Omega}$, and $N$ satisfy Assumption \ref%
{Ass_5} and where $N^{\ast }=N\cup \left\{ y\in \mathbb{R}^{n}\backslash
N:\det \check{\Omega}(y)=0\right\} $. We recall from Lemmata 5.15 and F.1 in 
\cite{PP2016} that $N^{\ast }$ is then a closed $\lambda _{\mathbb{R}^{n}}$%
-null set that is invariant under $G(\mathfrak{M})$, and that $T$ is
continuous on $\mathbb{R}^{n}\backslash N^{\ast }$ (and is obviously
Borel-measurable on $\mathbb{R}^{n}$). Furthermore, $T$ is $G(\mathfrak{M}%
_{0})$-invariant, i.e., $T(\delta (y-\mu _{0})+\mu _{0}^{\prime })=T(y)$
holds for every $y\in \mathbb{R}^{n}$, every $\delta \neq 0$, every $\mu
_{0}\in \mathfrak{M}_{0}$, and for every $\mu _{0}^{\prime }\in \mathfrak{M}%
_{0}$.

\begin{remark}
\label{rem_subclasses}(\emph{Important subclasses}) (i) Classical
autocorrelation robust test statistics (e.g., those considered in \cite{NW87}%
, \cite{A91} Sections 3-5, or in \cite{KVB2000}, \cite{KiefVogl2002, KV2002,
KV2005}) fall into this class: More precisely, denoting such a test
statistic by $T_{w}$ as in \cite{PP3}, it follows that $T_{w}$ is a
nonsphericity-corrected F-type test statistic with Assumption \ref{Ass_5}
above being satisfied, provided only Assumptions 1 and 2 of \cite{PP3} hold.
Here $\check{\beta}$ is given by the ordinary least squares estimator $\hat{%
\beta}$, $\check{\Omega}$ is given by $\hat{\Omega}_{w}$ defined in Section
3 of \cite{PP3}, and $N=\emptyset $ holds (see Remark 5.17 in \cite{PP3}).
Furthermore, $\check{\Omega}=\hat{\Omega}_{w}$ is then nonnegative definite
on all of $\mathbb{R}^{n}$ (see Section 3.2 of \cite{PP2016} or Section 3 of 
\cite{PP3}). We also recall from Section 5.3 of \cite{PP3} that in this case
the set $N^{\ast }$ can be shown to be a finite union of proper linear
subspaces of $\mathbb{R}^{n}$.

(ii) Classical autocorrelation robust test statistics like $T_{w}$, but
where the weights are now allowed to depend on the data (e.g., through
data-driven bandwidth choice or through prewithening, etc.) as considered,
e.g., in \cite{A91}, \cite{A92}, and \cite{NW94}, also fall into the class
of nonsphericity-corrected F-type tests under appropriate conditions (with
the set $N$ now typically being nonempty), see \cite{Prein2014b} for
details. The same is typically true for test statistics based on parametric
long-run variance estimators or test statistics based on feasible
generalized least squares (cf. Section 3.3 of \cite{PP2016}).

(iii) A statement completely analogous to (i) above applies to the more
general class of test statistics $T_{GQ}$ discussed in Section 3.4B of \cite%
{PP3}, provided Assumption 1 of \cite{PP3} is traded for the assumption that
the weighting matrix $\mathcal{W}_{n}^{\ast }$ appearing in the definition
of $T_{GQ}$ is positive definite (and $\check{\Omega}$ is of course now as
discussed in Section 3.4B of \cite{PP3}); see Remark 5.17 in \cite{PP3}.
Again, $\check{\Omega}$ is then nonnegative definite on all of $\mathbb{R}%
^{n}$ (see Section 3.2.1 of \cite{PP2016}), $N=\emptyset $ holds, and $%
N^{\ast }$ is a finite union of proper linear subspaces of $\mathbb{R}^{n}$
(see Section 5.3 of \cite{PP3}).

(iv) The (weighted) Eicker-test statistic $T_{E,\mathsf{W}}$ (cf.~\cite{E67}%
) as defined on pp.410-411 of \cite{PP3} is also a nonsphericity-corrected
F-type test statistic with Assumption \ref{Ass_5} above being satisfied,
where $\check{\beta}=\hat{\beta}$, $\check{\Omega}=\hat{\Omega}_{E,\mathsf{W}%
}$ defined on p.411 of \cite{PP3}, and $N=\emptyset $ holds. Again, $\check{%
\Omega}$ is nonnegative definite on all of $\mathbb{R}^{n}$, and $N^{\ast }=%
\mathrm{\limfunc{span}}(X)$ holds (see Sections 3 and 5.3 of \cite{PP3}). We
note that the classical (i.e., uncorrected) F-test statistic also falls into
this class as it coincides (up to a known constant) with $T_{E,\mathsf{W}}$
in case $\mathsf{W}$ is the identity matrix.

(v) Under the assumptions of Section 4 of \cite{PP2016} (including
Assumption 3 in that reference), usual heteroskedasticity-robust test
statistics considered in the literature (see \cite{LE2000} for an overview)
also fall into the class of nonsphericity-corrected F-type test statistics
with Assumption \ref{Ass_5} being satisfied. Again, the matrix $\check{\Omega%
}$ is then nonnegative definite everywhere, $N=\emptyset $ holds, and $%
N^{\ast }$ is a finite union of proper linear subspaces of $\mathbb{R}^{n}$
(the latter following from Lemma 4.1 in \cite{PP2016} combined with Lemma
5.18 of \cite{PP3}).
\end{remark}

We shall also encounter cases where $\check{\Omega}(y)$ may not be
nonnegative definite for some values of $y\in \mathbb{R}^{n}\backslash N$.
For these cases the following assumption, which is Assumption 7 in \cite%
{PP2016}, will turn out to be useful. For a discussion of this assumption
see p.~314 of that reference.

\begin{assumption}
\label{Ass_7} For every $v\in \mathbb{R}^{q}$ with $v\neq 0$ we have $%
\lambda _{\mathbb{R}^{n}}\left( \left\{ y\in \mathbb{R}^{n}\backslash
N^{\ast }:v^{\prime }\check{\Omega}^{-1}(y)v=0\right\} \right) =0$.
\end{assumption}

\section{Results on the size of nonsphericity-corrected F-type test
statistics\label{Sec_2}}

\subsection{A result for general covariance models\label{Sec_general_results}%
}

In this subsection we start with a negative result concerning the size of a
class of nonsphericity-corrected F-type test statistics that is central to
many of the results in the present paper. In particular, it allows us to
show that the sufficient conditions for size control obtained in \cite{PP3}
are often also necessary. The result complements negative results in \cite%
{PP2016} and is obtained by combining Lemmata \ref{associated_cov_model} and %
\ref{change cov-model} in Appendix \ref{App A} with Corollary 5.17 of \cite%
{PP2016}. Its relationship to negative results in \cite{PP2016} is further
discussed in Appendix \ref{sec_comments}. We recall the following definition
from \cite{PP3}.

\begin{definition}
Given a linear subspace $\mathcal{L}$ of $\mathbb{R}^{n}$ with $\dim (%
\mathcal{L})<n$ and a covariance model $\mathfrak{C}$, we let $\mathcal{L}(%
\mathfrak{C})=\left\{ \mathcal{L}(\Sigma ):\Sigma \in \mathfrak{C}\right\} $%
, where $\mathcal{L}(\Sigma )=\Pi _{\mathcal{L}^{\bot }}\Sigma \Pi _{%
\mathcal{L}^{\bot }}/\Vert {\Pi _{\mathcal{L}^{\bot }}\Sigma \Pi _{\mathcal{L%
}^{\bot }}}\Vert $. Furthermore, we define 
\begin{equation*}
\mathbb{J}(\mathcal{L},\mathfrak{C})=\left\{ \mathrm{\limfunc{span}}(\bar{%
\Sigma}):\bar{\Sigma}\in \limfunc{cl}(\mathcal{L}(\mathfrak{C})),\ \limfunc{%
rank}(\bar{\Sigma})<n-\dim (\mathcal{L})\right\} ,
\end{equation*}%
where the closure is here understood w.r.t. $\mathbb{R}^{n\times n}$. [The
symbol $\Vert {\cdot }\Vert $ here denotes a norm on $\mathbb{R}^{n\times n}$%
. Note that $\mathbb{J}(\mathcal{L},\mathfrak{C})$ does not depend on which
norm is chosen.]
\end{definition}

The space $\mathcal{L}$ figuring in this definition will always be an
appropriately chosen subspace related to invariance properties of the tests
under consideration. A leading case is when $\mathcal{L}=\mathfrak{M}%
_{0}^{lin}$. Loosely speaking, the linear spaces belonging to $\mathbb{J}(%
\mathcal{L},\mathfrak{C})$ are either (nontrivial) projections of
concentration spaces of the covariance model $\mathfrak{C}$ (in the sense of 
\cite{PP2016}) on $\mathcal{L}^{\bot }$, or are what one could call
\textquotedblleft higher-order\textquotedblright\ concentration spaces. For
a more detailed discussion see Appendix B.1 of \cite{PP3}.

\begin{theorem}
\label{size_one}Let $\mathfrak{C}$ be a covariance model. Let $T$ be a
nonsphericity-corrected F-type test statistic of the form (\ref{F-type})
based on $\check{\beta}$ and $\check{\Omega}$ satisfying Assumption \ref%
{Ass_5} with $N=\emptyset $. Furthermore, assume that $\check{\Omega}(y)$ is
nonnegative definite for every $y\in \mathbb{R}^{n}$. If an $\mathcal{S}\in 
\mathbb{J}(\mathfrak{M}_{0}^{lin},\mathfrak{C})$ satisfying $\mathcal{S}%
\subseteq \mathrm{\limfunc{span}}(X)$ exists, then 
\begin{equation}
\sup_{\Sigma \in \mathfrak{C}}P_{\mu _{0},\sigma ^{2}\Sigma }(T\geq C)=1
\label{eq_size_one}
\end{equation}%
holds for every critical value $C$, $-\infty <C<\infty $, for every $\mu
_{0}\in \mathfrak{M}_{0}$, and for every $\sigma ^{2}\in (0,\infty )$.
\end{theorem}

\begin{remark}
\emph{(Extensions)} (i) As noted in Section \ref{Classes}, any $T$ as in the
theorem is $G(\mathfrak{M}_{0})$-invariant. In some cases $T$ and its
associated set $N^{\ast }$ are additionally invariant w.r.t. addition of
elements from a linear space $\mathcal{V}\subseteq \mathbb{R}^{n}$. In such
a case $\mathcal{L}=\mathrm{\limfunc{span}}(\mathfrak{M}_{0}^{lin}\cup 
\mathcal{V})$ necessarily has dimension less than $n-1<n$, and the variant
of Theorem \ref{size_one} where $\mathbb{J}(\mathfrak{M}_{0}^{lin},\mathfrak{%
C})$ is replaced by $\mathbb{J}(\mathcal{L},\mathfrak{C})$ also holds.%
\footnote{\label{dimL}That $\dim (\mathcal{L})<n-1$ must hold is seen as
follows: Suppose $\dim (\mathcal{L})\geq n-1$. Then $T$ is $\lambda _{%
\mathbb{R}^{n}}$-almost everywhere constant (this is trivial if $\dim (%
\mathcal{L})=n$ and follows from Remark 5.14(i) in \cite{PP3} in case $\dim (%
\mathcal{L})=n-1$). However, this contradicts Part 2 of Lemma 5.16 of \cite%
{PP3}.}

(ii) A result similar to Theorem \ref{size_one}, operating under a weaker
condition than $\mathcal{S}\subseteq \mathrm{\limfunc{span}}(X)$ for some $%
\mathcal{S}\in \mathbb{J}(\mathfrak{M}_{0}^{lin},\mathfrak{C})$, is given in
Theorem \ref{size_one_extension} in Appendix \ref{App A}. This result also
allows for $N\neq \emptyset $, but is restricted to the case where $q$, the
number of restrictions tested, is equal to $1$ and where $\check{\beta}$ is
the least squares estimator in (\ref{lm}).
\end{remark}

The preceding theorem can now be used to show that the conditions for size
control obtained in Corollary 5.6 (and Remark 5.8) of \cite{PP3} are not
only sufficient, but are actually necessary, in some important scenarios.
This is formulated in the subsequent corollary; see also Remark \ref{nec_1}
below. [We note that $T$ in this corollary satisfies the assumptions of
Corollary 5.6 of \cite{PP3} (with $N^{\dag }=N^{\ast }$ and $\mathcal{V}%
=\{0\}$) in view of Lemma 5.16 in the same reference.]

\begin{corollary}
\label{suff&nec}Let $\mathfrak{C}$ be a covariance model. Let $T$ be a
nonsphericity-corrected F-type test statistic of the form (\ref{F-type})
based on $\check{\beta}$ and $\check{\Omega}$ satisfying Assumption \ref%
{Ass_5} with $N=\emptyset $. Furthermore, assume that $\check{\Omega}(y)$ is
nonnegative definite for every $y\in \mathbb{R}^{n}$, and that $N^{\ast }=%
\mathrm{\limfunc{span}}(X)$. Then $\mathcal{S}\nsubseteq \mathrm{\limfunc{%
span}}(X)$ for every $\mathcal{S}\in \mathbb{J}(\mathfrak{M}_{0}^{lin},%
\mathfrak{C})$ is necessary and sufficient for size-controllability (at any
significance level $\alpha \in (0,1)$), i.e., is necessary and sufficient
for the fact that for every $\alpha \in (0,1)$ there exists a real number $%
C(\alpha )$ such that%
\begin{equation}
\sup_{\mu _{0}\in \mathfrak{M}_{0}}\sup_{0<\sigma ^{2}<\infty }\sup_{\Sigma
\in \mathfrak{C}}P_{\mu _{0},\sigma ^{2}\Sigma }(T\geq C(\alpha ))\leq \alpha
\label{size_control}
\end{equation}%
holds.\footnote{%
For conditions under which a smallest size-controlling critical value exists
and when equality can be achieved in (\ref{size_control}) see \cite{PP3},
Section 5.2.}
\end{corollary}

\begin{remark}
\label{nec_1}\emph{(Special cases)} (i) Corollary \ref{suff&nec} applies, in
particular, to the (weighted) Eicker-test statistic $T_{E,\mathsf{W}}$ in
view of Remark \ref{rem_subclasses}(iv) above. Note that $N^{\ast }=\mathrm{%
\limfunc{span}}(X)$ is here always satisfied. By Remark \ref{rem_subclasses}%
(iv), Corollary \ref{suff&nec} also applies to the classical F-test
statistic.

(ii) Next consider the classical autocorrelation robust test statistic $%
T_{w} $ with Assumptions 1 and 2 of \cite{PP3} being satisfied. Then
Corollary \ref{suff&nec} also applies to $T_{w}$ in view of Remark \ref%
{rem_subclasses}(i) above, provided $N^{\ast }=\mathrm{\limfunc{span}}(X)$
holds. While the relation $N^{\ast }=\mathrm{\limfunc{span}}(X)$ need not
always hold for $T_{w}$ (see the discussion in Section 5.3 of \cite{PP3}),
it holds for many combinations of restriction matrix $R$ and design matrix $%
X $ (in fact, it holds generically in many universes of design matrices as a
consequence of Lemma A.3 in Appendix A of \cite{PP3}). Hence, for such
combinations of $R$ and $X$, Corollary \ref{suff&nec} applies to $T_{w}$.

(iii) For test statistics $T_{GQ}$ with positive definite weighting matrix $%
\mathcal{W}_{n}^{\ast }$ a statement completely analogous to (ii) above
holds in view of Remark \ref{rem_subclasses}(iii). The same is true for
heteroskedasticity-robust test statistics as discussed in Remark \ref%
{rem_subclasses}(v).
\end{remark}

\begin{remark}
While Theorem \ref{size_one} applies to any combination of test statistic $T$
and covariance model $\mathfrak{C}$ as long as they satisfy the assumptions
of the theorem, in a typical application the choice of the test statistic
used will certainly be dictated by properties of the covariance model $%
\mathfrak{C}$ one maintains. For example, in case $\mathfrak{C}$ models
stationary autocorrelated errors different test statistics will be employed
than in the case where $\mathfrak{C}$ models heteroskedasticity.
\end{remark}

\subsection{Results for covariance models obtained from stationary
autocorrelated errors\label{sec_stationary-results}}

We next specialize the results of the preceding section to the case of
stationary autocorrelated errors. i.e., to the case where the elements $%
\mathbf{u}_{1},\ldots ,\mathbf{u}_{n}$ of the $n\times 1$ disturbance vector 
$\mathbf{U}$ in model (\ref{lm}) are distributed like consecutive elements
of a zero mean weakly stationary Gaussian process with an unknown spectral
density, which is not almost everywhere equal to zero. Consequently, the
covariance matrix of the disturbance vector is positive definite and can be
written as $\sigma ^{2}\Sigma (f)$ where%
\begin{equation*}
\Sigma (f)=\left[ \int_{-\pi }^{\pi }e^{-\iota (j-l)\omega }f(\omega
)d\omega \right] _{j,l=1}^{n},
\end{equation*}%
with $f$ varying in $\mathfrak{F}$, a prescribed (nonempty) family of \emph{%
normalized} (i.e., $\int_{-\pi }^{\pi }f(\omega )d\omega =1$) spectral
densities, and where $0<\sigma ^{2}<\infty $ holds. Here $\iota $ denotes
the imaginary unit. We define the associated covariance model via $\mathfrak{%
C}(\mathfrak{F})=\left\{ \Sigma (f):f\in \mathfrak{F}\right\} $. Examples
for the set $\mathfrak{F}$ are (i) $\mathfrak{F}_{\mathrm{all}}$, the set of 
\emph{all} normalized spectral densities, or (ii) $\mathfrak{F}_{\mathrm{%
ARMA(}p,q\mathrm{)}}$, the set of all normalized spectral densities
corresponding to stationary autoregressive moving average models of order at
most $(p,q)$, or (iii) the set of normalized spectral densities
corresponding to (stationary) fractional autoregressive moving average
models, etc. We shall write $\mathfrak{F}_{\mathrm{AR(}p\mathrm{)}}$ for $%
\mathfrak{F}_{\mathrm{ARMA(}p,0\mathrm{)}}$.

We need to recall some more concepts and notation from \cite{PP3}; for
background see this reference. Let $\omega \in \lbrack 0,\pi ]$ and let $%
s\geq 0$ be an integer. Define $E_{n,s}(\omega )$ as the $n\times 2$%
-dimensional matrix with $j$-th row equal to $\left( j^{s}\cos (j\omega
),j^{s}\sin (j\omega )\right) $. Given a linear subspace $\mathcal{L}$ of $%
\mathbb{R}^{n}$ with $\limfunc{dim}(\mathcal{L})<n$, define for every $%
\omega \in \lbrack 0,\pi ]$%
\begin{equation}
\rho (\omega ,\mathcal{L})=\min \left\{ s\in \mathbb{N}\cup \{0\}:\limfunc{%
span}(E_{n,s}(\omega ))\nsubseteq \mathcal{L}\right\} .  \label{rho}
\end{equation}%
As discussed in Section 3.1 of \cite{PP3}, the set on the r.h.s. of (\ref%
{rho}) is nonempty for every $\omega \in \lbrack 0,\pi ]$. Thus $\rho $ is
well-defined and takes values in $\mathbb{N}\cup \{0\}$. Furthermore, $\rho
(\omega ,\mathcal{L})>0$ holds at most for finitely many $\omega \in \lbrack
0,\pi ]$ as shown in the same reference. We now define $\underline{\omega }(%
\mathcal{L})$ as the vector obtained by ordering the elements of $\{\omega
\in \lbrack 0,\pi ]:\rho (\omega ,\mathcal{L})>0\}$ from smallest to
largest, provided this set is nonempty, and we denote by $p(\mathcal{L})$
the dimension of this vector; furthermore, we set $d_{i}(\mathcal{L})=\rho
(\omega _{i}(\mathcal{L}),\mathcal{L})$ for every $i=1,\ldots ,p(\mathcal{L}%
) $, where $\omega _{i}(\mathcal{L})$ denotes the $i$-the coordinate of $%
\underline{\omega }(\mathcal{L})$, and we write $\underline{d}(\mathcal{L})$
for the vector with $i$-th coordinate equal to $d_{i}(\mathcal{L})$. If the
set $\{\omega \in \lbrack 0,\pi ]:\rho (\omega ,\mathcal{L})>0\}$ is empty,
we take $\underline{\omega }(\mathcal{L})$ as well as $\underline{d}(%
\mathcal{L})$ as the $0$-tuple and set $p(\mathcal{L})=0$. As in \cite{PP3},
for $d$ a natural number we define $\kappa (\omega ,d)=2d$ for $\omega \in
(0,\pi )$ and $\kappa (\omega ,d)=d$ for $\omega \in \left\{ 0,\pi \right\} $%
. Furthermore, we set $\kappa (\underline{\omega }(\mathcal{L}),\underline{d}%
(\mathcal{L}))=\sum \kappa (\omega _{i}(\mathcal{L}),d_{i}(\mathcal{L}))$
where the sum extends over $i=1,\ldots ,p(\mathcal{L})$, with the convention
that this sum is zero if $p(\mathcal{L})=0$. For ease of notation we shall
often simply write $\rho (\gamma )$ for $\rho (\gamma ,\mathfrak{M}%
_{0}^{lin})$.

The subsequent theorem specializes Theorem \ref{size_one} to the case where $%
\mathfrak{C}=\mathfrak{C}(\mathfrak{F})$. For a definition of the collection 
$\mathbb{S}(\mathfrak{F},\mathcal{L})$ of certain subsets of $[0,\pi ]$
figuring in this theorem see Definition 6.4 of \cite{PP3}.

\begin{theorem}
\label{size_one_auto}Let $\mathfrak{F}$ be a nonempty set of normalized
spectral densities, i.e., $\emptyset \neq \mathfrak{F}\subseteq \mathfrak{F}%
_{\mathrm{all}}$. Let $T$ be a nonsphericity-corrected F-type test statistic
of the form (\ref{F-type}) based on $\check{\beta}$ and $\check{\Omega}$
satisfying Assumption \ref{Ass_5} with $N=\emptyset $. Furthermore, assume
that $\check{\Omega}(y)$ is nonnegative definite for every $y\in \mathbb{R}%
^{n}$. Suppose there exists a linear subspace $\mathcal{S}$ of $\mathbb{R}%
^{n}$ that can be written as%
\begin{equation}
\mathcal{S}=\limfunc{span}\left( \Pi _{(\mathfrak{M}_{0}^{lin})^{\bot
}}\left( E_{n,\rho (\gamma _{1})}(\gamma _{1}),\ldots ,E_{n,\rho (\gamma
_{p})}(\gamma _{p})\right) \right) \text{ \ \ \ for some \ \ }\Gamma \in 
\mathbb{S}(\mathfrak{F},\mathfrak{M}_{0}^{lin}),  \label{aspt:S}
\end{equation}%
where the $\gamma _{i}$'s denote the elements of $\Gamma $ and $p=\limfunc{%
card}(\Gamma )$, such that $\mathcal{S}$ satisfies $\mathcal{S}\subseteq 
\mathrm{\limfunc{span}}(X)$ (or, equivalently, $\limfunc{span}(E_{n,\rho
(\gamma _{1})}(\gamma _{1}),\ldots ,E_{n,\rho (\gamma _{p})}(\gamma
_{p}))\subseteq \mathrm{\limfunc{span}}(X)$). Then $\limfunc{dim}(\mathcal{S}%
)<n-\limfunc{dim}(\mathfrak{M}_{0}^{lin})$ holds. Furthermore, 
\begin{equation*}
\sup_{f\in \mathfrak{F}}P_{\mu _{0},\sigma ^{2}\Sigma (f)}(T\geq C)=1
\end{equation*}%
holds for every critical value $C$, $-\infty <C<\infty $, for every $\mu
_{0}\in \mathfrak{M}_{0}$, and for every $\sigma ^{2}\in (0,\infty )$.
\end{theorem}

\begin{remark}
\label{rem_simpli}Suppose $\mathfrak{F}$ in Theorem \ref{size_one_auto} has
the property that $\gamma \in \bigcup \mathbb{S}(\mathfrak{F},\mathfrak{M}%
_{0}^{lin})$ implies $\{\gamma \}\in \mathbb{S}(\mathfrak{F},\mathfrak{M}%
_{0}^{lin})$ (as is, e.g., the case if $\mathfrak{F}\supseteq \mathfrak{F}_{%
\mathrm{AR(}2\mathrm{)}}$, cf. Lemma \ref{lemma:AR2} below). Then it is easy
to see that the set $\Gamma $ in the theorem can be chosen to be a singleton.
\end{remark}

This theorem is applicable to any nonempty set $\mathfrak{F}$ of normalized
spectral densities. In case more is known about the richness of $\mathfrak{F}
$, the sufficient condition in the preceding result can sometimes be
simplified substantially. Below we present such a result making use of the
subsequent lemma.

\begin{lemma}
\label{lemma:AR2}Let $\mathfrak{F}\subseteq \mathfrak{F}_{\mathrm{all}}$
satisfy $\mathfrak{F}\supseteq \mathfrak{F}_{\mathrm{AR(}2\mathrm{)}}$ and
let $\mathcal{L}$ be a linear subspace of $\mathbb{R}^{n}$ with $\limfunc{dim%
}(\mathcal{L})<n$. Let $\gamma \in \lbrack 0,\pi ]$. Then $\{\gamma \}\in 
\mathbb{S}(\mathfrak{F},\mathcal{L})$ if and only if $\kappa (\underline{%
\omega }(\mathcal{L}),\underline{d}(\mathcal{L}))+\kappa (\gamma ,1)<n$. And 
$\{\gamma \}\in \mathbb{S}(\mathfrak{F},\mathcal{L})$ holds for every $%
\gamma \in \lbrack 0,\pi ]$ if and only if $\kappa (\underline{\omega }(%
\mathcal{L}),\underline{d}(\mathcal{L}))+2<n$. Furthermore, $\gamma \in
\bigcup \mathbb{S}(\mathfrak{F},\mathcal{L})$ if and only if $\{\gamma \}\in 
\mathbb{S}(\mathfrak{F},\mathcal{L})$.
\end{lemma}

\begin{remark}
\label{rem:kappa and dim}(i) A sufficient condition for $\kappa (\underline{%
\omega }(\mathcal{L}),\underline{d}(\mathcal{L}))+\kappa (\gamma ,1)<n$ ($%
\kappa (\underline{\omega }(\mathcal{L}),\underline{d}(\mathcal{L}))+2<n$,
respectively) is given by $\limfunc{dim}(\mathcal{L})+\kappa (\gamma ,1)<n$ (%
$\limfunc{dim}(\mathcal{L})+2<n$, respectively). This follows from $\kappa (%
\underline{\omega }(\mathcal{L}),\underline{d}(\mathcal{L}))\leq \limfunc{dim%
}(\mathcal{L})$ established in Lemma D.1 in Appendix D of \cite{PP3}.

(ii) In the case $\mathcal{L}=\mathfrak{M}_{0}^{lin}$ the latter two
conditions become $k-q+\kappa (\gamma ,1)<n$ and $k-q+2<n$, respectively.
Note that the condition $k-q+\kappa (\gamma ,1)<n$ is always satisfied for $%
\gamma =0$ or $\gamma =\pi $ (as then $\kappa (\gamma ,1)=1$). For $\gamma
\in (0,\pi )$ this condition coincides with $k-q+2<n$, and is always
satisfied except if $k=n-1$ and $q=1$.
\end{remark}

Armed with the preceding lemma we can now establish the following
consequence of Theorem \ref{size_one_auto} provided $\mathfrak{F}$ is rich
enough to encompass $\mathfrak{F}_{\mathrm{AR(}2\mathrm{)}}$, which clearly
is a very weak condition in the context of autocorrelation robust testing.%
\footnote{%
Recall that a premise of autocorrelation robust testing is agnosticism about
the correlation structure of the error process.}

\begin{theorem}
\label{size_one_AR2}Let $\mathfrak{F}\subseteq \mathfrak{F}_{\mathrm{all}}$
satisfy $\mathfrak{F}\supseteq \mathfrak{F}_{\mathrm{AR(}2\mathrm{)}}$. Let $%
T$ be a nonsphericity-corrected F-type test statistic of the form (\ref%
{F-type}) based on $\check{\beta}$ and $\check{\Omega}$ satisfying
Assumption \ref{Ass_5} with $N=\emptyset $. Furthermore, assume that $\check{%
\Omega}(y)$ is nonnegative definite for every $y\in \mathbb{R}^{n}$. Suppose
there exists a $\gamma \in \lbrack 0,\pi ]$ such that $\limfunc{span}%
(E_{n,\rho (\gamma )}(\gamma ))\subseteq \mathrm{\limfunc{span}}(X)$. Then $%
\kappa (\underline{\omega }(\mathfrak{M}_{0}^{lin}),\underline{d}(\mathfrak{M%
}_{0}^{lin}))+\kappa (\gamma ,1)<n$ holds, and we have 
\begin{equation}
\sup_{f\in \mathfrak{F}}P_{\mu _{0},\sigma ^{2}\Sigma (f)}(T\geq C)=1
\label{size}
\end{equation}%
for every critical value $C$, $-\infty <C<\infty $, for every $\mu _{0}\in 
\mathfrak{M}_{0}$, and for every $\sigma ^{2}\in (0,\infty )$.
\end{theorem}

\begin{remark}
\emph{(Further comments on the necessity of the sufficient conditions for
size control in \cite{PP3})} (i) Suppose $T$ is as in Theorem \ref%
{size_one_auto}, additionally satisfying $N^{\ast }=\mathrm{\limfunc{span}}%
(X)$. Theorem \ref{size_one_auto} then shows that the sufficient conditions
for size control given in Part 1 of Theorem 6.5 in \cite{PP3} (or the
equivalent formulation given in Part 2 of that theorem) is also necessary.

(ii) Suppose $T$ is as in (i) and assume furthermore that $\mathfrak{F}$ is
as in Remark \ref{rem_simpli}. Then also the sufficient condition for size
control \textquotedblleft $\limfunc{span}(E_{n,\rho (\gamma )}(\gamma
))\nsubseteq \mathrm{\limfunc{span}}(X)$ for every $\gamma \in \bigcup 
\mathbb{S}(\mathfrak{F},\mathfrak{M}_{0}^{lin})$\textquotedblright\
mentioned in Part 2 of Theorem 6.5 of \cite{PP3} is necessary. [This is seen
as follows: Suppose not, i.e., $\limfunc{span}(E_{n,\rho (\gamma )}(\gamma
))\subseteq \mathrm{\limfunc{span}}(X)$ holds for some $\gamma \in \bigcup 
\mathbb{S}(\mathfrak{F},\mathfrak{M}_{0}^{lin})$. Now apply Theorem \ref%
{size_one_auto} with $\Gamma =\{\gamma \}$, which is possible because of
Remark \ref{rem_simpli}, resulting in size being equal to one, a
contradiction.]

(iii) Suppose $T$ is as in (i) and assume that $\mathfrak{F}\subseteq 
\mathfrak{F}_{\mathrm{all}}$ satisfies $\mathfrak{F}\supseteq \mathfrak{F}_{%
\mathrm{AR(}2\mathrm{)}}$. Then $\mathfrak{F}$ satisfies the property in
Remark \ref{rem_simpli} in view of Lemma \ref{lemma:AR2}, and thus (ii)
above applies. In this situation even more is true in view of Theorem \ref%
{size_one_AR2}: The further sufficient condition for size control
\textquotedblleft $\limfunc{span}(E_{n,\rho (\gamma )}(\gamma ))\nsubseteq 
\mathrm{\limfunc{span}}(X)$ for every $\gamma \in \lbrack 0,\pi ]$%
\textquotedblright\ given in Part 2 of Theorem 6.5 of \cite{PP3} is in fact
also necessary.

(iv) The discussion in (i)-(iii) covers (weighted) Eicker-test statistics $%
T_{E,\mathsf{W}}$ (including the classical F-test statistic) as well as
classical autocorrelation robust test statistics $T_{w}$ (the latter under
Assumptions 1 and 2 of \cite{PP3} and if $N^{\ast }=\mathrm{\limfunc{span}}%
(X)$ holds); it also covers the test statistics $T_{GQ}$ (provided the
weighting matrix $\mathcal{W}_{n}^{\ast }$ is positive definite and $N^{\ast
}=\mathrm{\limfunc{span}}(X)$ holds). In particular, the discussion in
(i)-(iii) thus applies to the sufficient conditions given in Theorem 6.6 in 
\cite{PP3} and its variants outlined in Remark 6.8 of that reference.
Furthermore, it transpires from this discussion that the sufficient
conditions for size control provided in Theorem 3.8 of \cite{PP3} are
actually necessary; and the same is true for Theorem 3.2 in that reference
(provided the set $\mathsf{B}$ given there coincides with $\mathrm{\limfunc{%
span}}(X)$).\footnote{%
Note that $\mathfrak{F}=\mathfrak{F}_{\mathrm{all}}$ in those two theorems.}
\end{remark}

The results so far have only concerned the size of nonsphericity-corrected
F-type test statistics for which the exceptional set $N$ is empty and $%
\check{\Omega}$ is nonnegative definite everywhere. We now provide a result
also for the case where this condition is not met.\footnote{%
Theorem \ref{size_one_extension} in Appendix \ref{App A} also allows for $%
N\neq \emptyset $, but requires $\check{\Omega}(y)$ to be nonnegative
definite for every $y\in \mathbb{R}^{n}\backslash N$ (implying that $\check{%
\Omega}$ is nonnegative definite $\lambda _{\mathbb{R}^{n}}$-a.e.). This
result also contains further assumptions such as $q=1$.}

\begin{definition}
Let $\mathfrak{F}_{\mathrm{AR(}2\mathrm{)}}^{ext}$ denote the set of all
normalized spectral densities of the form $c_{1}f+(2\pi )^{-1}c_{2}$ with $%
f\in \mathfrak{F}_{\mathrm{AR(}2\mathrm{)}}$ and $c_{1}+c_{2}=1$, $c_{1}\geq
0$, $c_{2}\geq 0$.
\end{definition}

Obviously, $\mathfrak{F}_{\mathrm{AR(}2\mathrm{)}}\subseteq \mathfrak{F}_{%
\mathrm{AR(}2\mathrm{)}}^{ext}\subseteq \mathfrak{F}_{\mathrm{ARMA(}2,2%
\mathrm{)}}$ holds. While the preceding result maintained that $\mathfrak{F}$
contains $\mathfrak{F}_{\mathrm{AR(}2\mathrm{)}}$, the next result maintains
the slightly stronger condition that $\mathfrak{F}\supseteq \mathfrak{F}_{%
\mathrm{AR(}2\mathrm{)}}^{ext}$.

\begin{theorem}
\label{theo:AR2+}Let $\mathfrak{F}\subseteq \mathfrak{F}_{\mathrm{all}}$
satisfy $\mathfrak{F}\supseteq \mathfrak{F}_{\mathrm{AR(}2\mathrm{)}}^{ext}$%
. Let $T$ be a nonsphericity-corrected F-type test statistic of the form (%
\ref{F-type}) based on $\check{\beta}$ and $\check{\Omega}$ satisfying
Assumption \ref{Ass_5}. Furthermore, assume that $\check{\Omega}$ also
satisfies Assumption \ref{Ass_7}. Suppose there exists a $\gamma \in \lbrack
0,\pi ]$ such that $\limfunc{span}(E_{n,\rho (\gamma )}(\gamma ))\subseteq 
\mathrm{\limfunc{span}}(X)$. Then for every critical value $C$, $-\infty
<C<\infty $, for every $\mu _{0}\in \mathfrak{M}_{0}$, and for every $\sigma
^{2}\in (0,\infty )$ it holds that 
\begin{equation}
P_{0,I_{n}}(\check{\Omega}\text{ is nonnegative definite})\leq K(\gamma
)\leq \sup_{f\in \mathfrak{F}}P_{\mu _{0},\sigma ^{2}\Sigma (f)}\left( T\geq
C\right) ,  \label{claim}
\end{equation}%
where $K(\gamma )$ is defined by%
\begin{equation*}
K(\gamma )=\int \Pr \left( \bar{\xi}_{\gamma }(x)\geq 0\right)
dP_{0,I_{\kappa (\gamma ,1)}}(x)
\end{equation*}%
with the random variable $\bar{\xi}_{\gamma }(x)$ given by 
\begin{equation*}
\bar{\xi}_{\gamma }(x)=(R\hat{\beta}_{X}(\bar{E}_{n,\rho (\gamma )}(\gamma
)x))^{\prime }\check{\Omega}^{-1}\left( \mathbf{G}\right) R\hat{\beta}_{X}(%
\bar{E}_{n,\rho (\gamma )}(\gamma )x)
\end{equation*}%
on the event where $\left\{ \mathbf{G}\in \mathbb{R}^{n}\backslash N^{\ast
}\right\} $ and by $\bar{\xi}_{\gamma }(x)=0$ otherwise. Here $\mathbf{G}$
is a standard normal $n$-vector, $\bar{E}_{n,\rho (\gamma )}(\gamma
)=E_{n,\rho (\gamma )}(\gamma )$ if $\gamma \in (0,\pi )$ and $\bar{E}%
_{n,\rho (\gamma )}(\gamma )$ denotes the first column of $E_{n,\rho (\gamma
)}(\gamma )$ otherwise. [Recall that $\hat{\beta}_{X}(y)=(X^{\prime
}X)^{-1}X^{\prime }y$.]
\end{theorem}

The significance of the preceding theorem is that it provides a lower bound
for the size of a large class of nonsphericity-corrected F-type tests,
including those with $N\neq \emptyset $ or with $\check{\Omega}$ not
necessarily nonnegative definite. In particular, it shows that size can not
be controlled at a given desired significance level $\alpha $, if $\alpha $
is below the threshold given by the lower bound in (\ref{claim}). Observe
that this threshold will typically be close to $1$, at least if $n$ is
sufficiently large, since (possibly after rescaling) $\check{\Omega}$ will
often approach a positive definite matrix as $n\rightarrow \infty $.

\begin{remark}
(i) There are at most finitely many $\gamma $ satisfying the assumption $%
\limfunc{span}(E_{n,\rho (\gamma )}(\gamma ))\subseteq \mathrm{\limfunc{span}%
}(X)$ in the preceding theorem. To see this note that any such $\gamma $
must coincide with a coordinate of $\underline{\omega }(\limfunc{span}(X))$
(since trivially $\limfunc{span}(E_{n,0}(\gamma ))\subseteq \limfunc{span}%
(X) $ in case $\rho (\gamma )=0$ by this assumption, and since $\limfunc{span%
}(E_{n,0}(\gamma ))\subseteq \mathfrak{M}_{0}^{lin}\subseteq \limfunc{span}%
(X) $ in case $\rho (\gamma )>0$), and that the dimension of the vector $%
\underline{\omega }(\limfunc{span}(X))$ is finite since $\rho (\omega ,%
\mathrm{\limfunc{span}}(X))>0$ can hold at most for finitely many $\omega $%
's as discussed subsequent to (\ref{rho}).

(ii) If $\digamma $ denotes the (finite) set of $\gamma $'s satisfying the
assumption $\limfunc{span}(E_{n,\rho (\gamma )}(\gamma ))\subseteq \mathrm{%
\limfunc{span}}(X)$ in the theorem, relation (\ref{claim}) in fact implies%
\begin{equation*}
P_{0,I_{n}}(\check{\Omega}\text{ is nonnegative definite})\leq \min_{\gamma
\in \digamma }K(\gamma )\leq \max_{\gamma \in \digamma }K(\gamma )\leq
\sup_{f\in \mathfrak{F}}P_{\mu _{0},\sigma ^{2}\Sigma (f)}\left( T\geq
C\right) .
\end{equation*}

(iii) Similar to Theorem \ref{size_one_AR2}, Theorem \ref{theo:AR2+} also
delivers (\ref{size}) in case $\check{\Omega}$ is nonnegative definite $%
\lambda _{\mathbb{R}^{n}}$-almost everywhere. However, note that the latter
theorem imposes a stronger condition on the set $\mathfrak{F}$.
\end{remark}

\begin{remark}
\emph{(Extensions) }Suppose $T$ is as in Theorem \ref{theo:AR2+}. If $T$ and
its associated set $N^{\ast }$ are not only $G(\mathfrak{M}_{0})$-invariant,
but are additionally invariant w.r.t. addition of elements from a linear
space $\mathcal{V}\subseteq \mathbb{R}^{n}$, then the variant of Theorem \ref%
{theo:AR2+}, where $\mathcal{L}$ replaces $\mathfrak{M}_{0}^{lin}$ and $\rho
(\gamma ,\mathcal{L})$ replaces $\rho (\gamma )$, can be seen to hold.
\end{remark}

\begin{remark}
Some results in this section are formulated for sets of spectral densities $%
\mathfrak{F}$ satisfying $\mathfrak{F}\supseteq \mathfrak{F}_{\mathrm{AR(}2%
\mathrm{)}}$ or $\mathfrak{F}\supseteq \mathfrak{F}_{\mathrm{AR(}2\mathrm{)}%
}^{ext}$, and thus for covariance models $\mathfrak{C}(\mathfrak{F})$
satisfying $\mathfrak{C}(\mathfrak{F})\supseteq \mathfrak{C}(\mathfrak{F}_{%
\mathrm{AR(}2\mathrm{)}})$ or $\mathfrak{C}(\mathfrak{F})\supseteq \mathfrak{%
C}(\mathfrak{F}_{\mathrm{AR(}2\mathrm{)}}^{ext})$, respectively. Trivially,
these results also hold for any covariance model $\mathfrak{C}$ (not
necessarily of the form $\mathfrak{C}(\mathfrak{F})$) that satisfies $%
\mathfrak{C}\supseteq \mathfrak{C}(\mathfrak{F}_{\mathrm{AR(}2\mathrm{)}})$
or $\mathfrak{C}\supseteq \mathfrak{C}(\mathfrak{F}_{\mathrm{AR(}2\mathrm{)}%
}^{ext})$, respectively. This observation also applies to other results in
this paper further below and will not be repeated.
\end{remark}

\section{Results concerning power\label{sec_power}}

We now show for a large class of test statistics, even larger than the class
of nonsphericity-corrected F-type test statistics, that -- under certain
conditions -- a choice of critical value leading to size less than one
necessarily implies that the test is severely biased and thus has bad power
properties in certain regions of the alternative hypothesis (cf. Part 3 of
Theorem 5.7 and Remark 5.5(iii) in \cite{PP2016}). The relevant conditions
essentially say that a collection $\mathbb{K}$ as in the subsequent lemma
can be found that is nonempty. It should be noted, however, that there are
important instances where (i) the relevant conditions are not satisfied
(that is, a nonempty $\mathbb{K}$ satisfying the properties required in the
lemma does not exist) and (ii) small size and good power properties coexist.
For results in that direction see Theorems 3.7, 5.10, 5.12, and 5.21 in \cite%
{PP2016} as well as Proposition 5.2 and Theorem 5.4 in \cite{Prein2014b}.

The subsequent lemma is a variant of Lemma 5.11 in \cite{PP3}. Recall that $%
\mathbb{H}$, defined in that lemma, certainly contains all one-dimensional $%
\mathcal{S}\in \mathbb{J}(\mathcal{L},\mathfrak{C})$ (provided such elements
exist).

\begin{lemma}
\label{improved_L_5.10}Let $\mathfrak{C}$ be a covariance model. Assume that
the test statistic $T:\mathbb{R}^{n}\rightarrow \mathbb{R}$ is
Borel-measurable and is continuous on the complement of a closed set $%
N^{\dag }$. Assume that $T$ and $N^{\dag }$ are $G(\mathfrak{M}_{0})$%
-invariant, and are also invariant w.r.t. addition of elements of a linear
subspace $\mathcal{V}$ of $\mathbb{R}^{n}$. Define $\mathcal{L}=\limfunc{span%
}(\mathfrak{M}_{0}^{lin}\cup \mathcal{V})$ and assume that $\dim \mathcal{L}%
<n$. Let $\mathbb{H}$ and $C(\mathcal{S})$ be defined as in Lemma 5.11 of 
\cite{PP3}. Let $\mathbb{K}$ be a subset of $\mathbb{H}$ and define $C_{\ast
}(\mathbb{K})=\inf_{\mathcal{S}\in \mathbb{K}}C(\mathcal{S})$ and $C^{\ast }(%
\mathbb{K})=\sup_{\mathcal{S}\in \mathbb{K}}C(\mathcal{S})$, with the
convention that $C_{\ast }(\mathbb{K})=\infty $ and $C^{\ast }(\mathbb{K}%
)=-\infty $ if $\mathbb{K}$ is empty. Suppose that $\mathbb{K}$ has the
property that for every $\mathcal{S}\in \mathbb{K}$ the set $N^{\dag }$ is a 
$\lambda _{\mu _{0}+\mathcal{S}}$-null set for some $\mu _{0}\in \mathfrak{M}%
_{0}$ (and hence for all $\mu _{0}\in \mathfrak{M}_{0})$. Then the following
holds:

\begin{enumerate}
\item For every $C\in (-\infty ,C^{\ast }(\mathbb{K}))$, every $\mu _{0}\in 
\mathfrak{M}_{0}$, and every $\sigma ^{2}\in (0,\infty )$ we have 
\begin{equation*}
\sup_{\Sigma \in \mathfrak{C}}P_{\mu _{0},\sigma ^{2}\Sigma }(T\geq C)=1.
\end{equation*}

\item For every $C\in (C_{\ast }(\mathbb{K}),\infty )$, every $\mu _{0}\in 
\mathfrak{M}_{0}$, and every $\sigma ^{2}\in (0,\infty )$ we have 
\begin{equation*}
\inf_{\Sigma \in \mathfrak{C}}P_{\mu _{0},\sigma ^{2}\Sigma }(T\geq C)=0.
\end{equation*}
\end{enumerate}
\end{lemma}

Part 1 of the lemma implies that the size of the test equals $1$ if $%
C<C^{\ast }(\mathbb{K})$. Part 2 shows that the test is severely biased for $%
C>C_{\ast }(\mathbb{K})$, which -- in view of the invariance properties of $%
T $ (cf. Part 3 of Theorem 5.7 and Remark 5.5(iii) in \cite{PP2016}) --
implies bad power properties such as (\ref{eqn:inf1}) and (\ref{eqn:inf2})
below. In particular, Part 2 implies that infimal power is zero for such
choices of $C$. [Needless to say, the lemma neither implies that $%
\sup_{\Sigma \in \mathfrak{C}}P_{\mu _{0},\sigma ^{2}\Sigma }(T\geq C)$ is
less than $1$ for $C>C^{\ast }(\mathbb{K})$ nor that $\inf_{\Sigma \in 
\mathfrak{C}}P_{\mu _{0},\sigma ^{2}\Sigma }(T\geq C)$ is positive for $%
C<C_{\ast }(\mathbb{K})$. For conditions implying that size is less than $1$
for appropriate choices of $C$ see \cite{PP3}.] The computation of the
constants $C^{\ast }(\mathbb{K})$ and $C_{\ast }(\mathbb{K})$ can sometimes
be simplified, see Lemma \ref{lem:simplify} in Appendix \ref{App C}. Before
proceeding, we want to note that the preceding lemma also provides a
negative size result (namely that the test based on $T$ has size equal to $1$
for \emph{every} $C$), if $C^{\ast }(\mathbb{K})=\infty $ holds for a
collection $\mathbb{K}$ satisfying the assumptions of that lemma.

The announced theorem is now as follows and builds on the preceding lemma.

\begin{theorem}
\label{neg_power}Let $\mathfrak{C}$ be a covariance model. Assume that the
test statistic $T:\mathbb{R}^{n}\rightarrow \mathbb{R}$ is Borel-measurable
and is continuous on the complement of a closed set $N^{\dag }$. Assume that 
$T$ and $N^{\dag }$ are $G(\mathfrak{M}_{0})$-invariant, and are also
invariant w.r.t. addition of elements of a linear subspace $\mathcal{V}$ of $%
\mathbb{R}^{n}$. Define $\mathcal{L}=\limfunc{span}(\mathfrak{M}%
_{0}^{lin}\cup \mathcal{V})$ and assume that $\dim \mathcal{L}<n$. Then the
following hold:

\begin{enumerate}
\item Suppose there exist two elements $\mathcal{S}_{1}$ and $\mathcal{S}%
_{2} $ of $\mathbb{H}$ such that $C(\mathcal{S}_{1})\neq C(\mathcal{S}_{2})$%
. Suppose further that for $i=1,2$ the set $N^{\dag }$ is a $\lambda _{\mu
_{0}+\mathcal{S}_{i}}$-null set for some $\mu _{0}\in \mathfrak{M}_{0}$ (and
hence for all $\mu _{0}\in \mathfrak{M}_{0})$. Then for any critical value $%
C $, $-\infty <C<\infty $, satisfying\footnote{\label{footnote_3}Because of $%
G(\mathfrak{M}_{0})$-invariance (cf. Remark 5.5(iii) in \cite{PP2016}), the
left-hand side of (\ref{eqn:sup0}) coincides with $\sup_{\Sigma \in 
\mathfrak{C}}P_{\mu _{0},\sigma ^{2}\Sigma }(T\geq C)$ for any $\mu _{0}\in 
\mathfrak{M}_{0}$ and any $\sigma ^{2}\in (0,\infty ).$ Similarly, the
left-hand side of (\ref{eqn:inf0}) coincides with $\inf_{\Sigma \in 
\mathfrak{C}}P_{\mu _{0},\sigma ^{2}\Sigma }(T\geq C)$ for any $\mu _{0}\in 
\mathfrak{M}_{0}$ and any $\sigma ^{2}\in (0,\infty ).$} 
\begin{equation}
\sup_{\mu _{0}\in \mathfrak{M}_{0}}\sup_{0<\sigma ^{2}<\infty }\sup_{\Sigma
\in \mathfrak{C}}P_{\mu _{0},\sigma ^{2}\Sigma }(T\geq C)<1,
\label{eqn:sup0}
\end{equation}%
we have 
\begin{equation}
\inf_{\mu _{0}\in \mathfrak{M}_{0}}\inf_{0<\sigma ^{2}<\infty }\inf_{\Sigma
\in \mathfrak{C}}P_{\mu _{0},\sigma ^{2}\Sigma }(T\geq C)=0.
\label{eqn:inf0}
\end{equation}

\item Suppose there exists an element $\mathcal{S}$ of $\mathbb{H}$ such
that $N^{\dag }$ is a $\lambda _{\mu _{0}+\mathcal{S}}$-null set for some $%
\mu _{0}\in \mathfrak{M}_{0}$ (and hence for all $\mu _{0}\in \mathfrak{M}%
_{0})$. Then (\ref{eqn:sup0}) implies that $C\geq C(\mathcal{S})$ must hold;
furthermore, (\ref{eqn:sup0}) implies (\ref{eqn:inf0}), except possibly if $%
C=C(\mathcal{S})$ holds.

\item Suppose (\ref{eqn:inf0}) holds for some $C$, $-\infty <C<\infty $. Then%
\begin{equation}
\inf_{0<\sigma ^{2}<\infty }\inf_{\Sigma \in \mathfrak{C}}P_{\mu _{1},\sigma
^{2}\Sigma }(T\geq C)=0  \label{eqn:inf1}
\end{equation}%
for every $\mu _{1}\in \mathfrak{M}_{1}$, and 
\begin{equation}
\inf_{\mu _{1}\in \mathfrak{M}_{1}}\inf_{\Sigma \in \mathfrak{C}}P_{\mu
_{1},\sigma ^{2}\Sigma }(T\geq C)=0  \label{eqn:inf2}
\end{equation}%
for every $\sigma ^{2}\in (0,\infty )$.
\end{enumerate}
\end{theorem}

In the important special case where $\mathcal{V}=\{0\}$, the assumptions on $%
T$ and the associated set $N^{\dag }$ in the second and third sentence of
the preceding theorem are satisfied, e.g., for nonsphericity-corrected
F-type test statistics (under Assumption \ref{Ass_5}), including the test
statistics $T_{w}$, $T_{GQ}$, and $T_{E,\mathsf{W}}$ given in Section \ref%
{Classes} above; see also Section 5.3 in \cite{PP3}. Furthermore, for the
class of test statistics $T$ such that Theorem \ref{size_one} applies (and
for which $N^{\dag }=N^{\ast }=\limfunc{span}(X)$ holds), it can be shown
that $N^{\dag }$ is a $\lambda _{\mu _{0}+\mathcal{S}}$-null set for any $%
\mathcal{S}\in \mathbb{H}$ (in fact, for any $\mathcal{S}\in \mathbb{J}(%
\mathcal{L},\mathfrak{C})$) provided (\ref{eqn:sup0}) holds. These
observations lead to the following corollary.

\begin{corollary}
\label{cor_neg_power}Let $\mathfrak{C}$ be a covariance model and let $T$ be
a nonsphericity-corrected F-type test statistic of the form (\ref{F-type})
based on $\check{\beta}$ and $\check{\Omega}$ satisfying Assumption \ref%
{Ass_5} with $N=\emptyset $. Furthermore, assume that $\check{\Omega}(y)$ is
nonnegative definite for every $y\in \mathbb{R}^{n}$ and that $N^{\ast }=%
\limfunc{span}(X)$.

\begin{enumerate}
\item Suppose there exist two elements $\mathcal{S}_{1}$ and $\mathcal{S}%
_{2} $ of $\mathbb{H}$ (where $\mathbb{H}$ is as in Theorem \ref{neg_power}
with $\mathcal{V}=\{0\}$) such that $C(\mathcal{S}_{1})\neq C(\mathcal{S}%
_{2})$. If a critical value $C$, $-\infty <C<\infty $, satisfies (\ref%
{eqn:sup0}), then it also satisfies (\ref{eqn:inf0}); and thus it also
satisfies (\ref{eqn:inf1}) and (\ref{eqn:inf2}).

\item Suppose that $\mathbb{H}$ is nonempty (where $\mathbb{H}$ is as in
Theorem \ref{neg_power} with $\mathcal{V}=\{0\}$) but $C(\mathcal{S})$ is
the same for all $\mathcal{S}\in \mathbb{H}$. Then (\ref{eqn:sup0}) implies
that $C\geq C(\mathcal{S})$ must hold; furthermore, (\ref{eqn:sup0}) implies
(\ref{eqn:inf0}) (and thus (\ref{eqn:inf1}) and (\ref{eqn:inf2})), except
possibly if $C=C(\mathcal{S})$ holds.
\end{enumerate}
\end{corollary}

Theorem \ref{neg_power} as well as the preceding corollary maintain
conditions that, in particular, require $\mathbb{H}$ to be nonempty. In view
of Lemma 5.11 in \cite{PP3}, $\mathbb{H}$ is certainly nonempty if a
one-dimensional $\mathcal{S}\in \mathbb{J}(\mathcal{L},\mathfrak{C})$
exists. The following lemma shows that for $\mathfrak{C}=\mathfrak{C}(%
\mathfrak{F})$ with $\mathfrak{F}\supseteq \mathfrak{F}_{AR(2)}$ this is
indeed the case; in fact, for such $\mathfrak{C}$ typically at least two
such spaces exist.\footnote{%
While the one-dimensional spaces given in the lemma typically will be
different, it is \emph{not} established in the lemma that this is
necessarily always the case.}

\begin{lemma}
\label{lemma:eleH}Let $\mathfrak{F}\subseteq \mathfrak{F}_{\mathrm{all}}$
satisfy $\mathfrak{F}\supseteq \mathfrak{F}_{\mathrm{AR(}2\mathrm{)}}$. Let $%
\mathcal{L}$ be a linear subspace of $\mathbb{R}^{n}$ satisfying $\dim (%
\mathcal{L})+1<n$. Then, for $\gamma \in \{0,\pi \}$, $\limfunc{span}\left(
\Pi _{\mathcal{L}^{\bot }}\left( E_{n,\rho (\gamma ,\mathcal{L})}(\gamma
)\right) \right) $ belongs to $\mathbb{J}(\mathcal{L},\mathfrak{C}(\mathfrak{%
F}))$ and is one-dimensional.
\end{lemma}

The preceding lemma continues to hold for any covariance model $\mathfrak{C}%
\supseteq \mathfrak{C}(\mathfrak{F}_{\mathrm{AR(}2\mathrm{)}})$ in a trivial
way, since $\mathbb{J}(\mathcal{L},\mathfrak{C})\supseteq \mathbb{J}(%
\mathcal{L},\mathfrak{C}(\mathfrak{F}_{\mathrm{AR(}2\mathrm{)}}))$ then
certainly holds. Also note that the condition $\dim (\mathcal{L})+1<n$ is
always satisfied in the important special case where $\mathcal{L}=\mathfrak{M%
}_{0}^{lin}$, since $\dim (\mathfrak{M}_{0}^{lin})=k-q<n-1$.

\section{Consequences for testing hypotheses on deterministic trends\label%
{sec_trends}}

In this section we discuss important consequences of the results obtained so
far for testing restrictions on coefficients of polynomial and cyclical
regressors when the errors are stationary, more precisely, have a covariance
model of the form $\mathfrak{C}(\mathfrak{F})$. Such testing problems have,
for obvious reasons, received a great deal of attention in econometrics, and
are relevant in many other fields such as, e.g., climate or ecological
research.\footnote{%
See, e.g., \cite{bence1995}, who finds substantial undercoverage of
confidence intervals derived from several tests corrected for
autocorrelation.} In particular, we show that a large class of
nonsphericity-corrected F-type test statistics leads to unsatisfactory test
procedures in this context. In Subsection \ref{subs:poly} we present results
concerning hypotheses on the coefficients of polynomial regressors. Results
concerning tests for hypotheses on the coefficients of cyclical regressors
are briefly discussed in Subsection \ref{subs:cycl}.

\subsection{Polynomial regressors\label{subs:poly}}

We consider here the case where one tests hypotheses that involve the
coefficient of a polynomial regressor as expressed in the subsequent
assumption:

\begin{assumption}
\label{as:poly} Suppose that $X=(F,\tilde{X})$, where $F$ is an $n\times
k_{F}$-dimensional matrix $(1\leq k_{F}\leq k$), the $j$-th column being
given by $(1^{j-1},\ldots ,n^{j-1})^{\prime }$, and where $\tilde{X}$ is an $%
n\times (k-k_{F})$-dimensional matrix such that $X$ has rank $k$ (here $%
\tilde{X}$ is the empty matrix if $k_{F}=k$). Furthermore, suppose that the
restriction matrix $R$ has a nonzero column $R_{\cdot i}$ for some $%
i=1,\ldots ,k_{F}$, i.e., the hypothesis involves coefficients of the
polynomial trend.
\end{assumption}

Under this assumption one obtains the subsequent theorem as a consequence of
Theorem \ref{size_one_AR2}.

\begin{theorem}
\label{size_one_polynomial}Let $\mathfrak{F}\subseteq \mathfrak{F}_{\mathrm{%
all}}$ satisfy $\mathfrak{F}\supseteq \mathfrak{F}_{\mathrm{AR(}2\mathrm{)}}$%
. Suppose that Assumption \ref{as:poly} holds. Let $T$ be a
nonsphericity-corrected F-type test statistic of the form (\ref{F-type})
based on $\check{\beta}$ and $\check{\Omega}$ satisfying Assumption \ref%
{Ass_5} with $N=\emptyset $. Furthermore, assume that $\check{\Omega}(y)$ is
nonnegative definite for every $y\in \mathbb{R}^{n}$. Then 
\begin{equation*}
\sup_{f\in \mathfrak{F}}P_{\mu _{0},\sigma ^{2}\Sigma (f)}(T\geq C)=1
\end{equation*}%
holds for every critical value $C$, $-\infty <C<\infty $, for every $\mu
_{0}\in \mathfrak{M}_{0}$, and for every $\sigma ^{2}\in (0,\infty )$.
\end{theorem}

The previous theorem relies in particular on the assumption that $%
N=\emptyset $ and that $\check{\Omega}$ is nonnegative definite everywhere.
While these two assumptions may appear fairly natural and are widely
satisfied, e.g., for the test statistics $T_{w}$, $T_{GQ}$, and $T_{E,%
\mathsf{W}}$ as discussed in Remark \ref{rem_subclasses}, we shall see in
Subsections \ref{subss:trad} and \ref{subss:rec} below that they are not
satisfied by some tests suggested in the literature. To obtain results also
for tests that are not covered by the previous theorem we can apply Theorem %
\ref{theo:AR2+}. The following result is then obtained.

\begin{theorem}
\label{size_bound_polynomial}Let $\mathfrak{F}\subseteq \mathfrak{F}_{%
\mathrm{all}}$ satisfy $\mathfrak{F}\supseteq \mathfrak{F}_{\mathrm{AR(}2%
\mathrm{)}}^{ext}$. Suppose that Assumption \ref{as:poly} holds. Let $T$ be
a nonsphericity-corrected F-type test statistic of the form (\ref{F-type})
based on $\check{\beta}$ and $\check{\Omega}$ satisfying Assumption \ref%
{Ass_5}. Furthermore, assume that $\check{\Omega}$ also satisfies Assumption %
\ref{Ass_7}. Then for every critical value $C$, $-\infty <C<\infty $, for
every $\mu _{0}\in \mathfrak{M}_{0}$, and for every $\sigma ^{2}\in
(0,\infty )$ it holds that%
\begin{equation}
P_{0,I_{n}}(\check{\Omega}\text{ is nonnegative definite})\leq
P_{0,I_{n}}(R_{\cdot i_{0}}^{\prime }\check{\Omega}^{-1}R_{\cdot i_{0}}\geq
0)\leq \sup_{f\in \mathfrak{F}}P_{\mu _{0},\sigma ^{2}\Sigma (f)}\left(
T\geq C\right) ,  \label{eq_lower_bound}
\end{equation}%
where $R_{\cdot i_{0}}$ denotes the first nonzero column of $R$. [Note that $%
\check{\Omega}$ is $P_{0,I_{n}}$-almost everywhere nonsingular in view of
Assumption \ref{Ass_5}.]
\end{theorem}

Theorem \ref{size_bound_polynomial} shows that under Assumption \ref{as:poly}
a large class of nonsphericity-corrected F-type tests, including cases with $%
N\neq \emptyset $ or with $N=\emptyset $ but where $\check{\Omega}$ is not
necessarily nonnegative definite everywhere, typically have large size. In
particular, size can not be controlled at a given desired significance level 
$\alpha $, if $\alpha $ is below the lower bound in (\ref{eq_lower_bound}).
Observe that this lower bound will typically be close to $1$, at least if $n$
is sufficiently large.

\begin{remark}
\label{rem:AR(1)}(i) In the special case where Assumption \ref{as:poly} is
satisfied with $R_{\cdot 1}\neq 0$, Theorem \ref{size_one_polynomial}
continues to hold even under the weaker assumption that only $\mathfrak{F}%
\supseteq \mathfrak{F}_{\mathrm{AR(}1\mathrm{)}}$ holds.\footnote{%
In fact, it holds more generally for any covariance model $\mathfrak{C}$
that has $\limfunc{span}(e_{+})$ as a concentration space in the sense of 
\cite{PP2016}.} This follows from Part 3 of Corollary 5.17 in \cite{PP2016}
upon noting that $\mathcal{Z}=\limfunc{span}(e_{+})$ is a concentration
space of $\mathfrak{C}(\mathfrak{F})$ by Lemma G.1 in the same reference,
that $\check{\Omega}$ vanishes on $\limfunc{span}(X)\supseteq \mathcal{Z}$
as a consequence of the assumption $N=\emptyset $ (see the discussion
following (27) in \cite{PP2016}), and that $R\check{\beta}(\lambda
e_{+})=\lambda R_{\cdot 1}\neq 0$ for all $\lambda \neq 0$.\footnote{\label%
{nezero}To see that $R\check{\beta}(\lambda e_{+})=\lambda R_{\cdot 1}$,
note that $\lambda e_{+}$ is of the form $X\gamma $ with $\gamma =\lambda
e_{1}(k)$, since $e_{+}$ is the first column of $X$. The equivariance
property of $\check{\beta}$ in Assumption \ref{Ass_5} gives $\check{\beta}%
(X\gamma )=\check{\beta}(0)+\gamma $ as well as $\check{\beta}(0)=\check{%
\beta}(\alpha 0)=\alpha \check{\beta}(0)$ for every $\alpha \neq 0$. This
implies $\check{\beta}(0)=0$, and hence $\check{\beta}(X\gamma )=\gamma $.}
Here $e_{+}$ denotes the $n\times 1$ vector of ones.

(ii) In the special case where Assumption \ref{as:poly} is satisfied with $%
R_{\cdot 1}\neq 0$, also Theorem \ref{size_bound_polynomial} continues to
hold under the weaker assumption that $\mathfrak{F}\supseteq \mathfrak{F}_{%
\mathrm{AR(}1\mathrm{)}}$ holds, provided the identity matrix $I_{n}$
appearing in (\ref{eq_lower_bound}) is replaced by the nonsingular matrix $%
\Phi (0)=e_{+}e_{+}^{\prime }+D(0)$, where $D(0)$ is the matrix $D$ given in
Part 3 of Lemma G.1 in \cite{PP2016}. This follows from Remark \ref%
{rem:final}(iii) further below, upon noting that the situation considered
here can be viewed as a special case of the situation described in Remark %
\ref{rem:final}(iii) with $\omega =0$.
\end{remark}

To illustrate the scope and applicability of Theorems \ref%
{size_one_polynomial} and \ref{size_bound_polynomial} above (beyond the test
statistics such as $T_{w}$, $T_{GQ}$, and $T_{E,\mathsf{W}}$ mentioned
before), we shall now apply them to some commonly used test statistics that
have been designed for testing polynomial trends. First, in Subsection \ref%
{subss:trad}, we shall derive properties of conventional tests for
polynomial trends. Such tests are based on long-run-variance estimators and
classical results due to \cite{grenander1954}. In Subsection \ref{subss:rec}
we shall discuss properties of tests that have been introduced more recently
by \cite{vogelsang1998} and \cite{bunzelvogelsang2005}. While our discussion
of methods is certainly not exhaustive (for example, we do not discuss tests
in \cite{Harvey2007} or \cite{Perron2009}, which have been suggested only
for the special case of testing a restriction on the slope in a
\textquotedblleft linear trend plus noise model\textquotedblright ), it
should also serve the purpose of presenting a general pattern how one can
check the reliability of polynomial trend tests. It might also help to avoid
pitfalls in the construction of novel tests for polynomial trends.

Before we proceed to a discussion of properties of specific tests, we would
like to emphasize the following: in the present section we provide, for some
commonly used tests, results on their maximal rejection probability over 
\begin{equation*}
\{P_{\mu _{0},\sigma ^{2}\Sigma (f)}:f\in \mathcal{F}\}
\end{equation*}%
for every $\mu _{0}\in \mathfrak{M}_{0}$ and every $\sigma ^{2}\in (0,\infty
)$. We establish these results under the weak assumption that $\mathfrak{F}$
contains at least $\mathfrak{F}_{\mathrm{AR(}2\mathrm{)}}$ or the slight
enlargement $\mathfrak{F}_{\mathrm{AR(}2\mathrm{)}}^{ext}\subseteq \mathfrak{%
F}_{\mathrm{ARMA(}2,2\mathrm{)}}$. The recent trend testing literature,
cf.~in particular Section 3.1 in \cite{vogelsang1998} and Assumption 1 in 
\cite{bunzelvogelsang2005}, studies tests for models induced by \emph{all}
regression errors $\mathbf{u}_{t}$ satisfying 
\begin{equation*}
\mathbf{u}_{t}=\delta \mathbf{u}_{t-1}+\mathbf{w}_{t},\quad t=2,\ldots
,n,\quad \mathbf{u}_{1}=\mathbf{w}_{1}\text{ \ \ }(\text{or }\mathbf{u}%
_{1}=\sum\nolimits_{j=0}^{\lfloor \tau n\rfloor }\delta ^{j}\mathbf{w}%
_{1-j}).
\end{equation*}%
Here $\delta \in (-1,1]$ is an additional unknown parameter and $\mathbf{w}%
_{t}$ is a weakly stationary linear process with martingale difference
innovations that have uniformly bounded fourth moments and conditional
variance $1$, and with coefficients $d_{i}$ for $i\in \mathbb{N\cup \{}0%
\mathbb{\}}$ satisfying $\sum_{i=0}^{\infty }d_{i}\neq 0$ and the
summability condition $\sum_{i=0}^{\infty }i|d_{i}|<\infty $. Also the
coefficients $d_{i}$ are unknown parameters. Obviously, the assumptions on
the innovations are satisfied for an i.i.d. sequence of standard normal
random variables. Hence, setting $\delta =0$ in the previous displayed
equation, we see that the model considered in \cite{vogelsang1998} or \cite%
{bunzelvogelsang2005} \emph{contains}, in particular, 
\begin{equation*}
\{P_{\mu _{0},\sigma ^{2}\Sigma (f)}:\mu _{0}\in \mathfrak{M}_{0},\sigma
^{2}\in (0,\infty ),f\in \mathfrak{F}_{\mathrm{ARMA(}2,2\mathrm{)}}\}.
\end{equation*}%
As a consequence, any lower bound for size obtained in our context for sets $%
\mathfrak{F}$ required only to satisfy $\mathfrak{F}\supseteq \mathfrak{F}_{%
\mathrm{AR(}2\mathrm{)}}$ (or $\mathfrak{F}\supseteq \mathfrak{F}_{\mathrm{%
AR(}2\mathrm{)}}^{ext}$) a fortiori provides a lower bound for the size in
the setting considered in \cite{vogelsang1998} and \cite{bunzelvogelsang2005}
(since $\mathfrak{F}_{\mathrm{AR(}2\mathrm{)}}\subseteq \mathfrak{F}_{%
\mathrm{AR(}2\mathrm{)}}^{ext}\subseteq \mathfrak{F}_{\mathrm{ARMA(}2,2%
\mathrm{)}}$).

\subsubsection{Properties of conventional tests for hypotheses on polynomial
trends\label{subss:trad}}

The structure of tests that have traditionally been used for testing
restrictions on coefficients of polynomial trends (i.e., when the design
matrix $X$ satisfies Assumption \ref{as:poly}, and in particular if $k_{F}=k$%
) is motivated by results concerning the asymptotic covariance matrix of the
OLS estimator (and its efficiency) in regression models with stationary
error processes and deterministic polynomial time trends by \cite%
{grenander1954} (cf. also the discussion in \cite{bunzelvogelsang2005} on
p.~383). The corresponding test statistics are nonsphericity-corrected
F-type test statistics as in (\ref{F-type}). They are based on the OLS
estimator $\hat{\beta}(=\hat{\beta}_{X})$ and a covariance matrix estimator 
\begin{equation}
\check{\Omega}_{\mathcal{W}}(y)=\hat{\omega}_{\mathcal{W}}(y)R(X^{\prime
}X)^{-1}R^{\prime }.  \label{eqn:covest}
\end{equation}%
Here the \textquotedblleft long-run-variance estimator\textquotedblright\ $%
\hat{\omega}_{_{\mathcal{W}}}$ is of the form 
\begin{equation}
\hat{\omega}_{\mathcal{W}}(y)=n^{-1}\hat{u}^{\prime }(y)\mathcal{W}(y)\hat{u}%
(y),  \label{eqn:lrve}
\end{equation}%
where $\mathcal{W}(y)$ is a symmetric, possibly data-dependent, $n\times n$%
-dimensional matrix that may not be well-defined on all of $\mathbb{R}^{n}$.%
\footnote{%
The matrix $\mathcal{W}$ may depend on $n$, a dependence not shown in the
notation. Furthermore, assuming symmetry of $\mathcal{W}$ entails no loss of
generality, since given a long-run-variance-estimator as in (\ref{eqn:lrve})
and based on a non-symmetric weights matrix $\mathcal{W}_{\ast }$, one can
always pass to an equivalent long-run-variance estimator by replacing $%
\mathcal{W}_{\ast }$ with the symmetric matrix $\mathcal{W}=(\mathcal{W}%
_{\ast }+\mathcal{W}_{\ast }^{\prime })/2$.} In many cases, however, $%
\mathcal{W}$ is constant, i.e., does not depend on $y$, and is also positive
definite. For example, this is so in the leading case where the $(i,j)$-th
element of $\mathcal{W}$ is of the form $\kappa (|i-j|/M)$ for some
(deterministic) $M>0$ (typically depending on $n$) and a kernel function $%
\kappa $ such as the Bartlett, Parzen, Quadratic-Spectral, or Daniell kernel
(positive definiteness does not hold, e.g., for the rectangular kernel with $%
M>1$). Note that in case $\mathcal{W}$ is given by a kernel $\kappa $ the
estimator $\hat{\omega}_{\mathcal{W}}$ in the previous display can be
written in the more familiar form 
\begin{equation*}
\hat{\omega}_{\mathcal{W}}(y)=\sum_{i=-(n-1)}^{n-1}\kappa (|i|/M)\hat{\gamma}%
_{i}(y),
\end{equation*}%
where $\hat{\gamma}_{i}(y)=\hat{\gamma}_{-i}(y)=n^{-1}\sum_{j=i+1}^{n}\hat{u}%
_{j}(y)\hat{u}_{j-i}(y)$ for $i\geq 0$. For trend tests based on the OLS
estimator $\hat{\beta}$ and a covariance estimator $\check{\Omega}_{\mathcal{%
W}}$ as in (\ref{eqn:covest}) we shall first obtain two corollaries from
Theorems \ref{size_one_polynomial} and \ref{size_bound_polynomial} that
cover the case where $\mathcal{W}$ is constant.\footnote{\label{ftnt_1}The
slightly more general case, where $\mathcal{W}$ is not constant in $y$ (and
is defined on all of $\mathbb{R}^{n}$) but $\mathcal{W}^{\ast }:=\Pi _{%
\limfunc{span}(X)^{\bot }}\mathcal{W}\Pi _{\limfunc{span}(X)^{\bot }}$ is
so, can immediately be subsumed under the present discussion, if one
observes that $\hat{\omega}_{\mathcal{W}}$ coincides with $\hat{\omega}_{%
\mathcal{W}^{\ast }}$ and $\mathcal{W}^{\ast }$ is constant.} Further below
we shall then address the case where $\mathcal{W}$ is allowed to depend on $%
y $. Note that the assumptions on $\mathcal{W}$ in the subsequent corollary
are certainly met if $\mathcal{W}$ is constant, symmetric, and positive
definite, and hence are satisfied in the leading case mentioned before
(provided $M$ is deterministic).

\begin{corollary}
\label{cor:polylrv1} Let $\mathfrak{F}\subseteq \mathfrak{F}_{\mathrm{all}}$
satisfy $\mathfrak{F}\supseteq \mathfrak{F}_{\mathrm{AR(}2\mathrm{)}}$ and
suppose that Assumption \ref{as:poly} holds. Suppose further that $\mathcal{W%
}$ is constant and symmetric, and that $\Pi _{\limfunc{span}(X)^{\bot }}%
\mathcal{W}\Pi _{\limfunc{span}(X)^{\bot }}$ is nonzero and nonnegative
definite. Then $\check{\beta}=\hat{\beta}$ and $\check{\Omega}=\check{\Omega}%
_{\mathcal{W}}$ satisfy Assumption \ref{Ass_5} with $N=\emptyset $. Let $T$
be of the form (\ref{F-type}) with $\check{\beta}=\hat{\beta}$, $\check{%
\Omega}=\check{\Omega}_{\mathcal{W}}$, and $N=\emptyset $. Then 
\begin{equation*}
\sup_{f\in \mathfrak{F}}P_{\mu _{0},\sigma ^{2}\Sigma (f)}(T\geq C)=1
\end{equation*}%
holds for every critical value $C$, $-\infty <C<\infty $, for every $\mu
_{0}\in \mathfrak{M}_{0}$, and for every $\sigma ^{2}\in (0,\infty )$.
\end{corollary}

We next consider the case where the matrix $\Pi _{\limfunc{span}(X)^{\bot }}%
\mathcal{W}\Pi _{\limfunc{span}(X)^{\bot }}$ is nonzero, but not
(necessarily) nonnegative definite, and thus the previous corollary is not
applicable. The subsequent corollary covers this case and is obtained under
the slightly stronger assumption that $\mathfrak{F}\supseteq \mathfrak{F}%
_{AR(2)}^{ext}$. [Note also that the case where $\mathcal{W}$ is constant
but $\Pi _{\limfunc{span}(X)^{\bot }}\mathcal{W}\Pi _{\limfunc{span}%
(X)^{\bot }}$ is equal to zero is of no interest as it leads to a
long-run-variance estimator that vanishes identically.]

\begin{corollary}
\label{cor:polylrv2} Let $\mathfrak{F}\subseteq \mathfrak{F}_{\mathrm{all}}$
satisfy $\mathfrak{F}\supseteq \mathfrak{F}_{\mathrm{AR(}2\mathrm{)}}^{ext}$
and suppose that Assumption \ref{as:poly} holds. Suppose further that $%
\mathcal{W}$ is constant and symmetric, and that $\Pi _{\limfunc{span}%
(X)^{\bot }}\mathcal{W}\Pi _{\limfunc{span}(X)^{\bot }}$ is nonzero. Then $%
\check{\beta}=\hat{\beta}$ and $\check{\Omega}=\check{\Omega}_{\mathcal{W}}$
satisfy Assumption \ref{Ass_5} with $N=\emptyset $. Let $T$ be of the form (%
\ref{F-type}) with $\check{\beta}=\hat{\beta}$, $\check{\Omega}=\check{\Omega%
}_{\mathcal{W}}$, and $N=\emptyset $. Then 
\begin{equation}
P_{0,I_{n}}(\hat{\omega}_{\mathcal{W}}\geq 0)\leq \sup_{f\in \mathfrak{F}%
}P_{\mu _{0},\sigma ^{2}\Sigma (f)}(T\geq C)  \label{eqn:lower1}
\end{equation}%
holds for every critical value $C$, $-\infty <C<\infty $, for every $\mu
_{0}\in \mathfrak{M}_{0}$, and for every $\sigma ^{2}\in (0,\infty )$.
Furthermore, for every $0\leq C<\infty $ the lower bound in the previous
display is an upper bound for the maximal power of the test under i.i.d.
errors, i.e., 
\begin{equation}
\sup_{\mu _{1}\in \mathfrak{M}_{1}}\sup_{0<\sigma ^{2}<\infty }P_{\mu
_{1},\sigma ^{2}I_{n}}(T\geq C)\leq P_{0,I_{n}}(\hat{\omega}_{\mathcal{W}%
}\geq 0).  \label{eqn:lower2}
\end{equation}
\end{corollary}

The previous corollary shows that the size of the test is bounded from below
by the probability that the long-run-variance estimator $\hat{\omega}_{_{%
\mathcal{W}}}$ used in the construction of the test statistic is
nonnegative, where the probability is taken under $N(0,I_{n})$-distributed
errors. For consistent long-run-variance estimators this probability
approaches $1$ as sample size increases, and hence the size of tests based
on such estimators $\hat{\omega}_{_{\mathcal{W}}}$ will exceed any
prescribed nominal significance level $\alpha \in (0,1)$ eventually.
Additionally, it is shown in that corollary that for nonnegative critical
values (the standard in applications) the probability $P_{0,I_{n}}(\hat{%
\omega}_{\mathcal{W}}\geq 0)$ also provides an upper bound on the \textit{%
maximal power} of the test under i.i.d. errors. Thus, if the lower bound in (%
\ref{eqn:lower1}) is small, and hence (\ref{eqn:lower1}) does not tell us
much about size, the inequality in (\ref{eqn:lower2}) shows that power must
then be small over a substantial subset of the parameter space (unless
perhaps one chooses a negative critical value). To get an idea of the
magnitude of the lower (upper) bound in (\ref{eqn:lower1}) ((\ref{eqn:lower2}%
)) in a special case, we computed $P_{0,I_{n}}(\hat{\omega}_{\mathcal{W}%
}\geq 0)$ numerically for the rectangular kernel, i.e., for $\mathcal{W}%
_{ij}=\mathbf{1}_{(-1,1)}((i-j)/M)$, for the cases when Assumption \ref%
{as:poly} is satisfied with $k_{F}=k\in \{2,3,4,\ldots ,10\}$, respectively,
sample size $n=150$, and bandwidth parameter $M=bn$ for $b\in
\{0.001,0.002,\ldots ,1\}$.\footnote{%
For $b\in \{0.994,\ldots ,1\}$ the matrix $\mathcal{W}$ has all entries
equal to one, implying that $\hat{\omega}_{_{\mathcal{W}}}$ and thus $\check{%
\Omega}_{\mathcal{W}}$ are identically zero. This is an uninteresting case
and falls outside the scope of Corollary \ref{cor:polylrv2}. [If one insists
on using the corresponding test statistic $T$ as defined in (\ref{F-type}), $%
T$ is then identically zero, leading to a useless testing procedure.] Of
course, for such values of $b$ the probability $P_{0,I_{n}}(\hat{\omega}_{%
\mathcal{W}}\geq 0)$ equals one, explaining the sharp increase of the graph
in Figure 1 for $b$ close to $1$.} The results are presented in Figure 1.%
\footnote{%
The corresponding figure in the previous versions of this paper was
incorrect due to a coding error. Furthermore, to emphasize that the
functions shown in the figure are step functions, we now use a finer grid
for $b$ in the computation; and the vertical connecting lines were added to
facilitate readability.} For all values of $b$ and $k$ the probability $%
P_{0,I_{n}}(\hat{\omega}_{\mathcal{W}}\geq 0)$ is quite large, in particular
is larger than $1/4$, and thus exceeds commonly used significance levels.
Thus, as a consequence of (\ref{eqn:lower1}), one has strong size
distortions regardless of the values of $b$ and $C$ chosen if one decides to
use a test based on the rectangular kernel. Together with (\ref{eqn:lower2}%
), Figure 1 also shows that for a large range of $b$'s the power of the
corresponding tests (with nonnegative critical value $C$) can nowhere exceed 
$0.6$, no matter how strong the deviation from the null hypothesis might be.
Note also that the probability $P_{0,I_{n}}(\hat{\omega}_{\mathcal{W}}\geq
0) $ can be easily obtained numerically in any other case, as it is the
probability that a quadratic form in a standard Gaussian random vector is
nonnegative (for the actual computation we used the algorithm by \cite%
{davies1980algorithm}).

\begin{figure}[tbp]
\centering
\includegraphics[width=\linewidth]{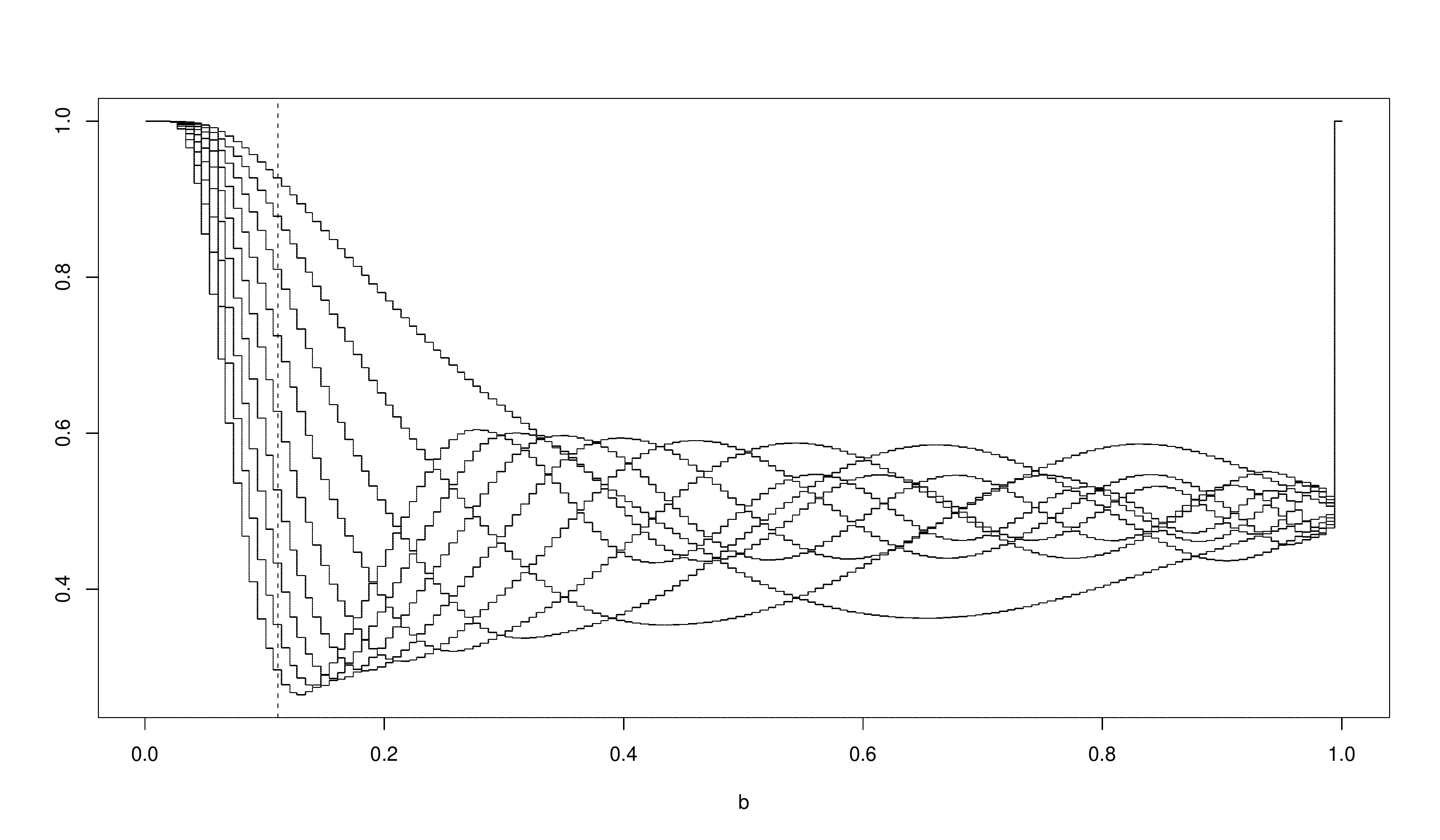}
\caption{Numerical values of $P_{0,I_{n}}(\hat{\protect\omega}_{\mathcal{W}%
}\geq 0)$ for $\mathcal{W}_{ij}=\mathbf{1}_{(-1,1)}((i-j)/(bn))$ as a
function of $b$. Sample size $n=150$ and Assumption \protect\ref{as:poly}
holds with $k_{F}=k$ and for different values of $k\in \{2,3,4,\ldots ,10\}$%
. The probabilities for $k=2$ correspond to the function with the largest
value at the dashed vertical line, the probabilities for $k=3$ correspond to
the function with the second largest value at the dashed vertical line, etc.}
\end{figure}

%

The assumption of $\mathcal{W}$ being data-independent, i.e., constant as a
function of $y\in \mathbb{R}^{n}$, in the previous two corollaries is not
satisfied for the important class of long-run-variance estimators that
incorporate prewhitening or data-dependent bandwidth parameters (e.g., \cite%
{A91}, \cite{A92}, and \cite{NW94}). An additional complication for such
estimators is that the corresponding weights matrix $\mathcal{W}(y)$, and
thus also $\check{\Omega}_{\mathcal{W}}$, are in general not well-defined
for every $y\in \mathbb{R}^{n}$. Nevertheless, after a careful structural
analysis of such estimators (similar to the results obtained in Section 3.3
of \cite{Prein2014b}), one can typically show that the resulting test
statistic satisfies the assumptions of Theorem \ref{size_bound_polynomial}
above and thus one can obtain suitable versions of the above corollaries
tailored towards test statistics based on specific classes of prewhitened
long-run-variance estimators with data-dependent bandwidth parameters. To
make this more compelling, we provide in the following such a result for a
widely used procedure in that class. We consider a version of the
AR(1)-prewhitened long-run-variance estimator based on auxiliary AR(1)
models for bandwidth selection and the Quadratic-Spectral kernel as
discussed in \cite{A92}. This is a long-run-variance estimator as in (\ref%
{eqn:lrve}), where the weights matrix is obtained as follows (the set where
all involved quantities are well-defined is given in (\ref{eq:defset})
further below): Let 
\begin{equation}
\hat{\rho}(y)=\frac{\sum_{i=2}^{n}\hat{u}_{i}(y)\hat{u}_{i-1}(y)}{%
\sum_{i=1}^{n-1}\hat{u}_{i}^{2}(y)},  \label{eqn:rhoest}
\end{equation}%
and define $\hat{v}_{i}(y)=\hat{u}_{i+1}(y)-\hat{\rho}(y)\hat{u}_{i}(y)$ for 
$i=1,\ldots ,n-1$, which one can write in an obvious way as $\hat{v}(y)=%
\mathsf{A}(\hat{\rho}(y))\hat{u}(y)$ with $\rho \mapsto \mathsf{A}(\rho )\in 
\mathbb{R}^{(n-1)\times n}$ a continuous function on $\mathbb{R}$. Define
the data-dependent bandwidth parameter $M_{\mathrm{AM}}$ via 
\begin{equation*}
M_{\mathrm{AM}}(y)=1.3221\left( n\frac{4\tilde{\rho}^{2}(y)}{(1-\tilde{\rho}%
(y))^{4}}\right) ^{1/5}\quad \text{ with }\quad \tilde{\rho}(y)=\frac{%
\sum_{i=2}^{n-1}\hat{v}_{i}(y)\hat{v}_{i-1}(y)}{\sum_{i=1}^{n-2}\hat{v}%
_{i}^{2}(y)}.
\end{equation*}%
The long-run-variance estimator $\hat{\omega}_{\mathcal{W}_{\mathrm{AM}}}$
is now obtained (granted the involved expressions are well-defined) by
choosing $\mathcal{W}$ in (\ref{eqn:lrve}) equal to 
\begin{equation*}
\mathcal{W}_{\mathrm{AM}}(y)=(1-\hat{\rho}(y))^{-2}\mathsf{A}^{\prime }(\hat{%
\rho}(y))\left[ \kappa _{\mathrm{QS}}(|i-j|/M_{\mathrm{AM}}(y))\right]
_{i,j=1}^{n-1}\mathsf{A}(\hat{\rho}(y)),
\end{equation*}%
where $[\kappa _{\mathrm{QS}}(|i-j|/M_{\mathrm{AM}}(y))]_{i,j=1}^{n-1}$ is
defined as $I_{n-1}$ in case $M_{\mathrm{AM}}(y)=0$ holds (cf., e.g., p.~821
in \cite{A91} for a definition of the Quadratic-Spectral kernel $\kappa _{%
\mathrm{QS}}$). The corresponding covariance matrix estimator $\check{\Omega}%
_{\mathcal{W}_{\mathrm{AM}}}$ is then given by plugging $\hat{\omega}_{%
\mathcal{W}_{\mathrm{AM}}}$ into (\ref{eqn:covest}). The set where $\mathcal{%
W}_{\mathrm{AM}}$ (and hence $\check{\Omega}_{\mathcal{W}_{\mathrm{AM}}}$)
is well-defined is easily seen to coincide with the set of all $y\in \mathbb{%
R}^{n}$ such that $\hat{\rho}(y)$ and $\tilde{\rho}(y)$ are both
well-defined and are not equal to $1$, i.e., with the set%
\begin{equation}
\left\{ y\in \mathbb{R}^{n}:\sum_{i=1}^{n-1}\hat{u}_{i}(y)(\hat{u}_{i+1}(y)-%
\hat{u}_{i}(y))\neq 0,~\sum_{i=1}^{n-2}\hat{v}_{i}(y)(\hat{v}_{i+1}(y)-\hat{v%
}_{i}(y))\neq 0\right\} .  \label{eq:defset}
\end{equation}%
Define $N_{\mathrm{AM}}$ as the complement of the set (\ref{eq:defset}) in $%
\mathbb{R}^{n}$. A result concerning size properties of polynomial trend
tests based on the long-run-variance estimator $\hat{\omega}_{\mathcal{W}_{%
\mathrm{AM}}}$ is now obtained by combining Theorem \ref%
{size_bound_polynomial} above with results obtained in Lemma \ref{lem:AM} in
Appendix \ref{App_D}, showing, in particular, that $\hat{\beta}$ and $\check{%
\Omega}_{\mathcal{W}_{\mathrm{AM}}}$ satisfy Assumptions \ref{Ass_5} with $%
N=N_{\mathrm{AM}}$, provided $N_{\mathrm{AM}}\neq \mathbb{R}^{n}$ holds.
Note that (i) the condition $N_{\mathrm{AM}}\neq \mathbb{R}^{n}$ only
depends on properties of the design matrix $X$ and hence can be checked, and
that (ii) in case $N_{\mathrm{AM}}=\mathbb{R}^{n}$, the matrix $\check{\Omega%
}_{\mathcal{W}_{\mathrm{AM}}}$ is \emph{nowhere} well-defined, and tests
based on this estimator hence break down in a trivial way.

\begin{corollary}
\label{cor:AM} Let $\mathfrak{F}\subseteq \mathfrak{F}_{\mathrm{all}}$
satisfy $\mathfrak{F}\supseteq \mathfrak{F}_{\mathrm{AR(}2\mathrm{)}}^{ext}$
and suppose Assumption \ref{as:poly} holds. Suppose further that $N_{\mathrm{%
AM}}\neq \mathbb{R}^{n}$. Then $\check{\beta}=\hat{\beta}$ and $\check{\Omega%
}=\check{\Omega}_{\mathcal{W}_{\mathrm{AM}}}$ satisfy Assumption \ref{Ass_5}
with $N=N_{\mathrm{AM}}$. Let $T$ be of the form (\ref{F-type}) with $\check{%
\beta}=\hat{\beta}$, $\check{\Omega}=\check{\Omega}_{\mathcal{W}_{\mathrm{AM}%
}}$, and $N=N_{\mathrm{AM}}$. Then 
\begin{equation*}
\sup_{f\in \mathfrak{F}}P_{\mu _{0},\sigma ^{2}\Sigma (f)}(T\geq C)=1
\end{equation*}%
holds for every critical value $C$, $-\infty <C<\infty $, for every $\mu
_{0}\in \mathfrak{M}_{0}$, and for every $\sigma ^{2}\in (0,\infty )$.
\end{corollary}

\begin{remark}
In the special case where Assumption \ref{as:poly} is satisfied with $%
R_{\cdot 1}\neq 0$, appropriate versions of Corollaries \ref{cor:polylrv1}, %
\ref{cor:polylrv2}, and \ref{cor:AM} maintaining only $\mathfrak{F}\supseteq 
\mathfrak{F}_{\mathrm{AR(}1\mathrm{)}}$ can be obtained by perusing Remark %
\ref{rem:AR(1)}. We abstain from spelling out details. A similar remark
applies to Corollaries \ref{cor:vogel}, \ref{cor:vogelbunzel}, and \ref%
{cor:vogelbunzelfeasible} given in the next subsection.
\end{remark}

\subsubsection{Properties of some recently suggested tests for hypotheses on
polynomial trends\label{subss:rec}}

In this subsection we discuss finite sample properties of classes of tests
for polynomial trends that have been suggested in \cite{vogelsang1998} and 
\cite{bunzelvogelsang2005}. We start with a discussion of the tests
introduced in the former article. \cite{vogelsang1998} introduces two
classes of tests for testing hypotheses on trends, in particular polynomial
trends. From Section 3.2 of \cite{vogelsang1998} it is not difficult to see
that these classes of test statistics (i.e., the classes referred to as $%
PS_{T}^{i}$ and $PSW_{T}^{i}$ in that reference) are (possibly up to a
constant positive multiplicative factor that can be absorbed into the
critical value) of the form (\ref{F-type}). More specifically, the test
statistics in \cite{vogelsang1998} are based on a combination of one of the
two estimators 
\begin{equation}
\check{\beta}_{V}(y)=\hat{\beta}_{VX}(Vy)=(X^{\prime }V^{\prime
}VX)^{-1}X^{\prime }V^{\prime }Vy\quad \text{ for }V\in \{A,I_{n}\},
\label{eqn:vobeta}
\end{equation}%
with a corresponding covariance estimator of the form 
\begin{equation}
\check{\Omega}_{c,U,i,V}^{\mathrm{Vo}}(y)=n^{j(V)}s_{A,X}^{2}(y)\exp
(cJ_{n,U}^{i}(y))R(X^{\prime }V^{\prime }VX)^{-1}R^{\prime },
\label{eqn:voomega}
\end{equation}%
for $i\in \{1,2\}$ and where $j(V)=1$ if $V=A$ and $j(V)=-1$ if $V=I_{n}$.
Here $A$ is the $n\times n$-dimensional matrix that has $0$ above the main
diagonal and $1$ on and below the main diagonal, $c\mathbb{\ }$is a real
number\footnote{%
We here also allow for the value $c=0$ in the formulation of the covariance
estimators because this turns out to be convenient in the proofs.}, $U$ is
an $n\times m$-dimensional matrix (with $m\geq 1$) such that $(X,U)$ is of
full column-rank $k+m<n$. [In \cite{vogelsang1998} the column vectors of $U$
correspond to polynomial trends of an order exceeding the polynomial trends
already contained in $\limfunc{span}(X)$.] Furthermore, 
\begin{equation}
J_{n,U}^{1}(y)=n^{-1}\hat{\beta}_{(X,U)}^{\prime }(y)G^{\prime }\left(
s_{I_{n},(X,U)}^{2}(y)G((X,U)^{\prime }(X,U))^{-1}G^{\prime }\right) ^{-1}G%
\hat{\beta}_{(X,U)}(y),  \label{eqn:J1}
\end{equation}%
and%
\begin{equation*}
J_{n,U}^{2}(y)=n^{-1}\hat{\beta}_{A(X,U)}^{\prime }(Ay)G^{\prime }\left(
s_{A,(X,U)}^{2}(y)G((X,U)^{\prime }A^{\prime }A(X,U))^{-1}G^{\prime }\right)
^{-1}G\hat{\beta}_{A(X,U)}(Ay),
\end{equation*}%
with $G=(0,I_{m})\in \mathbb{R}^{m\times (k+m)}$, where we use the notation%
\begin{equation*}
s_{D_{1},D_{2}}^{2}(y)=n^{-1}y^{\prime }D_{1}^{\prime }\Pi _{\limfunc{span}%
(D_{1}D_{2})^{\bot }}D_{1}y
\end{equation*}%
for nonsingular $D_{1}\in \mathbb{R}^{n\times n}$ and for $D_{2}\in \mathbb{R%
}^{n\times l}$ of rank $l\leq n$. It is obvious from the above expressions
that the covariance estimator $\check{\Omega}_{c,U,i,V}^{\mathrm{Vo}}$ is
not well-defined on all of $\mathbb{R}^{n}$. However, it is also not
difficult to see that the set where such an estimator is well-defined
coincides with $\mathbb{R}^{n}\backslash \limfunc{span}(X,U)$, see the proof
of Lemma \ref{lem:vogel} in Appendix \ref{App_D}. We stress once more that
the matrix $U$ used in the construction above is chosen in a particular way
in \cite{vogelsang1998}. We do not impose such a restriction here, because
it would unnecessarily complicate the presentation of the result below, and
because this restrictions is actually not necessary for establishing the
result. The following result now shows, in particular, that the tests
suggested in \cite{vogelsang1998} suffer from substantial size distortions
in case $\mathfrak{F}\supseteq \mathfrak{F}_{AR(2)}^{ext}$.

\begin{corollary}
\label{cor:vogel} Let $\mathfrak{F}\subseteq \mathfrak{F}_{\mathrm{all}}$
satisfy $\mathfrak{F}\supseteq \mathfrak{F}_{\mathrm{AR(}2\mathrm{)}}^{ext}$
and suppose Assumption \ref{as:poly} holds. Let $V\in \{A,I_{n}\}$, $c\in 
\mathbb{R}$, $i\in \{1,2\}$, and let $U$ be an $n\times m$-dimensional
matrix with $m\geq 1$, $k+m<n$, such that $(X,U)$ is of full column-rank.
Then $\check{\beta}=\check{\beta}_{V}$ and $\check{\Omega}=\check{\Omega}%
_{c,U,i,V}^{\mathrm{Vo}}$ satisfy Assumption \ref{Ass_5} with $N=\limfunc{%
span}(X,U)$. Let $T$ be of the form (\ref{F-type}) with $\check{\beta}=%
\check{\beta}_{V}$, $\check{\Omega}=\check{\Omega}_{c,U,i,V}^{\mathrm{Vo}}$,
and $N=\limfunc{span}(X,U)$. Then 
\begin{equation*}
\sup_{f\in \mathfrak{F}}P_{\mu _{0},\sigma ^{2}\Sigma (f)}(T\geq C)=1
\end{equation*}%
holds for every critical value $C$, $-\infty <C<\infty $, for every $\mu
_{0}\in \mathfrak{M}_{0}$, and for every $\sigma ^{2}\in (0,\infty )$.
\end{corollary}

Next we turn to the tests introduced in \cite{bunzelvogelsang2005}. We first
discuss tests introduced in that article with data-independent tuning
parameters and data-independent critical values: These tests are based on
the OLS estimator $\hat{\beta}$ and two classes of covariance matrix
estimators, both of which incorporate a tuning parameter $c\in \mathbb{R}$,
and which are defined as 
\begin{equation}
\check{\Omega}_{\mathcal{W},U,c}^{\mathrm{BV},J}(y)=\hat{\omega}_{\mathcal{W}%
}(y)\exp (cJ_{n,U}^{1}(y))R(X^{\prime }X)^{-1}R^{\prime },
\label{eqn:vogelbunzelcov1}
\end{equation}%
where $U$ is an $n\times m$-dimensional matrix with $m\geq 1$ such that $%
(X,U)$ is of full column-rank $k+m<n$ (note that $\hat{\omega}_{\mathcal{W}}$
and $J_{n,U}^{1}$ have been defined in (\ref{eqn:lrve}) and (\ref{eqn:J1})
above), and 
\begin{equation}
\check{\Omega}_{\mathcal{W},c}^{\mathrm{BV}}(y)=\hat{\omega}_{\mathcal{W}%
}(y)\exp \left( cn^{-2}\frac{\hat{u}^{\prime }(y)A^{\prime }A\hat{u}(y)}{%
\hat{u}^{\prime }(y)\hat{u}(y)}\right) R(X^{\prime }X)^{-1}R^{\prime }
\label{eqn:vogelbunzelcov2}
\end{equation}%
where $A$ has been defined below (\ref{eqn:voomega}).\ The subsequent result
applies, in particular, if $\mathcal{W}_{ij}=\kappa (|i-j|/M)$ where $M>0$
is a (fixed) real number and $\kappa $ is a kernel function such that $%
\mathcal{W}$ is positive definite, including the recommendation in \cite%
{bunzelvogelsang2005} to use the Daniell kernel. In that case, and more
generally whenever $\Pi _{\limfunc{span}(X)^{\bot }}\mathcal{W}\Pi _{%
\limfunc{span}(X)^{\bot }}$ is nonzero and nonnegative definite (with $%
\mathcal{W}$ constant\footnote{%
Cf. Footnote \ref{ftnt_1}} and symmetric), the subsequent corollary shows
that the above mentioned tests in \cite{bunzelvogelsang2005} have size equal
to one if $\mathfrak{F}\supseteq \mathfrak{F}_{AR(2)}^{ext}$; in case $\Pi _{%
\limfunc{span}(X)^{\bot }}\mathcal{W}\Pi _{\limfunc{span}(X)^{\bot }}$ is
nonzero but not nonnegative definite, a lower bound on the size is obtained,
which also provides an upper bound for the power in the case of i.i.d.
errors. A discussion similar to the discussion following Corollary \ref%
{cor:polylrv2} also applies here (cf. also Figure 1).

\begin{corollary}
\label{cor:vogelbunzel} Let $\mathfrak{F}\subseteq \mathfrak{F}_{\mathrm{all}%
}$ satisfy $\mathfrak{F}\supseteq \mathfrak{F}_{\mathrm{AR(}2\mathrm{)}%
}^{ext}$ and suppose Assumption \ref{as:poly} holds. Suppose that $\mathcal{W%
}$ is constant and symmetric, that $\Pi _{\limfunc{span}(X)^{\bot }}\mathcal{%
W}\Pi _{\limfunc{span}(X)^{\bot }}$ is nonzero, and that $c\in \mathbb{R}$.
Furthermore, for the statements that involve $U$, suppose $U$ is an $n\times
m$-dimensional matrix with $m\geq 1$ such that $(X,U)$ is of full
column-rank $k+m<n$. Then, $\check{\beta}=\hat{\beta}$ and $\check{\Omega}=%
\check{\Omega}_{\mathcal{W},c}^{\mathrm{BV}}$ ($\check{\beta}=\hat{\beta}$
and $\check{\Omega}=\check{\Omega}_{\mathcal{W},U,c}^{\mathrm{BV},J}$,
respectively) satisfy Assumption \ref{Ass_5} with $N=\limfunc{span}(X)$ ($N=%
\limfunc{span}(X,U)$, respectively). Let $T$ be of the form (\ref{F-type})
with $\check{\beta}=\hat{\beta}$, $\check{\Omega}=\check{\Omega}_{\mathcal{W}%
,c}^{\mathrm{BV}}$, and $N=\limfunc{span}(X)$, or with $\check{\beta}=\hat{%
\beta}$, $\check{\Omega}=\check{\Omega}_{\mathcal{W},U,c}^{\mathrm{BV},J}$,
and $N=\limfunc{span}(X,U)$. Then 
\begin{equation*}
P_{0,I_{n}}(\hat{\omega}_{\mathcal{W}}\geq 0)\leq \sup_{f\in \mathfrak{F}%
}P_{\mu _{0},\sigma ^{2}\Sigma (f)}(T\geq C)
\end{equation*}%
holds for every critical value $C$, $-\infty <C<\infty $, for every $\mu
_{0}\in \mathfrak{M}_{0}$, and for every $\sigma ^{2}\in (0,\infty )$. The
lower bound equals $1$ in case $\Pi _{\limfunc{span}(X)^{\bot }}\mathcal{W}%
\Pi _{\limfunc{span}(X)^{\bot }}$ is nonnegative definite. Furthermore, for
every $0\leq C<\infty $ the lower bound in the previous display is an upper
bound for the maximal power of the test under i.i.d. errors, i.e., 
\begin{equation}
\sup_{\mu _{1}\in \mathfrak{M}_{1}}\sup_{0<\sigma ^{2}<\infty }P_{\mu
_{1},\sigma ^{2}I_{n}}(T\geq C)\leq P_{0,I_{n}}(\hat{\omega}_{\mathcal{W}%
}\geq 0).  \label{eqn:lower4}
\end{equation}
\end{corollary}

We shall now turn to the approach \cite{bunzelvogelsang2005} suggest for
practical applications. This approach is based on a data-driven selection of
the weights matrix $\mathcal{W}$ and of the tuning parameter $c$, and on a
data-driven selection of the critical value $C$. Their approach is as
follows: \cite{bunzelvogelsang2005} focus on $\hat{\omega}_{\mathcal{W}}$
based on the Daniell kernel. More specifically, they set $\mathcal{W}%
_{ij}=\kappa _{D}(|i-j|/\max (bn,2))$ (cf. \cite{bunzelvogelsang2005},
Appendix B, for a definition of the Daniell kernel). Recall that, regardless
of the value of $b$, the matrix with elements $\mathcal{W}_{ij}=\kappa
_{D}(|i-j|/\max (bn,2))$ based on the Daniell kernel is positive definite.
The authors recommend to choose $b$ as a positive piecewise constant
function of $\hat{\rho}$ (which has been defined in (\ref{eqn:rhoest})
above), more precisely, for constants $a_{i}\in (0,\infty )$, $i=0,\ldots
,m^{\prime }$ ($m^{\prime }\in \mathbb{N}$), and $\bar{a}_{i}\in \mathbb{R}$%
, $i=1,\ldots ,m^{\prime }$, they suggest to use 
\begin{equation*}
b_{\mathrm{BV}}(y,a,\bar{a})=a_{0}+\sum_{i=1}^{{m^{\prime }}}a_{i}\mathbf{1}%
_{[\bar{a}_{i},\infty )}(\hat{\rho}(y)).
\end{equation*}%
For a recommendation concerning the choice of these constants see \cite%
{bunzelvogelsang2005}, p.~388. Furthermore, \cite{bunzelvogelsang2005}
suggest to choose their data-driven critical value $C$ and a data-driven
tuning parameter $c$ as a polynomial function of $b_{\mathrm{BV}}(y,a,\bar{a}%
)$, respectively. More precisely, for constants $h_{0},\ldots ,h_{m^{\prime
\prime }}\in \mathbb{R}$ ($m^{\prime \prime }\in \mathbb{N}$, $h_{m^{\prime
\prime }}\neq 0$) and $p_{0},\ldots ,p_{m^{^{\prime \prime \prime }}}\in 
\mathbb{R}$ ($m^{\prime \prime \prime }\in \mathbb{N}$, $p_{m^{\prime \prime
\prime }}\neq 0$) they suggest to use 
\begin{equation*}
C_{\mathrm{BV}}(y,h)=\sum_{i=0}^{m^{\prime \prime }}h_{i}(b_{\mathrm{BV}%
}(y,a,\bar{a}))^{i}\quad \text{ and }\quad c_{\mathrm{BV}}(y,p)=%
\sum_{i=0}^{m^{\prime \prime \prime }}p_{i}(b_{\mathrm{BV}}(y,a,\bar{a}%
))^{i}.
\end{equation*}%
Then they set 
\begin{equation*}
\mathcal{W}_{\mathrm{BV}}(y)=\left[ \kappa _{D}(|i-j|/\max (b_{\mathrm{BV}%
}(y,a,\bar{a})n,2))\right] _{i,j=1}^{n},
\end{equation*}%
and define, in correspondence with (\ref{eqn:vogelbunzelcov1}) and (\ref%
{eqn:vogelbunzelcov2}), the covariance estimators%
\begin{equation*}
\check{\Omega}_{U,a,\bar{a},h,p}^{\mathrm{BV},J}(y)=\hat{\omega}_{\mathcal{W}%
_{\mathrm{BV}}}(y)\exp \left( c_{\mathrm{BV}}(y,p)J_{n,U}^{1}(y)\right)
R(X^{\prime }X)^{-1}R^{\prime }
\end{equation*}%
and 
\begin{equation*}
\check{\Omega}_{a,\bar{a},h,p}^{\mathrm{BV}}(y)=\hat{\omega}_{\mathcal{W}_{%
\mathrm{BV}}}(y)\exp \left( c_{\mathrm{BV}}(y,p)n^{-2}\frac{\hat{u}^{\prime
}(y)A^{\prime }A\hat{u}(y)}{\hat{u}^{\prime }(y)\hat{u}(y)}\right)
R(X^{\prime }X)^{-1}R^{\prime }.
\end{equation*}%
The vectors of (constant) tuning parameters $a=(a_{0},\ldots ,a_{m^{\prime
}})^{\prime }$, $\bar{a}=(\bar{a}_{1},\ldots ,\bar{a}_{m^{\prime }})^{\prime
}$, $h=(h_{0},\ldots ,h_{m^{\prime \prime }})^{\prime }$, and $%
p=(p_{0},\ldots ,p_{m^{\prime \prime \prime }})^{\prime }$ this approach is
based on are tabulated in \cite{bunzelvogelsang2005} for certain cases, and
need to be obtained numerically, following the rationale in \cite%
{bunzelvogelsang2005}, for the cases not tabulated in that paper.
Furthermore, the data-driven tuning parameters $b_{\mathrm{BV}}$ and $c_{%
\mathrm{BV}}$ as well as the data-driven critical value $C_{\mathrm{BV}}$
are well-defined for a given $y\in \mathbb{R}^{n}$ if and only if $\hat{\rho}%
(y)$ is well-defined, i.e., these quantities are well-defined on the
complement of the closed set%
\begin{equation}
\tilde{N}:=\left\{ y\in \mathbb{R}^{n}:\sum_{i=1}^{n-1}\hat{u}%
_{i}^{2}(y)=0\right\} .  \label{eqn:excbv}
\end{equation}%
Clearly, $\limfunc{span}(X)$ is contained in $\tilde{N}$. Hence, it is not
difficult to see that the estimator $\check{\Omega}_{a,\bar{a},h,p}^{\mathrm{%
BV}}$ is well-defined on $\mathbb{R}^{n}\backslash \tilde{N}$ and that the
estimator $\check{\Omega}_{U,a,\bar{a},h,p}^{\mathrm{BV},J}$ is well-defined
on $\mathbb{R}^{n}\backslash (\limfunc{span}(X,U)\cup \tilde{N})$. In fact,
under Assumption \ref{as:poly} we have that $\tilde{N}=\limfunc{span}(X)$
(see the proof of the subsequent corollary). Consequently, under Assumption %
\ref{as:poly}, the estimator $\check{\Omega}_{a,\bar{a},h,p}^{\mathrm{BV}}$
is well defined on $\mathbb{R}^{n}\backslash \limfunc{span}(X)$ and $\check{%
\Omega}_{U,a,\bar{a},h,p}^{\mathrm{BV},J}$ is well-defined on $\mathbb{R}%
^{n}\backslash \limfunc{span}(X,U)$. [In order that the data-driven critical
value is also defined for every $y$, we set $C_{\mathrm{BV}}(y,h)$ equal to
an arbitrary value ($0$, say) on the null-set $\tilde{N}$. Of course, the
choice of assignment on this null-set is inconsequential for the result
below.]

The following corollary shows that the tests for hypotheses concerning
polynomial trends based on data-driven tuning parameters and a data-driven
critical value as suggested in \cite{bunzelvogelsang2005} have size one in
case $\mathfrak{F}\supseteq \mathfrak{F}_{AR(2)}^{ext}$. The proof of this
is based on a similar approach as used in the proof of Corollary \ref%
{cor:vogelbunzel} above, but has to deal with the fact that the choice of
the tuning parameters and the critical value is data-driven, and hence is
more involved. In particular, it turns out that in order for Assumption \ref%
{Ass_5} to be satisfied for the covariance estimators used here, one has to
work with null-sets $N_{\mathrm{BV},U}$ and $N_{\mathrm{BV}}$ that are
larger than $\limfunc{span}(X,U)$ and $\limfunc{span}(X)$, respectively.

\begin{corollary}
\label{cor:vogelbunzelfeasible} Let $\mathfrak{F}\subseteq \mathfrak{F}_{%
\mathrm{all}}$ satisfy $\mathfrak{F}\supseteq \mathfrak{F}_{\mathrm{AR(}2%
\mathrm{)}}^{ext}$ and suppose Assumption \ref{as:poly} holds. Let $a_{i}\in
(0,\infty )$ for $i=0,\ldots ,m^{\prime }$ ($m^{\prime }\in \mathbb{N}$), $%
\bar{a}_{i}\in \mathbb{R}$ for $i=1,\ldots ,m^{\prime }$, $h_{i}\in \mathbb{R%
}$ for $i=0,\ldots ,m^{\prime \prime }$ with $h_{m^{\prime \prime }}\neq 0$
and $m^{\prime \prime }\in \mathbb{N}$, and $p_{i}\in \mathbb{R}$ for $%
i=0,\ldots ,m^{\prime \prime \prime }$ with $p_{m^{\prime \prime \prime
}}\neq 0$ and $m^{\prime \prime \prime }\in \mathbb{N}$. Furthermore, for
the statements that involve $U$, suppose $U$ is an $n\times m$-dimensional
matrix with $m\geq 1$ such that $(X,U)$ is of full column-rank $k+m<n$.
Then, $\check{\beta}=\hat{\beta}$ and $\check{\Omega}=\check{\Omega}_{a,\bar{%
a},h,p}^{\mathrm{BV}}$ satisfy Assumption \ref{Ass_5} with $N=N_{\mathrm{BV}%
} $ (defined in Lemma \ref{lem:BVfeas} in Appendix \ref{App_D}), and $\check{%
\beta}=\hat{\beta}$ and $\check{\Omega}=\check{\Omega}_{U,a,\bar{a},h,p}^{%
\mathrm{BV},J}$ satisfy Assumption \ref{Ass_5} with $N=N_{\mathrm{BV},U}$
(defined in Lemma \ref{lem:BVfeas}). Let $T$ be of the form (\ref{F-type})
with $\check{\beta}=\hat{\beta}$, $\check{\Omega}=\check{\Omega}_{a,\bar{a}%
,h,p}^{\mathrm{BV}}$, and $N=N_{\mathrm{BV}}$, or with $\check{\beta}=\hat{%
\beta}$, $\check{\Omega}=\check{\Omega}_{U,a,\bar{a},h,p}^{\mathrm{BV},J}$,
and $N=N_{\mathrm{BV},U}$. Then 
\begin{equation*}
\sup_{f\in \mathfrak{F}}P_{\mu _{0},\sigma ^{2}\Sigma (f)}(\{y\in \mathbb{R}%
^{n}:T(y)\geq C_{\mathrm{BV}}(y,h)\})=1
\end{equation*}%
holds for every $\mu _{0}\in \mathfrak{M}_{0}$ and for every $\sigma ^{2}\in
(0,\infty )$.
\end{corollary}

\begin{remark}
Alternatively one can consider $T^{\ast }$, where 
\begin{equation*}
T^{\ast }(y)=(R\hat{\beta}(y)-r)^{\prime }\left( \check{\Omega}_{a,\bar{a}%
,h,p}^{\mathrm{BV}}(y)\right) ^{-1}(R\hat{\beta}(y)-r)
\end{equation*}%
for all $y\in \mathbb{R}^{n}\backslash \limfunc{span}(X)$ such that $\check{%
\Omega}_{a,\bar{a},h,p}^{\mathrm{BV}}(y)$ is nonsingular, and where $T^{\ast
}(y)=0$ else, (and we can similarly define a test statistic $T^{\ast \ast }$
with $\check{\Omega}_{U,a,\bar{a},h,p}^{\mathrm{BV},J}$ and $\limfunc{span}%
(X,U)$ in place of $\check{\Omega}_{a,\bar{a},h,p}^{\mathrm{BV}}$ and $%
\limfunc{span}(X)$, respectively). While $T^{\ast }$ and $T^{\ast \ast }$
are well-defined test statistics, we are not guaranteed that $\hat{\beta}$
and $\check{\Omega}_{a,\bar{a},h,p}^{\mathrm{BV}}$ ($\hat{\beta}$ and $%
\check{\Omega}_{U,a,\bar{a},h,p}^{\mathrm{BV},J}$, respectively) satisfy
Assumption \ref{Ass_5} with $N=\limfunc{span}(X)$ ($N=\limfunc{span}(X,U)$,
respectively). However, $T^{\ast }$ as well as $T^{\ast \ast }$ differ from
the corresponding test statistics considered in the preceding corollary at
most on a null-set, hence the conclusions of the corollary carry over to $%
T^{\ast }$ and $T^{\ast \ast }$.
\end{remark}

\subsection{Cyclical trends\label{subs:cycl}}

We here consider briefly the case when one tests hypotheses concerning a
cyclical trend, i.e., when the following assumption is satisfied:

\begin{assumption}
\label{as:cycl} Suppose that $X=(E_{n,0}(\omega ),\tilde{X})$ for some $%
\omega \in (0,\pi )$ where $\tilde{X}$ is an $n\times (k-2)$-dimensional
matrix such that $X$ has rank $k$ (here $\tilde{X}$ is the empty matrix if $%
k=2$). Furthermore, suppose that the restriction matrix $R$ has a nonzero
column $R_{\cdot i}$ for some $i=1,2$, i.e., the hypothesis involves
coefficients of the cyclical component.
\end{assumption}

Under this assumption we obtain the subsequent theorem from Theorem \ref%
{size_one_AR2}.

\begin{theorem}
\label{size_one_cyclical}Let $\mathfrak{F}\subseteq \mathfrak{F}_{\mathrm{all%
}}$ satisfy $\mathfrak{F}\supseteq \mathfrak{F}_{\mathrm{AR(}2\mathrm{)}}$
and suppose Assumption \ref{as:cycl} holds. Let $T$ be a
nonsphericity-corrected F-type test statistic of the form (\ref{F-type})
based on $\check{\beta}$ and $\check{\Omega}$ satisfying Assumption \ref%
{Ass_5} with $N=\emptyset $. Furthermore, assume that $\check{\Omega}(y)$ is
nonnegative definite for every $y\in \mathbb{R}^{n}$. Then 
\begin{equation*}
\sup_{f\in \mathfrak{F}}P_{\mu _{0},\sigma ^{2}\Sigma (f)}(T\geq C)=1
\end{equation*}%
holds for every critical value $C$, $-\infty <C<\infty $, for every $\mu
_{0}\in \mathfrak{M}_{0}$, and for every $\sigma ^{2}\in (0,\infty )$.
\end{theorem}

Under a slightly stronger condition on $\mathfrak{F}$, the following theorem
is applicable in case the assumption that $N=\emptyset $ or the nonnegative
definiteness assumption on $\check{\Omega}$ in the previous theorem are
violated.

\begin{theorem}
\label{size_bound_cyclical}Let $\mathfrak{F}\subseteq \mathfrak{F}_{\mathrm{%
all}}$ satisfy $\mathfrak{F}\supseteq \mathfrak{F}_{\mathrm{AR(}2\mathrm{)}%
}^{ext}$. Suppose Assumption \ref{as:cycl} holds. Let $T$ be a
nonsphericity-corrected F-type test statistic of the form (\ref{F-type})
based on $\check{\beta}$ and $\check{\Omega}$ satisfying Assumption \ref%
{Ass_5}. Furthermore, assume that $\check{\Omega}$ also satisfies Assumption %
\ref{Ass_7}. Then for every critical value $C$, $-\infty <C<\infty $, for
every $\mu _{0}\in \mathfrak{M}_{0}$, and for every $\sigma ^{2}\in
(0,\infty )$ it holds that 
\begin{equation*}
P_{0,I_{n}}(\check{\Omega}\text{ is nonnegative definite})\leq K(\omega
)\leq \sup_{f\in \mathfrak{F}}P_{\mu _{0},\sigma ^{2}\Sigma (f)}\left( T\geq
C\right) ,
\end{equation*}%
where $K(\omega )$ is defined in Theorem \ref{theo:AR2+}.
\end{theorem}

Using these results, one can now obtain similar results as in Subsection \ref%
{subss:rec} concerning the tests developed in \cite{vogelsang1998} and \cite%
{bunzelvogelsang2005} under Assumption \ref{as:cycl}. Due to space
constraints, however, we do not spell out the details.

\begin{remark}
\label{rem:final}\emph{(The cases }$\omega =0$ \emph{or }$\omega =\pi $\emph{%
)} (i) In case $\omega =0$ (or $\omega =\pi $) consider Assumption \ref%
{as:cycl} with the understanding that $X=(\bar{E}_{n,0}(\omega ),\tilde{X})$%
, that $\tilde{X}$ is now $n\times (k-1)$-dimensional, and that $R_{\cdot
1}\neq 0$, where $\bar{E}_{n,0}(\omega )$ denotes the first column of $%
E_{n,0}(\omega )$. Then Theorems \ref{size_one_cyclical} and \ref%
{size_bound_cyclical} continue to hold with this interpretation of
Assumption \ref{as:cycl}. Also note that the case $\omega =0$ can be
subsumed under the results of Subsection \ref{subs:poly} by setting $k_{F}=1$%
.

(ii) In case $\omega =0$ (or $\omega =\pi $), Theorem \ref{size_one_cyclical}
(with the before mentioned interpretation of Assumption \ref{as:cycl}) in
fact continues to hold under the weaker assumption that $\mathfrak{F}%
\supseteq \mathfrak{F}_{\mathrm{AR(}1\mathrm{)}}$.\footnote{%
In fact, it holds even more generally for any covariance model $\mathfrak{C}$
that has $\limfunc{span}(\bar{E}_{n,0}(\omega ))$ as a concentration space.}
This follows from Part 3 of Corollary 5.17 in \cite{PP2016} upon noting that 
$\mathcal{Z}=\limfunc{span}(\bar{E}_{n,0}(\omega ))$ is a concentration
space of the covariance model $\mathfrak{C}(\mathfrak{F})$, that $\check{%
\Omega}$ vanishes on $\limfunc{span}(X)\supseteq \mathcal{Z}$ as a
consequence of the assumption $N=\emptyset $ (see the discussion following
(27) in \cite{PP2016}), and that $R\check{\beta}(z)\neq 0$ for every $z\in 
\mathcal{Z}$ with $z\neq 0$.\footnote{%
This is proved similarly as in Footnote \ref{nezero}.}

(iii) In case $\omega =0$ (or $\omega =\pi $), Theorem \ref%
{size_bound_cyclical} (with the before mentioned interpretation of
Assumption \ref{as:cycl}) also continues to hold under the weaker assumption
that $\mathfrak{F}\supseteq \mathfrak{F}_{\mathrm{AR(}1\mathrm{)}}$ if $\bar{%
\xi}_{\omega }(x)$ in the definition of $K(\omega )$ is now replaced by $%
\breve{\xi}_{\omega }(x)$ defined as%
\begin{eqnarray*}
\breve{\xi}_{\omega }(x) &=&(R\hat{\beta}_{X}(\bar{E}_{n,0}(\omega
)x))^{\prime }\check{\Omega}^{-1}\left( \left( \left( \bar{E}_{n,0}(\omega )%
\bar{E}_{n,0}(\omega )^{\prime }\right) ^{1/2}+D(\omega )^{1/2}\right) 
\mathbf{G}\right) R\hat{\beta}_{X}(\bar{E}_{n,0}(\omega )x) \\
&=&x^{2}R_{\cdot 1}^{\prime }\check{\Omega}^{-1}\left( \left( \left( \bar{E}%
_{n,0}(\omega )\bar{E}_{n,0}(\omega )^{\prime }\right) ^{1/2}+D(\omega
)^{1/2}\right) \mathbf{G}\right) R_{\cdot 1}
\end{eqnarray*}%
on the event where $\{((\bar{E}_{n,0}(\omega )\bar{E}_{n,0}(\omega )^{\prime
})^{1/2}+D(\omega )^{1/2})\mathbf{G}\in \mathbb{R}^{n}\backslash N^{\ast }\}$
and by $\breve{\xi}_{\omega }(x)=0$ otherwise, and if the distribution $%
P_{0,I_{n}}$ appearing in the lower bound is replaced by $P_{0,\Phi (\omega
)}$ where $\Phi (\omega )=\bar{E}_{n,0}(\omega )\bar{E}_{n,0}(\omega
)^{\prime }+D(\omega )$ is nonsingular. Note that then $K(\omega )$ reduces
to $P_{0,\Phi (\omega )}(R_{\cdot 1}^{\prime }\check{\Omega}^{-1}R_{\cdot
1}\geq 0)$. Here $D(0)$ is the matrix $D$ given in Part 3 and $D(\pi )$ is
the matrix $D$ given in Part 4 of Lemma G.1 in \cite{PP2016}. This can be
proved by making use of Theorem 5.19 and Lemma G.1 in \cite{PP2016}.
\end{remark}

\appendix

\section{Appendix: Proofs and auxiliary results for Section \protect\ref%
{Sec_general_results} \label{App A}}

\begin{lemma}
\label{associated_cov_model}Let $\mathfrak{C}$ be a covariance model and let 
$\mathcal{L}$ be a linear subspace of $\mathbb{R}^{n}$ with $\dim (\mathcal{L%
})=l<n$. Let $\mathfrak{C}^{\sharp }=\left\{ \Sigma ^{\sharp }:\Sigma \in 
\mathfrak{C}\right\} $ and $\mathfrak{C}^{\natural }=\left\{ \Sigma
^{\natural }:\Sigma \in \mathfrak{C}\right\} $, where $\Sigma ^{\sharp }=%
\mathcal{L}(\Sigma )+\lambda _{l+1}(\mathcal{L}(\Sigma ))\Pi _{\mathcal{L}}$
and where $\Sigma ^{\natural }=\mathcal{L}(\Sigma )+\Pi _{\mathcal{L}}$;
here $\lambda _{l+1}(\mathcal{L}(\Sigma ))$ denotes the $(l+1)$-th
eigenvalue of $\mathcal{L}(\Sigma )$ counting (with multiplicity) from
smallest to largest. Then $\mathfrak{C}^{\sharp }$ and $\mathfrak{C}%
^{\natural }$ are covariance models. Furthermore, the collection of
concentration spaces of $\mathfrak{C}^{\sharp }$ coincides with $\mathbb{J}(%
\mathcal{L},\mathfrak{C})$, and the collection of concentration spaces of $%
\mathfrak{C}^{\natural }$ coincides with the collection $\left\{ \mathcal{S}+%
\mathcal{L}:\mathcal{S}\in \mathbb{J}(\mathcal{L},\mathfrak{C})\right\} $.
\end{lemma}

\textbf{Proof:} 1. That $\mathfrak{C}^{\sharp }$ and $\mathfrak{C}^{\natural
}$ are covariance models is obvious since the elements of these two
collections are clearly symmetric and positive definite matrices (as $%
\lambda _{l+1}(\mathcal{L}(\Sigma ))>0$ by construction).

2. Suppose $\mathcal{S}\in \mathbb{J}(\mathcal{L},\mathfrak{C})$. Then $%
\mathcal{S}=\limfunc{span}(\bar{\Sigma})$ for some $\bar{\Sigma}\in \limfunc{%
cl}(\mathcal{L}(\mathfrak{C}))$ with $\limfunc{rank}(\bar{\Sigma})<n-l$. In
particular, $\bar{\Sigma}$ is the limit of $\mathcal{L}(\Sigma _{m})$ for a
sequence $\Sigma _{m}\in \mathfrak{C}$. But then $\Sigma _{m}^{\sharp }=%
\mathcal{L}(\Sigma _{m})+\lambda _{l+1}(\mathcal{L}(\Sigma _{m}))\Pi _{%
\mathcal{L}}$ belongs to $\mathfrak{C}^{\sharp }$ and converges to $\bar{%
\Sigma}$ for $m\rightarrow \infty $, since $\lambda _{l+1}(\mathcal{L}%
(\Sigma _{m}))$ converges to $\lambda _{l+1}(\bar{\Sigma})$, which equals
zero as a consequence of $\limfunc{rank}(\bar{\Sigma})<n-l$. This shows that 
$\limfunc{span}(\bar{\Sigma})$, and hence $\mathcal{S}$, is a concentration
space of $\mathfrak{C}^{\sharp }$. Conversely, suppose $\mathcal{Z}$ is a
concentration space of $\mathfrak{C}^{\sharp }$. Then $\mathcal{Z}=\limfunc{%
span}(\breve{\Sigma})$ for some singular matrix that is the limit of some
sequence $\Sigma _{m}^{\sharp }\in \mathfrak{C}^{\sharp }$. In particular, $%
\Sigma _{m}^{\sharp }=\mathcal{L}(\Sigma _{m})+\lambda _{l+1}(\mathcal{L}%
(\Sigma _{m}))\Pi _{\mathcal{L}}$ holds for some sequence $\Sigma _{m}\in 
\mathfrak{C}$. Since the matrices $\mathcal{L}(\Sigma _{m})$ reside in the
unit sphere in $\mathbb{R}^{n\times n}$, we have convergence of $\mathcal{L}%
(\Sigma _{m_{i}})$ to a limit $\bar{\Sigma}\in \mathbb{R}^{n\times n}$ along
an appropriate subsequence $m_{i}$; in particular, $\bar{\Sigma}\in \limfunc{%
cl}(\mathcal{L}(\mathfrak{C}))$ follows. Furthermore, we conclude that $%
\Sigma _{m_{i}}^{\sharp }$ converges to $\bar{\Sigma}+\lambda _{l+1}(\bar{%
\Sigma})\Pi _{\mathcal{L}}$, and hence obtain the equality $\breve{\Sigma}=%
\bar{\Sigma}+\lambda _{l+1}(\bar{\Sigma})\Pi _{\mathcal{L}}$. Since $\bar{%
\Sigma}$ is certainly symmetric and nonnegative definite, we have that $%
\lambda _{l+1}(\bar{\Sigma})\geq 0$. Note that $\bar{\Sigma}x=0$ for every $%
x\in \mathcal{L}$ by construction of $\bar{\Sigma}$. Hence $\limfunc{rank}(%
\bar{\Sigma})\leq n-l$ must hold. If $\limfunc{rank}(\bar{\Sigma})=n-l$
would hold we would have $\lambda _{l+1}(\bar{\Sigma})>0$, implying that $%
\bar{\Sigma}+\lambda _{l+1}(\bar{\Sigma})\Pi _{\mathcal{L}}$ is nonsingular,
contradicting singularity of $\breve{\Sigma}$. Consequently, $\limfunc{rank}(%
\bar{\Sigma})<n-l$ and $\lambda _{l+1}(\bar{\Sigma})=0$ must hold, implying
that $\mathcal{S}=\limfunc{span}(\bar{\Sigma})$ belongs to $\mathbb{J}(%
\mathcal{L},\mathfrak{C})$ and that $\breve{\Sigma}=\bar{\Sigma}$ holds. But
this shows $\mathcal{Z}=\mathcal{S}\in \mathbb{J}(\mathcal{L},\mathfrak{C})$.

3. Suppose $\mathcal{S}\in \mathbb{J}(\mathcal{L},\mathfrak{C})$. Then $%
\mathcal{S}=\limfunc{span}(\bar{\Sigma})$ for some $\bar{\Sigma}\in \limfunc{%
cl}(\mathcal{L}(\mathfrak{C}))$ with $\limfunc{rank}(\bar{\Sigma})<n-l$. In
particular, $\bar{\Sigma}$ is the limit of $\mathcal{L}(\Sigma _{m})$ for a
sequence $\Sigma _{m}\in \mathfrak{C}$. But then $\Sigma _{m}^{\natural }=%
\mathcal{L}(\Sigma _{m})+\Pi _{\mathcal{L}}$ belongs to $\mathfrak{C}%
^{\natural }$ and converges to $\bar{\Sigma}+\Pi _{\mathcal{L}}$ for $%
m\rightarrow \infty $. Now $\bar{\Sigma}+\Pi _{\mathcal{L}}$ is singular
since $\limfunc{rank}(\bar{\Sigma})<n-l$. Hence, $\limfunc{span}(\bar{\Sigma}%
+\Pi _{\mathcal{L}})$ is a concentration space of $\mathfrak{C}^{\natural }$
and $\limfunc{span}(\bar{\Sigma}+\Pi _{\mathcal{L}})=\limfunc{span}(\bar{%
\Sigma})+\mathcal{L}=\mathcal{S}+\mathcal{L}$ clearly holds. This proves one
direction. Conversely, suppose $\mathcal{Z}$ is a concentration space of $%
\mathfrak{C}^{\natural }$. Then $\mathcal{Z}=\limfunc{span}(\breve{\Sigma})$
for some singular matrix that is the limit of some sequence $\Sigma
_{m}^{\natural }\in \mathfrak{C}^{\natural }$, where $\Sigma _{m}^{\natural
}=\mathcal{L}(\Sigma _{m})+\Pi _{\mathcal{L}}$ for some $\Sigma _{m}\in 
\mathfrak{C}$. By the same compactness argument as before, we have $\mathcal{%
L}(\Sigma _{m_{i}})\rightarrow \bar{\Sigma}$ implying that $\bar{\Sigma}\in 
\limfunc{cl}(\mathcal{L}(\mathfrak{C}))$. Furthermore, we immediately arrive
at $\breve{\Sigma}=\bar{\Sigma}+\Pi _{\mathcal{L}}$. As before it follows
that $\limfunc{rank}(\bar{\Sigma})<n-l$ must hold and hence that $\mathcal{S}%
=\limfunc{span}(\bar{\Sigma})\in \mathbb{J}(\mathcal{L},\mathfrak{C})$. But
then $\mathcal{Z}=\limfunc{span}(\breve{\Sigma})=\limfunc{span}(\bar{\Sigma}%
+\Pi _{\mathcal{L}})=\limfunc{span}(\bar{\Sigma})+\mathcal{L}$ holds,
implying the result. $\blacksquare $

\begin{remark}
(i) By construction $\mathbb{J}(\mathcal{L},\mathfrak{C})=\mathbb{J}(%
\mathcal{L},\mathfrak{C}^{\sharp })=\mathbb{J}(\mathcal{L},\mathfrak{C}%
^{\natural })$. Furthermore, all three collections coincide with the
collection of all concentration spaces of $\mathfrak{C}^{\sharp }$ (the
union over which is $J(\mathfrak{C}^{\sharp })$ in the notation of \cite%
{PP2016}).

(ii) The sum $\mathcal{S}+\mathcal{L}$ is an orthogonal sum and hence $%
\mathcal{S}$ is uniquely determined.

(iii) The map $\Sigma \mapsto \Sigma ^{\sharp }$ is surjective from $%
\mathfrak{C}$ to $\mathfrak{C}^{\sharp }$ by definition, and the analogous
statement holds for the map $\Sigma \mapsto \Sigma ^{\natural }$. But these
maps need not be injective.
\end{remark}

\begin{lemma}
\label{change cov-model}Let $\mathfrak{C}$ be a covariance model and let $%
\mathcal{L}$ be a linear subspace of $\mathbb{R}^{n}$ with $\dim (\mathcal{L}%
)<n$. Furthermore, let $W\subseteq \mathbb{R}^{n}$ be a rejection region of
a test, which is $G(a+\mathcal{L})$-invariant for some $a\in \mathbb{R}^{n}$%
. Then for every $\sigma $, $0<\sigma <\infty $, and every $\Sigma \in 
\mathfrak{C}$ we have%
\begin{equation*}
P_{a,\sigma ^{2}\Sigma }(W)=P_{a,\sigma ^{2}\mathcal{L}(\Sigma
)}(W)=P_{a,\sigma ^{2}\Sigma ^{\sharp }}(W)=P_{a,\sigma ^{2}\Sigma
^{\natural }}(W).
\end{equation*}%
Furthermore, these probabilities do not depend on $\sigma $ and they are
unaffected if $a$ is replaced by an arbitrary element of $a+\mathcal{L}$.
\end{lemma}

\textbf{Proof:} The first claim is essentially proved by the argument
establishing (B.1) in Appendix B of \cite{PP3}. The second claim is an
immediate consequence of the assumed invariance (cf. also Proposition 5.4 in 
\cite{PP2016}). $\blacksquare $

\textbf{Proof of Theorem \ref{size_one}:} By monotonicity w.r.t. $C$ we may
assume $C>0$. Note that $\dim (\mathfrak{M}_{0}^{lin})=k-q<n$ by our general
model assumptions. Since $T$ is $G(\mathfrak{M}_{0})$-invariant by Lemma
5.16 in \cite{PP2016}, the preceding Lemma \ref{change cov-model}, applied
with $\mathcal{L}=\mathfrak{M}_{0}^{lin}$ and $a=\mu _{0}$, hence shows that
it suffices to prove the theorem with $\mathfrak{C}$ replaced by $\mathfrak{C%
}^{\sharp }$. By Lemma \ref{associated_cov_model},\ also applied with $%
\mathcal{L}=\mathfrak{M}_{0}^{lin}$, the space $\mathcal{S}$ appearing in
the formulation of the theorem is a concentration space of $\mathfrak{C}%
^{\sharp }$. We now apply Part 3 of Corollary 5.17 of \cite{PP2016} to the
linear model (\ref{lm}) considered in the present paper, but with $\mathfrak{%
C}$ replaced by $\mathfrak{C}^{\sharp }$. All assumptions of that result,
except for the assumption that $\check{\Omega}(z)=0$ and $R\check{\beta}%
(z)\neq 0$ simultaneously hold $\lambda _{\mathcal{S}}$-almost everywhere,
are easily seen to be satisfied. We verify the remaining assumption now as
follows: The discussion following (27) in Section 5.4 of \cite{PP2016} shows
that in case $N=\emptyset $ (which is assumed here) $\check{\Omega}(z)=0$
holds for every $z\in \mathrm{\limfunc{span}}(X)$, and thus for every $z\in 
\mathcal{S}$ (since $\mathcal{S}\subseteq \mathrm{\limfunc{span}}(X)$ has
been assumed). Hence, $\check{\Omega}(z)=0$ $\lambda _{\mathcal{S}}$-almost
everywhere follows (note that $\lambda _{\mathcal{S}}(\mathbb{R}%
^{n}\backslash \mathcal{S})=0$ trivially holds). Furthermore, Assumption \ref%
{Ass_5} together with $N=\emptyset $ imply that $\check{\beta}(X\gamma )=%
\check{\beta}(\varepsilon \cdot 0+X\gamma )=\varepsilon \check{\beta}%
(0)+\gamma $ for every $\gamma \in \mathbb{R}^{k}$ and every $\varepsilon
\neq 0$, which of course implies $\check{\beta}(X\gamma )=\gamma $ for every 
$\gamma \in \mathbb{R}^{k}$. Since we have assumed $\mathcal{S}\subseteq 
\mathrm{\limfunc{span}}(X)$, it follows on the one hand that for every $z\in 
\mathcal{S}$ we have $R\check{\beta}(z)=0$ if and only if $z\in \mathfrak{M}%
_{0}^{lin}$. On the other hand, by construction $\mathcal{S}\subseteq (%
\mathfrak{M}_{0}^{lin})^{\bot }$ holds, showing that $R\check{\beta}(z)\neq
0 $ must hold for all nonzero $z\in \mathcal{S}$ in view of the fact that $%
\mathcal{S}\subseteq \mathrm{\limfunc{span}}(X)$ has been assumed. Since $%
\mathcal{S}$ can not be zero-dimensional in view of its definition (cf. the
discussion in \cite{PP3} following Definition 5.1), $\lambda _{\mathcal{S}%
}(\left\{ 0\right\} )=0$ follows, which completes the proof (since $\lambda
_{\mathcal{S}}(\mathbb{R}^{n}\backslash \mathcal{S})=0$ trivially holds). $%
\blacksquare $

\textbf{Proof of Corollary \ref{suff&nec}:} Necessity follows immediately
from Theorem \ref{size_one}. For sufficiency we apply Corollary 5.6 in \cite%
{PP3} with $\mathcal{V}=\{0\}$, i.e., with $\mathcal{L}=\mathfrak{M}%
_{0}^{lin}$: Observe that $\dim (\mathcal{L})=k-q<n$ holds, and that $T$ and 
$N^{\dag }=N^{\ast }$ satisfy the assumptions of this corollary in view of
Lemma 5.16 in the same reference. Since $N^{\ast }=\mathrm{\limfunc{span}}%
(X) $ is assumed, the condition $\mathcal{S}\nsubseteq \mathrm{\limfunc{span}%
}(X) $ for every $\mathcal{S}\in $ $\mathbb{J}(\mathfrak{M}_{0}^{lin},%
\mathfrak{C})$ implies $\mu _{0}+\mathcal{S}\nsubseteq N^{\ast }=N^{\dag }$
for every $\mu _{0}\in \mathfrak{M}_{0}$ (as $\mathrm{\limfunc{span}}(X)$ is
obviously invariant under addition of elements $\mu _{0}\in \mathfrak{M}_{0}$%
) and for every $\mathcal{S}\in $ $\mathbb{J}(\mathfrak{M}_{0}^{lin},%
\mathfrak{C})$. An application of Corollary 5.6 in \cite{PP3} now delivers (%
\ref{size_control}). $\blacksquare $

\begin{theorem}
\label{size_one_extension}Let $\mathfrak{C}$ be a covariance model. Let $T$
be a nonsphericity-corrected F-type test statistic of the form (\ref{F-type}%
) based on $\check{\beta}$ and $\check{\Omega}$ satisfying Assumption \ref%
{Ass_5}. Assume further that $q=1$, that $\check{\beta}=\hat{\beta}_{X}$,
and that $\check{\Omega}(y)$ is nonnegative definite for every $y\in \mathbb{%
R}^{n}\backslash N$. Suppose there exists an $\mathcal{S}\in \mathbb{J}(%
\mathfrak{M}_{0}^{lin},\mathfrak{C})$ with the property that $s\in \mathbb{R}%
^{n}\backslash N$ and $s\in N^{\ast }$ hold for $\lambda _{\mathcal{S}}$%
-almost all $s\in \mathcal{S}$. Furthermore, assume that $\mathcal{S}$ is
not orthogonal to $\mathrm{\limfunc{span}}(X)$. Then (\ref{eq_size_one})
holds for every critical value $C$, $-\infty <C<\infty $, for every $\mu
_{0}\in \mathfrak{M}_{0}$, and for every $\sigma ^{2}\in (0,\infty )$.
\end{theorem}

\textbf{Proof:} The proof proceeds as the proof of Theorem \ref{size_one} up
to the point where Part 3 of Corollary 5.17 of \cite{PP2016} is applied to
the linear model (\ref{lm}), but with $\mathfrak{C}$ replaced by $\mathfrak{C%
}^{\sharp }$. Here now all assumptions of this result in \cite{PP2016} are
easily seen to be satisfied, except for (i) $\check{\Omega}(s)=0$ $\lambda _{%
\mathcal{S}}$-almost everywhere, and (ii) $R\check{\beta}(s)\neq 0$ $\lambda
_{\mathcal{S}}$-almost everywhere. Since $s\in N^{\ast }$ hold for $\lambda
_{\mathcal{S}}$-almost all $s\in \mathcal{S}$ by assumption, we have that $%
\check{\Omega}(s)$ is singular for $\lambda _{\mathcal{S}}$-almost all $s\in 
\mathcal{S}$. But this implies $\check{\Omega}(s)=0$ for $\lambda _{\mathcal{%
S}}$-almost all $s\in \mathcal{S}$ since $q=1$ has been assumed. Since
trivially $\lambda _{\mathcal{S}}(\mathbb{R}^{n}\backslash \mathcal{S})=0$,
this verifies (i). We turn to (ii): Let $s\in \mathcal{S}$. Note that then $%
s\in (\mathfrak{M}_{0}^{lin})^{\bot }$ by construction of $\mathcal{S}$. But
then%
\begin{equation*}
\Pi _{\mathrm{\limfunc{span}}(X)}s=s-\Pi _{(\mathrm{\limfunc{span}}%
(X))^{\bot }}s
\end{equation*}%
belongs to $(\mathfrak{M}_{0}^{lin})^{\bot }$ since $\Pi _{(\mathrm{\limfunc{%
span}}(X))^{\bot }}s\in (\mathrm{\limfunc{span}}(X))^{\bot }\subseteq (%
\mathfrak{M}_{0}^{lin})^{\bot }$. Now,%
\begin{eqnarray*}
R\check{\beta}(s) &=&R\hat{\beta}_{X}(s)=R\hat{\beta}_{X}(\Pi _{(\mathrm{%
\limfunc{span}}(X))^{\bot }}s)+R\hat{\beta}_{X}(\Pi _{\mathrm{\limfunc{span}}%
(X)}s) \\
&=&R(X^{\prime }X)^{-1}X^{\prime }\Pi _{(\mathrm{\limfunc{span}}(X))^{\bot
}}s+R\hat{\beta}_{X}(\Pi _{\mathrm{\limfunc{span}}(X)}s)=R\hat{\beta}%
_{X}(\Pi _{\mathrm{\limfunc{span}}(X)}s).
\end{eqnarray*}%
Hence, $R\check{\beta}(s)=0$ if and only if $R\hat{\beta}_{X}(\Pi _{\mathrm{%
\limfunc{span}}(X)}s)=0$, which in turn is equivalent to $\Pi _{\mathrm{%
\limfunc{span}}(X)}s\in \mathfrak{M}_{0}^{lin}$ (since $\Pi _{\mathrm{%
\limfunc{span}}(X)}s\in \mathrm{\limfunc{span}}(X)$). But since $\Pi _{%
\mathrm{\limfunc{span}}(X)}s$ also belongs to $(\mathfrak{M}%
_{0}^{lin})^{\bot }$ as shown before, we conclude that $R\check{\beta}(s)=0$
holds if and only if $\Pi _{\mathrm{\limfunc{span}}(X)}s=0$. As a
consequence,%
\begin{equation*}
\left\{ s\in \mathcal{S}:R\check{\beta}(s)=0\right\} =\left\{ s\in \mathcal{S%
}:\Pi _{\mathrm{\limfunc{span}}(X)}s=0\right\} =\mathcal{S}\cap \ker (\Pi _{%
\mathrm{\limfunc{span}}(X)}).
\end{equation*}%
This is a proper linear subspace of $\mathcal{S}$ except in case $\mathcal{S}%
\subseteq \ker (\Pi _{\mathrm{\limfunc{span}}(X)})$, which, however, is
impossible by the assumption that $\mathcal{S}$ is not orthogonal to $%
\mathrm{\limfunc{span}}(X)$. Hence, $R\check{\beta}(s)=0$ only occurs on a
proper linear subspace of $\mathcal{S}$, and hence on a subset of $\mathcal{S%
}$ that has $\lambda _{\mathcal{S}}$-measure zero. Since trivially $\lambda
_{\mathcal{S}}(\mathbb{R}^{n}\backslash \mathcal{S})=0$, this proves (ii)
and completes the proof. $\blacksquare $

\subsection{Some comments on Lemmata \protect\ref{associated_cov_model} and 
\protect\ref{change cov-model}\label{sec_comments}}

Lemmata \ref{associated_cov_model} and \ref{change cov-model} allow one to
derive results regarding the rejection probabilities under a covariance
model $\mathfrak{C}$ by working with a different, though related, covariance
model $\mathfrak{C}^{\sharp }$. [By Lemma \ref{associated_cov_model} this
related covariance model has the property that its concentration spaces in
the sense of \cite{PP2016} are precisely given by the elements $\mathcal{S}$
of $\mathbb{J}(\mathcal{L},\mathfrak{C})$.] A case in point is Theorem \ref%
{size_one} in Section \ref{Sec_general_results}, which provides a
\textquotedblleft size one\textquotedblright\ result for the covariance
model $\mathfrak{C}$, and which has been derived by applying Part 3 of
Corollary 5.17 in \cite{PP2016} to the covariance model $\mathfrak{C}%
^{\sharp }$, after an appeal to the aforementioned lemmata. In a similar
vein, one can combine other results of \cite{PP2016} with these lemmata, but
we do not spell this out here. Often this will lead to improvements over
what one obtains from a direct application of the respective result of \cite%
{PP2016} to the covariance model $\mathfrak{C}$. We illustrate this in the
following by comparing the result in Theorem \ref{size_one} with what one
gets if instead one works with the originally given $\mathfrak{C}$ and
directly applies Part 3 of Corollary 5.17 in \cite{PP2016} to $\mathfrak{C}$.

Suppose $\mathfrak{C}$ and $T$ are as in Theorem \ref{size_one} (again with $%
N=\emptyset $ and nonnegative definiteness of $\check{\Omega}(y)$ for every $%
y\in \mathbb{R}^{n}$). Applying Part 3 of Corollary 5.17 in \cite{PP2016} to
the originally given covariance model $\mathfrak{C}$ allows one to obtain
the following result: If a concentration space $\mathcal{Z}$ of $\mathfrak{C}
$ exists that satisfies $\mathcal{Z}\subseteq \mathrm{\limfunc{span}}(X)$
and $\mathcal{Z}\nsubseteq \mathfrak{M}_{0}^{lin}$, then (\ref{eq_size_one})
holds (for every $C$, every $\mu _{0}\in \mathfrak{M}_{0}$, and every $%
\sigma ^{2}\in (0,\infty )$). [To see this note that by Corollary 5.17 in 
\cite{PP2016} one only has to verify that $\check{\Omega}(z)=0$ and $R\check{%
\beta}(z)\neq 0$ hold $\lambda _{\mathcal{Z}}$-almost everywhere. The
argument for $\check{\Omega}(z)=0$ $\lambda _{\mathcal{Z}}$-a.e. is
identical to the corresponding argument given in the proof of Theorem \ref%
{size_one}. For the second claim a similar argument as in the proof of
Theorem \ref{size_one} shows that for $z\in \mathcal{Z}$ we have $R\check{%
\beta}(z)=0$ if and only if $z\in \mathfrak{M}_{0}^{lin}$. In other words, $R%
\check{\beta}(z)=0$ for $z\in \mathcal{Z}$ only occurs when $z\in \mathcal{Z}%
\cap \mathfrak{M}_{0}^{lin}$, which is a $\lambda _{\mathcal{Z}}$-null set,
since $\mathcal{Z}\nsubseteq \mathfrak{M}_{0}^{lin}$.]

We now show that Theorem \ref{size_one} is indeed at least as good a result
as the result obtained in the preceding paragraph. For this it suffices to
show that a concentration space $\mathcal{Z}$ of $\mathfrak{C}$ satisfying $%
\mathcal{Z}\subseteq \mathrm{\limfunc{span}}(X)$ and $\mathcal{Z}\nsubseteq 
\mathfrak{M}_{0}^{lin}$ gives rise to an element $\mathcal{S}\in \mathbb{J}(%
\mathfrak{M}_{0}^{lin},\mathfrak{C})$ satisfying the assumptions of Theorem %
\ref{size_one}: To see this, set $\mathcal{S}=\Pi _{(\mathfrak{M}%
_{0}^{lin})^{\bot }}\mathcal{Z}$ and observe that $\mathcal{S}\in $ $\mathbb{%
J}(\mathfrak{M}_{0}^{lin},\mathfrak{C})$ by Part 1 of Lemma B.3 in Appendix
B.1 of \cite{PP3} (since $\Pi _{(\mathfrak{M}_{0}^{lin})^{\bot }}\mathcal{Z}%
\neq \left\{ 0\right\} $ in view of $\mathcal{Z}\nsubseteq \mathfrak{M}%
_{0}^{lin}$, and since $\Pi _{(\mathfrak{M}_{0}^{lin})^{\bot }}\mathcal{Z}%
\neq (\mathfrak{M}_{0}^{lin})^{\bot }$ in view of$\ \mathcal{Z}\subseteq 
\mathrm{\limfunc{span}}(X)$, $\mathfrak{M}_{0}^{lin}\subseteq \mathrm{%
\limfunc{span}}(X)$, and $\limfunc{rank}(X)<n$). Furthermore, observe that $%
\mathcal{S}\subseteq \mathrm{\limfunc{span}}(X)$ must also hold, since $%
\mathcal{Z}\subseteq \mathrm{\limfunc{span}}(X)$ and $\mathfrak{M}%
_{0}^{lin}\subseteq \mathrm{\limfunc{span}}(X)$.

Theorem \ref{size_one} will sometimes actually give a strictly better result
for the following reason (at least for covariance models $\mathfrak{C}$ that
are bounded, an essentially costfree assumption in view of Remark 5.1(ii) in 
\cite{PP3}): Concentration spaces $\mathcal{Z}$ of $\mathfrak{C}$, that
satisfy $\mathcal{Z}\subseteq \mathrm{\limfunc{span}}(X)$ but also $\mathcal{%
Z}\subseteq \mathfrak{M}_{0}^{lin}$, can not be used in a direct application
of Part 3 of Corollary 5.17 in \cite{PP2016} since such spaces do not
satisfy the relevant assumptions (note that $R\check{\beta}(z)=0$ for all $%
z\in \mathcal{Z}$ holds for such spaces $\mathcal{Z}$); hence they do not
help in establishing a result of the form (\ref{eq_size_one}) via a direct
application of Part 3 of Corollary 5.17 in \cite{PP2016}. Nevertheless, such
concentration spaces can have associated with them spaces $\mathcal{S}\in $ $%
\mathbb{J}(\mathfrak{M}_{0}^{lin},\mathfrak{C})$ in the way as described in
Part 2 of Lemma B.3 in Appendix B.1 of \cite{PP3}, that then may allow one
to establish (\ref{eq_size_one}) via an application of Theorem \ref{size_one}
(provided the condition $\mathcal{S}\subseteq \mathrm{\limfunc{span}}(X)$
can be shown to hold).

\section{Appendix: Proofs and auxiliary results for Section \protect\ref%
{sec_stationary-results}}

\textbf{Proof of Theorem \ref{size_one_auto}: }First, that $\mathcal{S}%
\subseteq \mathrm{\limfunc{span}}(X)$ is equivalent to $\mathcal{A}\subseteq 
\mathrm{\limfunc{span}}(X)$ where $\mathcal{A}:=\limfunc{span}(E_{n,\rho
(\gamma _{1})}(\gamma _{1}),\ldots ,E_{n,\rho (\gamma _{p})}(\gamma _{p}))$
is obvious since any element of $\mathcal{A}$ is the sum of an element of $%
\mathcal{S}$ and an element of $\mathfrak{M}_{0}^{lin}\subseteq \mathrm{%
\limfunc{span}}(X)$. Second, $\mathcal{S}\subseteq \mathrm{\limfunc{span}}%
(X) $, $\mathfrak{M}_{0}^{lin}\subseteq \mathrm{\limfunc{span}}(X)$, and the
fact that $\mathcal{S}$ is certainly orthogonal to $\mathfrak{M}_{0}^{lin}$
imply $\limfunc{dim}(\mathcal{S})+\limfunc{dim}(\mathfrak{M}_{0}^{lin})\leq 
\limfunc{dim}(\mathrm{\limfunc{span}}(X))=k$. Since we always maintain $k<n$
we can conclude that $\limfunc{dim}(\mathcal{S})<n-\limfunc{dim}(\mathfrak{M}%
_{0}^{lin})$ must hold. This together with Proposition 6.1 of \cite{PP3} now
shows that the linear subspace $\mathcal{S}$ figuring in the theorem belongs
to $\mathbb{J}(\mathfrak{M}_{0}^{lin},\mathfrak{C}(\mathfrak{F}))$ as
clearly $\limfunc{dim}(\mathfrak{M}_{0}^{lin})=k-q<n$ holds. An application
of Theorem \ref{size_one} with $\mathfrak{C}=\mathfrak{C}(\mathfrak{F})$
then completes the proof. $\blacksquare $

\textbf{Proof of Lemma \ref{lemma:AR2}:} If $\{\gamma \}\in \mathbb{S}(%
\mathfrak{F},\mathcal{L})$ holds, the definition of $\mathbb{S}(\mathfrak{F},%
\mathcal{L})$ (Definition 6.4 in \cite{PP3}) immediately implies that $%
\kappa (\underline{\omega }(\mathcal{L}),\underline{d}(\mathcal{L}))+\kappa
(\gamma ,1)<n$ must hold. To prove the converse, we first claim that there
exists a sequence of spectral densities $f_{m}$ in $\mathfrak{F}$ so that
the sequence of spectral measures $\mathsf{m}_{g_{m}}$ defined by their
spectral densities 
\begin{equation*}
g_{m}(\nu )=|\Delta _{\underline{\omega }(\mathcal{L}),\underline{d}(%
\mathcal{L})}(e^{\iota \nu })|^{2}f_{m}(\nu )/\int_{-\pi }^{\pi }|\Delta _{%
\underline{\omega }(\mathcal{L}),\underline{d}(\mathcal{L})}(e^{\iota \nu
})|^{2}f_{m}(\nu )d\nu
\end{equation*}%
converges weakly to a spectral measure $\mathsf{m}$ that satisfies $\limfunc{%
supp}(\mathsf{m})\cap \lbrack 0,\pi ]=\{\gamma \}$. Here $\Delta _{%
\underline{\omega }(\mathcal{L}),\underline{d}(\mathcal{L})}$ is a certain
differencing operator given in Definition 6.3 of \cite{PP3} and $\limfunc{%
supp}(\mathsf{m})$ denotes the support of $\mathsf{m}$. To prove this claim,
let $\rho _{m}\in (0,1)$ converge to $1$ as $m\rightarrow \infty $, and let $%
\xi _{j}$ for $j\in \mathbb{N}$ be a sequence in $[0,\pi ]\backslash
\{0,\omega _{1}(\mathcal{L}),\ldots ,\omega _{p(\mathcal{L})}(\mathcal{L}%
),\pi \}$, where $\underline{\omega }(\mathcal{L})=(\omega _{1}(\mathcal{L}%
),\ldots ,\omega _{p(\mathcal{L})}(\mathcal{L}))$, that converges to $\gamma 
$ as $j\rightarrow \infty $. Now for every fixed $j\in \mathbb{N}$ the
sequence of spectral measures $\mathsf{m}_{h_{m,j}}$ with spectral density 
\begin{equation*}
h_{m,j}(\nu )=(2\pi )^{-1}\frac{(1-\rho _{m}^{2})((1+\rho
_{m}^{2})^{2}-4\rho _{m}^{2}\cos ^{2}(\xi _{j}))}{1+\rho _{m}^{2}}\left\vert
1-\rho _{m}e^{-\iota \xi _{j}}e^{-\iota \nu }\right\vert ^{-2}\left\vert
1-\rho _{m}e^{\iota \xi _{j}}e^{-\iota \nu }\right\vert ^{-2}
\end{equation*}%
converges weakly to $(\delta _{-\xi _{j}}+\delta _{\xi _{j}})/2$ as $%
m\rightarrow \infty $ (cf., e.g., the argument given in the proof of Lemma
G.2 in \cite{PP2016}). Note that $h_{m,j}\in \mathfrak{F}_{\mathrm{AR(}2%
\mathrm{)}}$ and thus $h_{m,j}\in \mathfrak{F}$. Since $\xi _{j}\notin
\{\omega _{1}(\mathcal{L}),\ldots ,\omega _{p}(\mathcal{L})\}$, we can
conclude that the map $\nu \mapsto \Delta _{\underline{\omega }(\mathcal{L}),%
\underline{d}(\mathcal{L})}(e^{\iota \nu })$ does not vanish on $\{-\xi
_{j},\xi _{j}\}$. It follows that the spectral measures $\mathsf{m}%
_{g_{m,j}} $ with spectral densities%
\begin{equation*}
g_{m,j}(\nu )=|\Delta _{\underline{\omega }(\mathcal{L}),\underline{d}(%
\mathcal{L})}(e^{\iota \nu })|^{2}h_{m,j}(\nu )/\int_{-\pi }^{\pi }|\Delta _{%
\underline{\omega }(\mathcal{L}),\underline{d}(\mathcal{L})}(e^{\iota \nu
})|^{2}h_{m,j}(\nu )d\nu
\end{equation*}%
also converge weakly to $(\delta _{-\xi _{j}}+\delta _{\xi _{j}})/2$, for
fixed $j$ and for $m\rightarrow \infty $. Since $(\delta _{-\xi _{j}}+\delta
_{\xi _{j}})/2$ certainly converges weakly to $(\delta _{-\gamma }+\delta
_{\gamma })/2$ as $j\rightarrow \infty $, a standard diagonal argument now
delivers a sequence $f_{m}=h_{m,j(m)}$ as required above, for $j(m)$ a
suitable subsequence of $j$. Together with the condition $\kappa (\underline{%
\omega }(\mathcal{L}),\underline{d}(\mathcal{L}))+\kappa (\gamma ,1)<n$ we
see that $\{\gamma \}\in \mathbb{S}(\mathfrak{F},\mathcal{L})$ follows. This
proves the first claim. The second claim is a trivial consequence of the
first claim, since $\kappa (\gamma ,1)=1$ for $\gamma =0,\pi $ and $\kappa
(\gamma ,1)=2$ for $\gamma \in (0,\pi )$. The third claim is seen as
follows: If $\{\gamma \}\in \mathbb{S}(\mathfrak{F},\mathcal{L})$, then
certainly $\gamma \in \bigcup \mathbb{S}(\mathfrak{F},\mathcal{L})$.
Conversely, let $\gamma \in \bigcup \mathbb{S}(\mathfrak{F},\mathcal{L})$.
Then $\gamma \in \Gamma $ for some $\Gamma \in \mathbb{S}(\mathfrak{F},%
\mathcal{L})$. By definition of $\mathbb{S}(\mathfrak{F},\mathcal{L})$, see
Definition 6.4 in \cite{PP3}, we have 
\begin{equation*}
\sum\nolimits_{\gamma ^{\prime }\in \Gamma }\kappa (\gamma ^{\prime
},1)<n-\kappa (\underline{\omega }(\mathcal{L}),\underline{d}(\mathcal{L})),
\end{equation*}%
implying that $\kappa (\gamma ,1)<n-\kappa (\underline{\omega }(\mathcal{L}),%
\underline{d}(\mathcal{L}))$ holds. But then $\{\gamma \}\in \mathbb{S}(%
\mathfrak{F},\mathcal{L})$ follows from the already established first claim. 
$\blacksquare $

\textbf{Proof of Theorem \ref{size_one_AR2}:} Since $\limfunc{span}%
(E_{n,\rho (\gamma )}(\gamma ))\subseteq \mathrm{\limfunc{span}}(X)$ but $%
\limfunc{span}(E_{n,\rho (\gamma )}(\gamma ))\nsubseteq \mathfrak{M}%
_{0}^{lin}\subseteq \mathrm{\limfunc{span}}(X)$ in view of the definition of 
$\rho (\gamma )$, it easily follows that%
\begin{equation*}
\kappa (\underline{\omega }(\mathfrak{M}_{0}^{lin}),\underline{d}(\mathfrak{M%
}_{0}^{lin}))+\kappa (\gamma ,1)\leq \kappa (\underline{\omega }(\mathrm{%
\limfunc{span}}(X)),\underline{d}(\mathrm{\limfunc{span}}(X)))
\end{equation*}%
must hold. The r.h.s. of the above inequality is now not larger than $k$ in
view of Lemma D.1 in Appendix D of \cite{PP3}. As we always maintain $k<n$,
the first claim follows. Because of the claim just established and since $%
\mathfrak{F}\supseteq \mathfrak{F}_{\mathrm{AR(}2\mathrm{)}}$, we conclude
from Lemma \ref{lemma:AR2} that $\{\gamma \}\in \mathbb{S}(\mathfrak{F},%
\mathfrak{M}_{0}^{lin})$ (note that $\dim (\mathfrak{M}_{0}^{lin})=k-q<n$
always holds). Set $\mathcal{S}=\limfunc{span}(\Pi _{(\mathfrak{M}%
_{0}^{lin})^{\bot }}E_{n,\rho (\gamma )}(\gamma ))$ and observe that $%
\mathcal{S}$ satisfies all the conditions of Theorem \ref{size_one_auto}
(recall that $\mathcal{S}\subseteq \mathrm{\limfunc{span}}(X)$ if and only
if $\limfunc{span}(E_{n,\rho (\gamma )}(\gamma ))\subseteq \mathrm{\limfunc{%
span}}(X)$ holds as noted in that theorem). An application of Theorem \ref%
{size_one_auto} then establishes (\ref{size}). $\blacksquare $

\begin{lemma}
\label{lemma:AR2+}For every $\gamma \in \lbrack 0,\pi ]$ and every $c>0$
there exists a sequence $h_{m}\in \mathfrak{F}_{\mathrm{AR(}2\mathrm{)}%
}^{ext}$ and a sequence $\sigma _{m}^{2}$ of positive real numbers such that 
\begin{equation}
\sigma _{m}^{2}\Pi _{(\mathfrak{M}_{0}^{lin})^{\bot }}\Sigma (h_{m})\Pi _{(%
\mathfrak{M}_{0}^{lin})^{\bot }}\rightarrow \Pi _{(\mathfrak{M}%
_{0}^{lin})^{\bot }}\left( E_{n,\rho (\gamma )}(\gamma )E_{n,\rho (\gamma
)}^{\prime }(\gamma )+cI_{n}\right) \Pi _{(\mathfrak{M}_{0}^{lin})^{\bot }}%
\text{ as }m\rightarrow \infty .  \label{hm-limit}
\end{equation}
\end{lemma}

\textbf{Proof:} Let $\gamma \in \lbrack 0,\pi ]$ and $c>0$ be given. For
ease of notation we set $\mathcal{L}=\mathfrak{M}_{0}^{lin}$ in the
remainder of the proof. We can use the argument in the proof of Lemma \ref%
{lemma:AR2} to obtain a sequence of spectral densities $f_{m}$ in $\mathfrak{%
F}_{\mathrm{AR(}2\mathrm{)}}$ so that the sequence $\mathsf{m}_{g_{m}}$ with
spectral density given by 
\begin{equation*}
g_{m}(\nu )=|\Delta _{\underline{\omega }(\mathcal{L}),\underline{d}(%
\mathcal{L})}(e^{\iota \nu })|^{2}f_{m}(\nu )/\int_{-\pi }^{\pi }|\Delta _{%
\underline{\omega }(\mathcal{L}),\underline{d}(\mathcal{L})}(e^{\iota \nu
})|^{2}f_{m}(\nu )d\nu
\end{equation*}%
converges weakly to the spectral measure $(\delta _{-\gamma }+\delta
_{\gamma })/2$. Now, set $e_{m}:=\int_{-\pi }^{\pi }|\Delta _{\underline{%
\omega }(\mathcal{L}),\underline{d}(\mathcal{L})}(e^{\iota \nu
})|^{2}f_{m}(\nu )d\nu $, which is a sequence of positive real numbers
(since $\Delta _{\underline{\omega }(\mathcal{L}),\underline{d}(\mathcal{L}%
)} $ is a polynomial and $f_{m}$ is nonzero a.e.). By Lemma D.2 in Appendix
D of \cite{PP3} we have%
\begin{align*}
& e_{m}^{-1}\Pi _{\mathcal{L}^{\bot }}\Sigma (f_{m})\Pi _{\mathcal{L}^{\bot
}} \\[7pt]
=& e_{m}^{-1}\Pi _{\mathcal{L}^{\bot }}H_{n}(\underline{\omega }(\mathcal{L}%
),\underline{d}(\mathcal{L}))\Sigma (\Delta _{\underline{\omega }(\mathcal{L}%
),\underline{d}(\mathcal{L})}\odot \mathsf{m}_{f_{m}},n-\kappa (\underline{%
\omega }(\mathcal{L}),\underline{d}(\mathcal{L})))H_{n}^{\prime }(\underline{%
\omega }(\mathcal{L}),\underline{d}(\mathcal{L}))\Pi _{\mathcal{L}^{\bot }}
\\
=& \Pi _{\mathcal{L}^{\bot }}H_{n}(\underline{\omega }(\mathcal{L}),%
\underline{d}(\mathcal{L}))\Sigma (\mathsf{m}_{g_{m}},n-\kappa (\underline{%
\omega }(\mathcal{L}),\underline{d}(\mathcal{L})))H_{n}^{\prime }(\underline{%
\omega }(\mathcal{L}),\underline{d}(\mathcal{L}))\Pi _{\mathcal{L}^{\bot }}
\\[7pt]
\rightarrow & \Pi _{\mathcal{L}^{\bot }}H_{n}(\underline{\omega }(\mathcal{L}%
),\underline{d}(\mathcal{L}))E_{n-\kappa (\underline{\omega }(\mathcal{L}),%
\underline{d}(\mathcal{L})),0}(\gamma )E_{n-\kappa (\underline{\omega }(%
\mathcal{L}),\underline{d}(\mathcal{L})),0}^{\prime }(\gamma )H_{n}^{\prime
}(\underline{\omega }(\mathcal{L}),\underline{d}(\mathcal{L}))\Pi _{\mathcal{%
L}^{\bot }}
\end{align*}%
as $m\rightarrow \infty $, where the convergence is due to weak convergence
of $m_{g_{m}}$ to $(\delta _{-\gamma }+\delta _{\gamma })/2$; see Appendix D
and Definition C.3 in Appendix C of \cite{PP3} for a definition of $H_{n}$, $%
\Sigma (\cdot ,\cdot )$, as well as $\odot $. Lemma D.3 in Appendix D of the
same reference now shows that the limit in the preceding display can be
written as 
\begin{equation*}
a\Pi _{\mathcal{L}^{\bot }}E_{n,\rho (\gamma )}(\gamma )E_{n,\rho (\gamma
)}^{\prime }(\gamma )\Pi _{\mathcal{L}^{\bot }}
\end{equation*}%
for some positive real number $a=a(\gamma )$. Now set $\sigma
_{m}^{2}=e_{m}^{-1}\left( a^{-1}+ce_{m}\right) $ and set%
\begin{equation*}
h_{m}=\left( a^{-1}f_{m}+(2\pi )^{-1}ce_{m}\right) /\left(
a^{-1}+ce_{m}\right) .
\end{equation*}%
Observe that $h_{m}\in \mathfrak{F}_{\mathrm{AR(}2\mathrm{)}}^{ext}$ holds.
But then 
\begin{equation*}
\sigma _{m}^{2}\Pi _{\mathcal{L}^{\bot }}\Sigma (h_{m})\Pi _{\mathcal{L}%
^{\bot }}=a^{-1}e_{m}^{-1}\Pi _{\mathcal{L}^{\bot }}\Sigma (f_{m})\Pi _{%
\mathcal{L}^{\bot }}+c\Pi _{\mathcal{L}^{\bot }}
\end{equation*}%
obtains, implying (\ref{hm-limit}). $\blacksquare $

\textbf{Proof of Theorem \ref{theo:AR2+}:} It suffices to prove the result
for $C>0$, which we henceforth assume. For ease of notation we set $\mathcal{%
L}=\mathfrak{M}_{0}^{lin}$ in the remainder of the proof. Let $\gamma \in
\lbrack 0,\pi ]$ satisfy $\limfunc{span}(E_{n,\rho (\gamma )}(\gamma
))\subseteq \mathrm{\limfunc{span}}(X)$. Observe that for $\mu _{0}\in 
\mathfrak{M}_{0}$, $0<\tau ^{2}<\infty $, and $h\in \mathfrak{F}_{\mathrm{AR(%
}2\mathrm{)}}^{ext}$ it holds that 
\begin{equation}
P_{\mu _{0},\tau ^{2}\Sigma (h)}(T\geq C)=P_{\mu _{0},\tau ^{2}\Pi _{%
\mathcal{L}^{\bot }}\Sigma (h)\Pi _{\mathcal{L}^{\bot }}}(T\geq C)=P_{\mu
_{0},\tau ^{2}\left[ \Pi _{\mathcal{L}^{\bot }}\Sigma (h)\Pi _{\mathcal{L}%
^{\bot }}+\Pi _{\mathcal{L}}\right] }(T\geq C).  \label{before}
\end{equation}%
This follows from $G(\mathfrak{M}_{0})$-invariance of $T$ and is proved in
the same way as is relation (B.1) in Appendix B of \cite{PP3}. Let now $c>0$
and fix $\mu _{0}\in \mathfrak{M}_{0}$, $0<\sigma ^{2}<\infty $. By Lemma %
\ref{lemma:AR2+} there exists a sequence $h_{m}\in \mathfrak{F}_{\mathrm{AR(}%
2\mathrm{)}}^{ext}$ and a sequence $\sigma _{m}^{2}$ of positive real
numbers such that 
\begin{equation*}
\sigma _{m}^{2}\Pi _{\mathcal{L}^{\bot }}\Sigma (h_{m})\Pi _{\mathcal{L}%
^{\bot }}+\Pi _{\mathcal{L}}\rightarrow \Pi _{\mathcal{L}^{\bot }}\left(
E_{n,\rho (\gamma )}(\gamma )E_{n,\rho (\gamma )}^{\prime }(\gamma
)+cI_{n}\right) \Pi _{\mathcal{L}^{\bot }}+\Pi _{\mathcal{L}},
\end{equation*}%
where the limit matrix is obviously nonsingular. Consequently,%
\begin{equation*}
P_{\mu _{0},\sigma _{m}^{2}\left[ \Pi _{\mathcal{L}^{\bot }}\Sigma
(h_{m})\Pi _{\mathcal{L}^{\bot }}+\Pi _{\mathcal{L}}\right] }\rightarrow
P_{\mu _{0},\Pi _{\mathcal{L}^{\bot }}E_{n,\rho (\gamma )}(\gamma )E_{n,\rho
(\gamma )}^{\prime }(\gamma )\Pi _{\mathcal{L}^{\bot }}+c\Pi _{\mathcal{L}%
^{\bot }}+\Pi _{\mathcal{L}}}
\end{equation*}%
for $m\rightarrow \infty $ in total variation norm (by an application of
Scheff\'{e}'s Lemma). By $G(\mathfrak{M}_{0})$-invariance of $T$ we also
have 
\begin{equation*}
P_{\mu _{0},\sigma ^{2}\Sigma (h_{m})}(T\geq C)=P_{\mu _{0},\sigma
_{m}^{2}\Sigma (h_{m})}(T\geq C),
\end{equation*}%
cf. Remark 5.5(iii) in \cite{PP2016}. Using (\ref{before}), the preceding
displays now imply that 
\begin{eqnarray*}
P_{\mu _{0},\sigma ^{2}\Sigma (h_{m})}(T\geq C) &=&P_{\mu _{0},\sigma
_{m}^{2}\left[ \Pi _{\mathcal{L}^{\bot }}\Sigma (h_{m})\Pi _{\mathcal{L}%
^{\bot }}+\Pi _{\mathcal{L}}\right] }(T\geq C) \\
&\rightarrow &P_{\mu _{0},\Pi _{\mathcal{L}^{\bot }}E_{n,\rho (\gamma
)}(\gamma )E_{n,\rho (\gamma )}^{\prime }(\gamma )\Pi _{\mathcal{L}^{\bot
}}+c\Pi _{\mathcal{L}^{\bot }}+\Pi _{\mathcal{L}}}(T\geq C).
\end{eqnarray*}%
The limit in the preceding display coincides -- using again $G(\mathfrak{M}%
_{0})$-invariance of $T$ similarly as in (\ref{before}) -- with 
\begin{equation*}
P_{\mu _{0},\sigma ^{2}\left[ E_{n,\rho (\gamma )}(\gamma )E_{n,\rho (\gamma
)}^{\prime }(\gamma )+cI_{n}\right] }(T\geq C).
\end{equation*}%
Since $\mathfrak{F}\supseteq \mathfrak{F}_{\mathrm{AR(}2\mathrm{)}}^{ext}$
has been assumed and since $c>0$ was arbitrary in the above discussion, it
follows that $\sup_{f\in \mathfrak{F}}P_{\mu _{0},\sigma ^{2}\Sigma
(f)}\left( T\geq C\right) $ is not smaller than $\sup_{\Sigma \in \mathfrak{C%
}(\gamma )}P_{\mu _{0},\sigma ^{2}\Sigma }(T\geq C)$ where $\mathfrak{C}%
(\gamma )$ denotes the auxiliary covariance model 
\begin{equation*}
\mathfrak{C}(\gamma )=\{E_{n,\rho (\gamma )}(\gamma )E_{n,\rho (\gamma
)}^{\prime }(\gamma )+cI_{n}:c>0\}.
\end{equation*}%
To prove the right-most inequality in (\ref{claim}) it hence suffices to
verify that for every $\mu _{0}\in \mathfrak{M}_{0}$ and every $0<\sigma
^{2}<\infty $ it holds that 
\begin{equation}
K(\gamma )\leq \sup_{\Sigma \in \mathfrak{C}(\gamma )}P_{\mu _{0},\sigma
^{2}\Sigma }(T\geq C).  \label{intermediate claim}
\end{equation}%
To this end, we shall use Theorem 5.19 of \cite{PP2016} applied to the
linear model (\ref{lm}) together with the covariance model $\mathfrak{C}%
(\gamma )$. Let $c_{m}$ be a sequence of positive real numbers satisfying $%
c_{m}\rightarrow 0$, and consider the corresponding sequence $\Sigma
_{m}=E_{n,\rho (\gamma )}(\gamma )E_{n,\rho (\gamma )}^{\prime }(\gamma
)+c_{m}I_{n}$ in $\mathfrak{C}(\gamma )$. Obviously $\Sigma _{m}\rightarrow
E_{n,\rho (\gamma )}(\gamma )E_{n,\rho (\gamma )}^{\prime }(\gamma )=:\bar{%
\Sigma}$ and $\limfunc{span}(\bar{\Sigma})=\limfunc{span}(E_{n,\rho (\gamma
)}(\gamma ))$ is $\kappa (\gamma ,1)$-dimensional. Note that $\kappa (\gamma
,1)$ is positive and that the $n\times n$-matrix $\bar{\Sigma}$ is singular
because the assumption $\limfunc{span}(E_{n,\rho (\gamma )}(\gamma
))\subseteq \mathrm{\limfunc{span}}(X)$ implies $\kappa (\gamma ,1)\leq k<n$%
. Next, observe that 
\begin{equation*}
\Pi _{\limfunc{span}(\bar{\Sigma})^{\bot }}\Sigma _{m}\Pi _{\limfunc{span}(%
\bar{\Sigma})^{\bot }}=\Pi _{\limfunc{span}(E_{n,\rho (\gamma )}(\gamma
))^{\bot }}\Sigma _{m}\Pi _{\limfunc{span}(E_{n,\rho (\gamma )}(\gamma
))^{\bot }}=c_{m}\Pi _{\limfunc{span}(\bar{\Sigma})^{\bot }},
\end{equation*}%
and that 
\begin{equation*}
\Pi _{\limfunc{span}(\bar{\Sigma})^{\bot }}\Sigma _{m}\Pi _{\limfunc{span}(%
\bar{\Sigma})}=0.
\end{equation*}%
Hence the additional assumption on $\Sigma _{m}$ appearing in Theorem 5.19
of \cite{PP2016} is satisfied with $s_{m}=c_{m}$ and $D=\Pi _{\limfunc{span}%
(E_{n,\rho (\gamma )}(\gamma ))^{\bot }}$. Note also that $\limfunc{span}(%
\bar{\Sigma})\subseteq \mathfrak{M}=\limfunc{span}(X)$ holds by our
assumption on $\gamma $.$\mathbb{\ }$Furthermore, since $\limfunc{span}(\bar{%
\Sigma})=\limfunc{span}(E_{n,\rho (\gamma )}(\gamma ))$ is not contained in $%
\mathcal{L}=\mathfrak{M}_{0}^{lin}$ in view of the definition of $\rho
(\gamma )$, it follows that there exists a $z\in \limfunc{span}(\bar{\Sigma}%
) $ so that $z\notin \mathcal{L}$. As both spaces are linear it even follows
that $z\notin \mathcal{L}$ is true for $\lambda _{\limfunc{span}(\bar{\Sigma}%
)}$-almost all $z\in \limfunc{span}(\bar{\Sigma})$. In view of the $\limfunc{%
span}(\bar{\Sigma})\subseteq \limfunc{span}(X)$, this implies that $R\hat{%
\beta}(z)\neq 0$ holds $\lambda _{\limfunc{span}(\bar{\Sigma})}$-almost
everywhere. Thus Theorem 5.19 of \cite{PP2016} is applicable, and delivers
(setting $Z=\bar{E}_{n,\rho (\gamma )}(\gamma )$ in that theorem) the claim (%
\ref{intermediate claim}), upon observing that in the definition of $\bar{\xi%
}(\gamma )$ in Theorem 5.19 of \cite{PP2016} and in the event following that
definition given in Theorem 5.19 of \cite{PP2016} one can replace $\bar{%
\Sigma}^{1/2}$ by $\Pi _{\limfunc{span}(\bar{\Sigma})}$ due to $\limfunc{span%
}(\bar{\Sigma})\subseteq \mathfrak{M}$, due to the equivariance property of $%
\check{\Omega}$ expressed in Assumption \ref{Ass_5}, and due to $G(\mathfrak{%
M})$-invariance of $N^{\ast }$ (and noting that in the case considered here $%
\Pi _{\limfunc{span}(\bar{\Sigma})}+D^{1/2}$ translates into $I_{n}$). It
remains to show the left-most inequality in (\ref{claim}). But this is
obvious upon noting that the event where $\check{\Omega}(\mathbf{G})$ is
nonnegative definite is contained in the event $\left\{ \bar{\xi}_{\gamma
}(x)\geq 0\right\} $ for every $x$. $\blacksquare $

\section{Appendix: Proofs for Section \protect\ref{sec_power}\label{App C}}

\textbf{Proof of Lemma \ref{improved_L_5.10}:} In view of $G(\mathfrak{M}%
_{0})$-invariance of $T$ we may set $\sigma ^{2}=1$. In case $\mathbb{K}$ is
empty there is nothing to prove. Hence assume $\mathbb{K\neq \varnothing }$.
To prove Part 1, observe that then $C^{\ast }(\mathbb{K})>-\infty $. Choose $%
C\in (-\infty ,C^{\ast }(\mathbb{K}))$. Since $C<C^{\ast }(\mathbb{K})$,
there exists an $\mathcal{S}\in \mathbb{K}$ with $C<C(\mathcal{S})\leq
C^{\ast }(\mathbb{K})$. Now repeat, with obvious modifications, the
arguments in the proof of Part 2 of Lemma 5.11 of \cite{PP3} that establish
(25) in that reference. To prove Part 2, observe that $C_{\ast }(\mathbb{K}%
)<\infty $, and choose $C\in (C_{\ast }(\mathbb{K}),\infty )$. Then there
exists an $\mathcal{S}\in \mathbb{K}$ with $C_{\ast }(\mathbb{K})\leq C(%
\mathcal{S})<C$. Now repeat, with obvious modifications, the arguments in
the proof of Part 3 of Lemma 5.11 of \cite{PP3}. $\blacksquare $

\begin{lemma}
\label{lem:simplify} Suppose the assumptions of Lemma \ref{improved_L_5.10}
are satisfied and suppose that $\mathbb{G}$ is a subset of $\mathbb{K}$ with
the property that for any $\mathcal{S}\in \mathbb{K}$ there is an element $%
\mathcal{S}^{\prime }\in \mathbb{G}$ such that $\mathcal{S}^{\prime
}\subseteq \mathcal{S}$ or $\mathcal{S}^{\prime }\supseteq \mathcal{S}$
holds. Then $C_{\ast }(\mathbb{K})=C_{\ast }(\mathbb{G})$ and $C^{\ast }(%
\mathbb{K})=C^{\ast }(\mathbb{G})$.
\end{lemma}

\textbf{Proof:} If $\mathbb{K}$ is empty, so is $\mathbb{G}$, and there is
nothing to prove. Hence, assume that $\mathbb{K}$ is nonempty. Then also $%
\mathbb{G}$ is nonempty. The claim will follow if we can show that for $%
\mathcal{S}_{1}\subseteq \mathcal{S}_{2}$, $\mathcal{S}_{1}\in \mathbb{K}$, $%
\mathcal{S}_{2}\in \mathbb{K}$, we have $C(\mathcal{S}_{1})=C(\mathcal{S}%
_{2})$. To this end fix $\mu _{0}\in \mathfrak{M}_{0}$ arbitrary. Since $%
\mathcal{S}_{1}\in \mathbb{K\subseteq H}$, we then have that $T(\mu
_{0}+s)=C(\mathcal{S}_{1})$ for $\lambda _{\mu _{0}+\mathcal{S}_{1}}$-almost
all $s\in \mathcal{S}_{1}$. Since also $\lambda _{\mu _{0}+\mathcal{S}%
_{1}}(N^{\dag })=0$ is assumed in Lemma \ref{improved_L_5.10}, we can find
an element $s_{1}\in \mathcal{S}_{1}$ such that $\mu _{0}+s_{1}\notin
N^{\dag }$ and such that $T(\mu _{0}+s_{1})=C(\mathcal{S}_{1})$. Since $%
\mathcal{S}_{2}\in \mathbb{K\subseteq H}$, the set%
\begin{equation*}
A=\left\{ \mu _{0}+s:s\in \mathcal{S}_{2},T(\mu _{0}+s)=C(\mathcal{S}%
_{2})\right\}
\end{equation*}%
is the complement in $\mu _{0}+\mathcal{S}_{2}$ of an $\lambda _{\mu _{0}+%
\mathcal{S}_{2}}$-null set. Hence, it intersects each neighborhood of $\mu
_{0}+s_{1}\in \mu _{0}+\mathcal{S}_{1}\subseteq \mu _{0}+\mathcal{S}_{2}$,
the neighborhood being relative to $\mu _{0}+\mathcal{S}_{2}$. Thus we may
choose a sequence $\mu _{0}+s(m)\in A$, such that $\mu _{0}+s(m)$ converges
to $\mu _{0}+s_{1}$ for $m\rightarrow \infty $. Since $\mu _{0}+s_{1}\notin
N^{\dag }$, it is a continuity point of $T$. Consequently, $T(\mu _{0}+s(m))$
converges to $T(\mu _{0}+s_{1})=C(\mathcal{S}_{1})$ for $m\rightarrow \infty 
$. But $T(\mu _{0}+s(m))=C(\mathcal{S}_{2})$ by the definition of $A$,
showing that $C(\mathcal{S}_{1})=C(\mathcal{S}_{2})$ must hold. $%
\blacksquare $

\begin{remark}
An example of such a collection $\mathbb{G}$ is provided by the set of all
minimal (maximal) elements of $\mathbb{K}$ w.r.t. inclusion. Note that this
set is well-defined as $\mathbb{K}$ is a collection of linear subspaces of $%
\mathbb{R}^{n}$.
\end{remark}

\textbf{Proof of Theorem \ref{neg_power}:} 1. Applying Part 1 of Lemma \ref%
{improved_L_5.10} with $\mathbb{K}=\{\mathcal{S}_{1},\mathcal{S}_{2}\}$
shows that $C$ satisfying (\ref{eqn:sup0}) must also satisfy $C\in \lbrack
C^{\ast }(\mathbb{K}),\infty )$. Since $C(\mathcal{S}_{1})\neq C(\mathcal{S}%
_{2})$ by assumption, it follows that $C_{\ast }(\mathbb{K})<C^{\ast }(%
\mathbb{K})$. Hence, we arrive at $C>C_{\ast }(\mathbb{K})$, which in view
of Part 2 of Lemma \ref{improved_L_5.10} implies (\ref{eqn:inf0}).

2. The same reasoning, but now with $\mathbb{K}=\{\mathcal{S}\}$, where $%
\mathcal{S}$ is as in the theorem, yields $C\geq C(\mathcal{S})$($=C_{\ast }(%
\mathbb{K})=C^{\ast }(\mathbb{K})$). Furthermore, note that $C>C(\mathcal{S}%
) $ obviously implies $C>C_{\ast }(\mathbb{K})$ and thus (\ref{eqn:inf0})
follows from Part 2 of Lemma \ref{improved_L_5.10}.

3. From $G(\mathfrak{M}_{0})$-invariance of $T$ (cf. Footnote \ref%
{footnote_3}) we know that (\ref{eqn:inf0}) implies%
\begin{equation}
\inf_{\Sigma \in \mathfrak{C}}P_{\mu _{0},\sigma ^{2}\Sigma }(T\geq C)=0
\label{star}
\end{equation}%
for every $\mu _{0}\in \mathfrak{M}_{0}$ and every $\sigma ^{2}$, $0<\sigma
^{2}<\infty $. Since $G(\mathfrak{M}_{0})$-invariance of $T$ implies $%
G(\left\{ \mu _{0}\right\} )$-almost invariance of $T$ for every $\mu
_{0}\in \mathfrak{M}_{0}$, (\ref{eqn:inf1}) now follows from (\ref{star})
together with Part 3 of Theorem 5.7 in \cite{PP2016}.\footnote{%
We note that the assumption in Theorem 5.7 of \cite{PP2016} that $\mathcal{Z}
$ is a concentration space is nowhere used in Part 3 of that theorem and its
proof.} Finally, (\ref{eqn:inf2}) follows immediately from (\ref{star}) by
noting that for every $\Sigma \in \mathfrak{C}$ and every $\sigma ^{2}\in
(0,\infty )$ the measures $P_{\mu _{1},\sigma ^{2}\Sigma }$ converge to $%
P_{\mu _{0},\sigma ^{2}\Sigma }$ in the total variation distance when $\mu
_{1}$ converges to $\mu _{0}$ (cf. the proof of Theorem 5.7, Part 2, in \cite%
{PP2016}). $\blacksquare $

\textbf{Proof of Corollary \ref{cor_neg_power}:} Set $\mathcal{V}=\{0\}$.
The assumptions on $T$ and on $N^{\dag }=N^{\ast }$ in the second and third
sentence of Theorem \ref{neg_power} are satisfied in view of Lemma 5.16 in 
\cite{PP3}. The assumption on the dimension of $\mathcal{L}:=\mathfrak{M}%
_{0}^{lin}$ is also satisfied since we always maintain $k<n$. If (\ref%
{eqn:sup0}) holds for a given $C$, Theorem \ref{size_one} implies that any $%
\mathcal{S}\in \mathbb{J}(\mathfrak{M}_{0}^{lin},\mathfrak{C})$ must satisfy 
$\mathcal{S}\nsubseteq \limfunc{span}(X)$; and thus $\mathcal{S}\nsubseteq
N^{\ast }$, since $N^{\ast }=\limfunc{span}(X)$ is assumed in the corollary.
Since $N^{\ast }$ is $G(\mathfrak{M})$-invariant (see Section \ref{Classes}%
), we also have $\mu _{0}+\mathcal{S}\nsubseteq N^{\ast }$ for every $\mu
_{0}\in \mathfrak{M}_{0}$. As $\mu _{0}+\mathcal{S}$ and $N^{\ast }$ are
affine subspaces of $\mathbb{R}^{n}$, this implies $\lambda _{\mu _{0}+%
\mathcal{S}}(N^{\ast })=0$ for every $\mu _{0}\in \mathfrak{M}_{0}$. Since $%
N^{\dag }$ coincides with $N^{\ast }$ for the class of test statistics
considered, we obtain that $N^{\dag }$ is a $\lambda _{\mu _{0}+\mathcal{S}}$%
-null set for every $\mu _{0}\in \mathfrak{M}_{0}$ and for every $\mathcal{S}%
\in \mathbb{J}(\mathfrak{M}_{0}^{lin},\mathfrak{C})$, and thus a fortiori
for every $\mathcal{S}\in \mathbb{H}$. We now see that Part 1 (Part 2,
respectively) follows from the corresponding parts of Theorem \ref{neg_power}
together with Part 3 of that theorem. $\blacksquare $

\textbf{Proof of Lemma \ref{lemma:eleH}:} Because of the assumption that $%
\mathfrak{F}$ contains $\mathfrak{F}_{\mathrm{AR(}2\mathrm{)}}$ and that $%
\dim (\mathcal{L})+1<n$, Lemma \ref{lemma:AR2} implies (cf. Remark \ref%
{rem:kappa and dim}(i)) that $\{\gamma \}\in \mathbb{S}(\mathfrak{F},%
\mathcal{L})$ for every $\gamma \in \{0,\pi \}$ (recall that $\kappa (\gamma
,1)=1$ for these $\gamma $'s). Furthermore, the dimension of 
\begin{equation*}
\mathcal{S}:=\limfunc{span}\left( \Pi _{\mathcal{L}^{\bot }}\left( E_{n,\rho
(\gamma ,\mathcal{L})}(\gamma )\right) \right)
\end{equation*}%
is $1$ (since the dimension of $\limfunc{span}(E_{n,\rho (\gamma ,\mathcal{L}%
)}(\gamma ))$ is $1$ for $\gamma \in \{0,\pi \}$ and since $E_{n,\rho
(\gamma ,\mathcal{L})}(\gamma )\nsubseteq \mathcal{L}$ in view of the
definition of $\rho (\gamma ,\mathcal{L})$). Therefore the dimension of $%
\mathcal{S}$ is smaller than $n-\dim (\mathcal{L})$, and it follows from
Proposition 6.1 in \cite{PP3} that $\mathcal{S}\in \mathbb{J}(\mathcal{L},%
\mathfrak{C}(\mathfrak{F}))$. $\blacksquare $

\section{Appendix: Auxiliary results and proofs for Section \protect\ref%
{sec_trends} \label{App_D}}

\textbf{Proof of Theorem \ref{size_one_polynomial}:} We first show that $%
\func{span}(E_{n,\rho (0)}(0))\subseteq \func{span}(X)$ is satisfied: For
any $i=1,\ldots ,k_{F}$ with $R_{\cdot i}\neq 0$, the $i$-th column of $F$
does not belong to $\mathfrak{M}_{0}^{lin}$. Observe that the $i$-th column
of $F$ spans $\func{span}(E_{n,i-1}(0))$. Hence $\rho (0)$ must satisfy $%
0\leq \rho (0)\leq k_{F}-1$. But then clearly $\func{span}(E_{n,\rho
(0)}(0))\subseteq \func{span}(F)\subseteq \func{span}(X)$. All the other
assumptions being obviously satisfied, Theorem \ref{size_one_AR2} completes
the proof. $\blacksquare $

\textbf{Proof of Theorem \ref{size_bound_polynomial}:} We apply Theorem \ref%
{theo:AR2+}. It suffices to verify that $\gamma =0$ satisfies the assumption 
$\limfunc{span}(E_{n,\rho (\gamma )}(\gamma ))\subseteq \mathrm{\limfunc{span%
}}(X)$ in that theorem. But this can be established exactly in the same way
as in the proof of Theorem \ref{size_one_polynomial}. It remains to verify
that $K(0)=P_{0,I_{n}}(R_{\cdot i_{0}}^{\prime }\check{\Omega}^{-1}R_{\cdot
i_{0}}\geq 0)$: Recall that $\kappa (0,1)=1$, and note that 
\begin{equation*}
\bar{\xi}_{0}(x)=x^{2}\bar{\xi}_{0}(1)\quad \text{ for every }x\in {\mathbb{R%
}}.
\end{equation*}%
This is trivial on the event $\{\mathbf{G}\in N^{\ast }\}$. On the
complement of this event, it follows from $\bar{E}_{n,\rho (0)}(0)$ being $%
n\times 1$-dimensional, and by using that $\hat{\beta}_{X}(\bar{E}_{n,\rho
(0)}(0)x)=x\hat{\beta}_{X}(\bar{E}_{n,\rho (0)}(0))$ holds for every $x\in {%
\mathbb{R}}$. From the equation in the previous display, we now obtain $%
K(0)=\Pr (\bar{\xi}_{0}(1)\geq 0)$. To prove the statement, we thus need to
show that $R\hat{\beta}_{X}(\bar{E}_{n,\rho (0)}(0))$ coincides with $%
R_{\cdot i_{0}}$, the first nonzero column of $R$. From a similar reasoning
as in the proof of Theorem 5.1, we see that $\bar{E}_{n,i}(0)=F_{\cdot
(i+1)} $ holds for $i=0,\ldots ,\rho (0)$. Hence, $\hat{\beta}_{X}(\bar{E}%
_{n,\rho (0)}(0))=e_{\rho (0)+1}(k)$ holds. Furthermore, from the definition
of $\rho (0)$, it follows that the first $\rho (0)$ columns of $R$ are zero,
and that the ($\rho (0)+1$)-th column of $R$ is nonzero. The statement
follows. $\blacksquare $

\begin{lemma}
\label{lem:auxA5} Let $H\in \mathbb{R}^{n\times n}$ be nonsingular and
define $\check{\beta}(y)=\hat{\beta}_{HX}(Hy)=(X^{\prime }H^{\prime
}HX)^{-1}X^{\prime }H^{\prime }Hy$. Let $\nu :\mathbb{R}^{n}\backslash
N^{\prime }\rightarrow \mathbb{R}$, for $N^{\prime }$ a subset of $\mathbb{R}%
^{n}$, and set 
\begin{equation*}
\check{\Omega}(y)=\nu (y)R(X^{\prime }H^{\prime }HX)^{-1}R^{\prime }\quad 
\text{ for every }y\notin N^{\prime }.
\end{equation*}%
Suppose that the following holds:

(a) $N^{\prime }$ is closed and $\lambda _{\mathbb{R}^{n}}(N^{\prime })=0$,

(b) $\delta y+X\eta \in \mathbb{R}^{n}\backslash N^{\prime }$ and $\nu
(\delta y+X\eta )=\delta ^{2}\nu (y)$ holds for every $y\in \mathbb{R}%
^{n}\backslash N^{\prime }$, every $\delta \neq 0$, and every $\eta \in 
\mathbb{R}^{k}$,

(c) $\nu $ is continuous on $\mathbb{R}^{n}\backslash N^{\prime }$,

(d) $\nu $ is $\lambda _{\mathbb{R}^{n}}$-almost everywhere nonzero on $%
\mathbb{R}^{n}\backslash N^{\prime }$.

Then $\check{\beta}$ and $\check{\Omega}$ satisfy Assumption \ref{Ass_5}
with $N=N^{\prime }$, and $\check{\Omega}$ satisfies Assumption \ref{Ass_7}.
Furthermore, if $\nu $ is nonnegative (positive) everywhere on $\mathbb{R}%
^{n}\backslash N^{\prime }$, then $\check{\Omega}$ is nonnegative (positive)
definite everywhere on $\mathbb{R}^{n}\backslash N^{\prime }$.
\end{lemma}

\textbf{Proof:} Obviously $\check{\beta}$ is well-defined and continuous on
all of $\mathbb{R}^{n}$, and thus also when restricted to $\mathbb{R}%
^{n}\backslash N^{\prime }$. Furthermore, $\check{\Omega}$ is clearly
well-defined and symmetric on $\mathbb{R}^{n}\backslash N^{\prime }$, and is
continuous on $\mathbb{R}^{n}\backslash N^{\prime }$ in view of (c). Since $%
N^{\prime }$ is a closed $\lambda _{\mathbb{R}^{n}}$-null set by (a), we
have verified Part (i) of Assumption \ref{Ass_5} with $N=N^{\prime }$. Part
(ii) of this assumption is contained in (b). That $\check{\beta}$ satisfies
the required equivariance property in Part (iii) of Assumption \ref{Ass_5}
is obvious. That $\check{\Omega}$ satisfies the required equivariance
property in that assumption follows immediately from (b), completing the
verification of Part (iii) of Assumption \ref{Ass_5}. Part (iv) in that
assumption follows from (d) together with $R(X^{\prime }H^{\prime
}HX)^{-1}R^{\prime }$ being positive definite. The same argument also shows
that $\check{\Omega}$ satisfies Assumption \ref{Ass_7}. The final statement
is trivial. $\blacksquare $

\begin{lemma}
\label{lem:classic} Suppose $\mathcal{W}$ is constant and symmetric, and
that $\Pi _{\limfunc{span}(X)^{\bot }}\mathcal{W}\Pi _{\limfunc{span}%
(X)^{\bot }}$ is nonzero. Then the estimators $\hat{\beta}$ and $\check{%
\Omega}_{\mathcal{W}}$ satisfy Assumption \ref{Ass_5} with $N=\emptyset $,
and $\check{\Omega}_{\mathcal{W}}$ satisfies Assumption \ref{Ass_7}. If,
additionally, $\Pi _{\limfunc{span}(X)^{\bot }}\mathcal{W}\Pi _{\limfunc{span%
}(X)^{\bot }}$ is nonnegative definite, then $\check{\Omega}_{\mathcal{W}%
}(y) $ is nonnegative definite for every $y\in \mathbb{R}^{n}$.
\end{lemma}

\textbf{Proof:} We verify (a)-(d) in Lemma \ref{lem:auxA5} for $H=I_{n}$, $%
\nu =\hat{\omega}_{\mathcal{W}}$, and $N^{\prime }=\emptyset $. Obviously
(a) is satisfied, and (c) follows immediately from the constancy assumption
on $\mathcal{W}$, since $\nu =\hat{\omega}_{\mathcal{W}}$ can clearly be
written as a quadratic form in $y$. Concerning (d), note that $\hat{\omega}_{%
\mathcal{W}}(y)=0$ is equivalent to $y^{\prime }\Pi _{\limfunc{span}%
(X)^{\bot }}\mathcal{W}\Pi _{\limfunc{span}(X)^{\bot }}y=0$. In view of the
constancy assumption on $\mathcal{W}$, the subset of $\mathbb{R}^{n}$ on
which $\hat{\omega}_{\mathcal{W}}$ vanishes is the zero set of a
multivariate polynomial, in fact of a quadratic form, on $\mathbb{R}^{n}$.
Since the (constant) matrix $\Pi _{\limfunc{span}(X)^{\bot }}\mathcal{W}\Pi
_{\limfunc{span}(X)^{\bot }}$ is symmetric and nonzero, the polynomial under
consideration does not vanish everywhere on $\mathbb{R}^{n}$, implying that
the zero set is a $\lambda _{\mathbb{R}^{n}}$-null set. This completes the
verification of (d). That (b) is satisfied follows immediately from $\nu (y)=%
\hat{\omega}_{\mathcal{W}}(y)=n^{-1}y^{\prime }\Pi _{\limfunc{span}(X)^{\bot
}}\mathcal{W}\Pi _{\limfunc{span}(X)^{\bot }}y$, the constancy of $\mathcal{W%
}$, and from $\Pi _{\limfunc{span}(X)^{\bot }}(\delta y+X\eta )=\delta \Pi _{%
\limfunc{span}(X)^{\bot }}(y)$ for every $\delta \in \mathbb{R}$, every $%
y\in \mathbb{R}^{n}$ and every $\eta \in \mathbb{R}^{k}$. Now apply Lemma %
\ref{lem:auxA5}. Note that the final statement concerning nonnegative
definiteness follows from the last part of Lemma \ref{lem:auxA5}, since
nonnegative definiteness of $\Pi _{\limfunc{span}(X)^{\bot }}\mathcal{W}\Pi
_{\limfunc{span}(X)^{\bot }}$ obviously implies nonnegativity of $\hat{\omega%
}_{\mathcal{W}}$ on $\mathbb{R}^{n}$. $\blacksquare $

\textbf{Proof of Corollary \ref{cor:polylrv1}:} The statement follows upon
combining Lemma \ref{lem:classic} with Theorem \ref{size_one_polynomial}. $%
\blacksquare $

\textbf{Proof of Corollary \ref{cor:polylrv2}:} The first part of the
corollary follows upon combining Lemma \ref{lem:classic} with Theorem \ref%
{size_bound_polynomial} noting that $\check{\Omega}_{\mathcal{W}}(z)$ is
nonnegative definite if and only if $\hat{\omega}_{\mathcal{W}}(z)\geq 0$.
For the second statement, note that $R(X^{\prime }X)^{-1}R^{\prime }$ is
positive definite, and hence%
\begin{equation*}
\{T\geq 0\}=\{\hat{\omega}_{\mathcal{W}}\geq 0\}\cup \{R\hat{\beta}=r\},
\end{equation*}%
from which it follows (note that $\{y:R\hat{\beta}(y)=r\}$ is an affine
subspace of $\mathbb{R}^{n}$ that does not coincide with $\mathbb{R}^{n}$,
and is hence a $\lambda _{\mathbb{R}^{n}}$-null set) that $P_{\mu ,\sigma
^{2}I_{n}}(T\geq 0)$ coincides with $P_{\mu ,\sigma ^{2}I_{n}}(\hat{\omega}_{%
\mathcal{W}}\geq 0)$. For $C\geq 0$ we then have (using monotonicity w.r.t. $%
C$)%
\begin{equation}
\sup_{\mu \in \mathfrak{M}_{1}}\sup_{0<\sigma ^{2}<\infty }P_{\mu ,\sigma
^{2}I_{n}}(T\geq C)\leq \sup_{\mu \in \mathfrak{M}_{1}}\sup_{0<\sigma
^{2}<\infty }P_{\mu ,\sigma ^{2}I_{n}}(\hat{\omega}_{\mathcal{W}}\geq 0).
\label{eqn:powerbound}
\end{equation}%
But from the equivariance property $\hat{\omega}_{\mathcal{W}}(\delta
y+X\eta )=\delta ^{2}\hat{\omega}_{\mathcal{W}}(y)$ for $\delta \neq 0$, $%
y\in \mathbb{R}^{n}$ and $\eta \in \mathbb{R}^{k}$, which was established in
the proof of Lemma \ref{lem:classic}, it follows straightforwardly that $%
P_{\mu ,\sigma ^{2}I_{n}}(\hat{\omega}_{\mathcal{W}}\geq 0)=P_{0,I_{n}}(\hat{%
\omega}_{\mathcal{W}}\geq 0)$ holds for every $\mu \in \mathfrak{M}$ and
every $0<\sigma <\infty $. This completes the proof. $\blacksquare $

\begin{lemma}
\label{lem:AM} If $N_{\mathrm{AM}}\neq \mathbb{R}^{n}$, then the estimators $%
\hat{\beta}$ and $\check{\Omega}_{\mathcal{W}_{\mathrm{AM}}}$ satisfy
Assumption \ref{Ass_5} with $N=N_{\mathrm{AM}}$, and $\check{\Omega}_{%
\mathcal{W}_{\mathrm{AM}}}$ satisfies Assumption \ref{Ass_7}; furthermore $%
\check{\Omega}_{\mathcal{W}_{\mathrm{AM}}}(z)$ is positive definite for
every $z\in \mathbb{R}^{n}\backslash N_{\mathrm{AM}}$.
\end{lemma}

\textbf{Proof:} Observe that $\hat{\rho}$, $\tilde{\rho}$, $M_{\mathrm{AM}}$%
, $\mathcal{W}_{\mathrm{AM}}$, and $\hat{\omega}_{\mathcal{W}_{\mathrm{AM}}}$
are well-defined on $\mathbb{R}^{n}\backslash N_{\mathrm{AM}}$. We next
verify (a)-(d) in Lemma \ref{lem:auxA5} for $H=I_{n}$, $\nu =\hat{\omega}_{%
\mathcal{W}_{\mathrm{AM}}}$, and $N^{\prime }=N_{\mathrm{AM}}$. We start
with (a): Using arguments as in the proof of Lemma 3.9 in \cite{Prein2014b},
or in the proof of Lemma B.1 in \cite{PP2016}, it is not difficult to verify
that $N_{\mathrm{AM}}$ is an algebraic set. We leave the details to the
reader. This, and the assumption $N_{\mathrm{AM}}\neq \mathbb{R}^{n}$,
implies that $N_{\mathrm{AM}}$ is a closed $\lambda _{\mathbb{R}^{n}}$-null
set. To verify (c) in Lemma \ref{lem:auxA5} it suffices to establish
continuity of $\mathcal{W}_{\mathrm{AM}}$ on $\mathbb{R}^{n}\backslash N_{%
\mathrm{AM}}$, since $\hat{u}(y)$ is certainly continuous on $\mathbb{R}^{n}$%
. To achieve this note that, since $\hat{\rho}$ is obviously continuous on $%
\mathbb{R}^{n}\backslash N_{\mathrm{AM}}$, since $\hat{\rho}(y)\neq 1$ for $%
y\in \mathbb{R}^{n}\backslash N_{\mathrm{AM}}$, and since $\mathsf{A}(\cdot
) $ is continuous on $\mathbb{R}$, it suffices to verify that $[\kappa _{%
\mathrm{QS}}(|i-j|/M_{\mathrm{AM}})]_{i,j=1}^{n-1}$ is continuous on $%
\mathbb{R}^{n}\backslash N_{\mathrm{AM}}$. Now, $M_{\mathrm{AM}}$ is
certainly continuous on $\mathbb{R}^{n}\backslash N_{\mathrm{AM}}$ and $%
\kappa _{\mathrm{QS}}$ is continuous on $\mathbb{R}$. Hence, $[\kappa _{%
\mathrm{QS}}(|i-j|/M_{\mathrm{AM}})]_{i,j=1}^{n-1}$ is easily seen to be
continuous at every $y\in \mathbb{R}^{n}\backslash N_{\mathrm{AM}}$ that
satisfies $M_{\mathrm{AM}}(y)\neq 0$. For $y\in \mathbb{R}^{n}\backslash N_{%
\mathrm{AM}}$ satisfying $M_{\mathrm{AM}}(y)=0$ continuity of $[\kappa _{%
\mathrm{QS}}(|i-j|/M_{\mathrm{AM}})]_{i,j=1}^{n-1}$ follows from continuity
of $M_{\mathrm{AM}}$ on $\mathbb{R}^{n}\backslash N_{\mathrm{AM}}$ together
with $\kappa _{\mathrm{QS}}(x)\rightarrow 0$ as $|x|\rightarrow \infty $, $%
\kappa _{\mathrm{QS}}(0)=1$, and the convention $[\kappa _{\mathrm{QS}%
}(|i-j|/M_{\mathrm{AM}}(y))]_{i,j=1}^{n-1}=I_{n-1}$ for $y$ so that $M_{%
\mathrm{AM}}(y)=0$. That (b) in Lemma \ref{lem:auxA5} holds is easily seen
to follow from $\hat{u}(\delta y+X\eta )=\delta \hat{u}(y)$ for every $%
\delta \in \mathbb{R}$, every $y\in \mathbb{R}^{n}$ and every $\eta \in 
\mathbb{R}^{k}$, which in particular implies $\hat{\rho}(\delta y+X\eta )=%
\hat{\rho}(y)$ and $\tilde{\rho}(\delta y+X\eta )=\tilde{\rho}(y)$ for every 
$\delta \neq 0$, every $y\in \mathbb{R}^{n}\backslash N_{\mathrm{AM}}$ and
every $\eta \in \mathbb{R}^{k}$. Finally, note that (d) in Lemma \ref%
{lem:auxA5} is satisfied, because $\hat{\omega}_{\mathcal{W}_{\mathrm{AM}%
}}(y)>0$ holds if $y\in \mathbb{R}^{n}\backslash N_{\mathrm{AM}}$. The
latter follows from the well-known fact that $\left[ \kappa _{\mathrm{QS}%
}(|i-j|/M_{\mathrm{AM}}(y))\right] _{i,j=1}^{n-1}$ is positive definite in
case $M_{\mathrm{AM}}(y)$ is well-defined (recall that this matrix is
defined as $I_{n-1}$ in case $M_{\mathrm{AM}}(y)=0$), together with the
observation that $y\in \mathbb{R}^{n}\backslash N_{\mathrm{AM}}$ implies $%
\mathsf{A}(\hat{\rho}(y))\hat{u}(y)=\hat{v}(y)\neq 0$. Now apply Lemma \ref%
{lem:auxA5}. Note that the just established fact, that $\hat{\omega}_{%
\mathcal{W}_{\mathrm{AM}}}(y)>0$ holds if $y\in \mathbb{R}^{n}\backslash N_{%
\mathrm{AM}}$, also shows that the last part of Lemma \ref{lem:auxA5}
applies, and hence shows that $\check{\Omega}_{\mathcal{W}_{\mathrm{AM}}}(y)$
is positive definite for every $y\in \mathbb{R}^{n}\backslash N_{\mathrm{AM}%
} $. $\blacksquare $

\textbf{Proof of Corollary \ref{cor:AM}:} This follows upon combining Lemma %
\ref{lem:AM} and Theorem \ref{size_bound_polynomial}, noting that the lower
bound obtained via Theorem \ref{size_bound_polynomial} equals $1$ due to
nonnegative definiteness of $\check{\Omega}_{\mathcal{W}_{\mathrm{AM}}}(y)$
for every $y\in \mathbb{R}^{n}\backslash N_{\mathrm{AM}}$, which is the
complement of a $\lambda _{\mathbb{R}^{n}}$-null set. $\blacksquare $

\begin{lemma}
\label{lem:vogel} Let $V\in \{A,I_{n}\}$, $c\in \mathbb{R}$, let $i\in
\{1,2\}$, and let $U$ be an $n\times m$-dimensional matrix with $m\geq 1$
such that $(X,U)$ is of full column-rank $k+m<n$. Then the estimators $%
\check{\beta}_{V}$ and $\check{\Omega}_{c,U,i,V}^{\mathrm{Vo}}$ satisfy
Assumption \ref{Ass_5} with $N=\limfunc{span}(X,U)$, and $\check{\Omega}%
_{c,U,i,V}^{\mathrm{Vo}}$ also satisfies Assumption \ref{Ass_7}; furthermore 
$\check{\Omega}_{c,U,i,V}^{\mathrm{Vo}}$ is positive definite on $\mathbb{R}%
^{n}\backslash \limfunc{span}(X,U)$.
\end{lemma}

\textbf{Proof:} We verify (a)-(d) in Lemma \ref{lem:auxA5} for $H=V$ (which
is invertible), $\nu =n^{j(V)}s_{A,X}^{2}\exp (cJ_{n,U}^{i})$, and $%
N^{\prime }=\limfunc{span}(X,U)$. By assumption, $k+m<n$, hence $\limfunc{%
span}(X,U)$ is a closed $\lambda _{\mathbb{R}^{n}}$-null set, showing that
(a) in Lemma \ref{lem:auxA5} is satisfied. Next, note that $s_{A,X}^{2}$, $%
s_{I_{n},(X,U)}^{2}$, and $s_{A,(X,U)}^{2}$ are well-defined and continuous
on $\mathbb{R}^{n}$; and that $J_{n,U}^{1}$ and $J_{n,U}^{2}$ are
well-defined and continuous on the set where $s_{I_{n},(X,U)}^{2}$ and $%
s_{A,(X,U)}^{2}$ are nonzero, respectively. Obviously, $%
s_{I_{n},(X,U)}^{2}(y)=0$ if and only if $y\in \limfunc{span}(X,U)$.
Similarly, $s_{A,(X,U)}^{2}(y)=0$ if and only if $Ay\in \limfunc{span}%
(A(X,U))$, or equivalently, $y\in \limfunc{span}(X,U)$. Hence (c) in Lemma %
\ref{lem:auxA5} follows. For (b) note first that obviously $\delta y+X\eta
\notin \limfunc{span}(X,U)$ holds for every $y\notin \limfunc{span}(X,U)$,
every $\delta \neq 0$ and every $\eta \in \mathbb{R}^{k}$. Second,
concerning the equivariance property of $\nu $, we note that for every $y\in 
\mathbb{R}^{n}$, every $\delta \in \mathbb{R}$, and every $\eta \in \mathbb{R%
}^{k}$ 
\begin{align}
& s_{A,X}^{2}(\delta y+X\eta )=\delta ^{2}s_{A,X}^{2}(y)  \label{eqn:s1} \\
& s_{I_{n},(X,U)}^{2}(\delta y+X\eta )=\delta ^{2}s_{I_{n},(X,U)}^{2}(y)
\label{eqn:s2} \\
& s_{A,(X,U)}^{2}(\delta y+X\eta )=\delta ^{2}s_{A,(X,U)}^{2}(y).
\label{eqn:s3}
\end{align}%
From Equations (\ref{eqn:s1})-(\ref{eqn:s3}) we hence see that the required
equivariance property follows if we can show that 
\begin{equation}
J_{n,U}^{i}(\delta y+X\eta )=J_{n,U}^{i}(y)\text{ for every }y\in \mathbb{R}%
^{n}\backslash \limfunc{span}(X,U),\text{ every }\delta \neq 0,\text{ and
every }\eta \in \mathbb{R}^{k}.  \label{eqn:Jinv}
\end{equation}%
To see this let $y\in \mathbb{R}^{n}\backslash \limfunc{span}(X,U)$, $\delta
\neq 0$, and $\eta \in \mathbb{R}^{k}$. We consider first the case where $%
i=1 $. Note that $G\hat{\beta}_{(X,U)}(\delta y+X\eta )=\delta G\hat{\beta}%
_{(X,U)}(y)$, and recall from (\ref{eqn:s2}) that $s_{I_{n},(X,U)}^{2}(%
\delta y+X\eta )=\delta ^{2}s_{I_{n},(X,U)}^{2}(y)>0$ (positivity following
from $y\notin \limfunc{span}(X,U)$), showing that $J_{n,U}^{1}(\delta
y+X\eta )=J_{n,U}^{1}(y)$. For $i=2$, note that $G\hat{\beta}%
_{A(X,U)}(A(\delta y+X\eta ))=\delta G\hat{\beta}_{A(X,U)}(Ay)$, and recall
from (\ref{eqn:s3}) that $s_{A,(X,U)}^{2}(\delta y+X\eta )=\delta
^{2}s_{A,(X,U)}^{2}(y)>0$ (positivity following from $y\notin \limfunc{span}%
(X,U)$), showing that $J_{n,U}^{2}(\delta y+X\eta )=J_{n,U}^{2}(y)$. This
verifies the statement in (\ref{eqn:Jinv}) and thus (b). Concerning (d) (and
the final claim in the lemma) note that for $y\notin \limfunc{span}(X,U)$ it
holds that $s_{A,X}^{2}(y)\exp (cJ_{n,U}^{i}(y))>0$. $\blacksquare $

\textbf{Proof of Corollary \ref{cor:vogel}:} This follows upon combining
Lemma \ref{lem:vogel} and Theorem \ref{size_bound_polynomial}, noting that
the lower bound obtained via Theorem \ref{size_bound_polynomial} equals $1$
due to nonnegative definiteness of $\check{\Omega}_{c,U,i,V}^{\mathrm{Vo}}$
on the complement of the $\lambda _{\mathbb{R}^{n}}$-null set $\limfunc{span}%
(X,U)$. $\blacksquare $

\begin{lemma}
\label{lem:vogelbunzel} Suppose that $\mathcal{W}$ is constant and
symmetric, that $\Pi _{\limfunc{span}(X)^{\bot }}\mathcal{W}\Pi _{\limfunc{%
span}(X)^{\bot }}$ is nonzero, and that $c\in \mathbb{R}$. Then the
following holds:

\begin{enumerate}
\item If $U$ is an $n\times m$-dimensional matrix with $m\geq 1$ such that $%
(X,U)$ is of full column-rank $k+m<n$, then the estimators $\hat{\beta}$ and 
$\check{\Omega}_{\mathcal{W},U,c}^{\mathrm{BV},J}$ satisfy Assumption \ref%
{Ass_5} with $N=\limfunc{span}(X,U)$, and $\check{\Omega}_{\mathcal{W},U,c}^{%
\mathrm{BV},J}$ satisfies Assumption \ref{Ass_7}. If, additionally, $\Pi _{%
\limfunc{span}(X)^{\bot }}\mathcal{W}\Pi _{\limfunc{span}(X)^{\bot }}$ is
nonnegative definite, then $\check{\Omega}_{\mathcal{W},U,c}^{\mathrm{BV},J}$
is nonnegative definite on $\mathbb{R}^{n}\backslash \limfunc{span}(X,U)$.

\item The estimators $\hat{\beta}$ and $\check{\Omega}_{\mathcal{W},c}^{%
\mathrm{BV}}$ satisfy Assumption \ref{Ass_5} with $N=\limfunc{span}(X)$, and 
$\check{\Omega}_{\mathcal{W},c}^{BV}$ satisfies Assumption \ref{Ass_7}. If,
additionally, $\Pi _{\limfunc{span}(X)^{\bot }}\mathcal{W}\Pi _{\limfunc{span%
}(X)^{\bot }}$ is nonnegative definite, then $\check{\Omega}_{\mathcal{W}%
,c}^{\mathrm{BV}}$ is nonnegative definite on $\mathbb{R}^{n}\backslash 
\limfunc{span}(X)$.
\end{enumerate}
\end{lemma}

\textbf{Proof:} 1. We verify (a)-(d) in Lemma \ref{lem:auxA5} for $H=I_{n}$, 
$\nu =\hat{\omega}_{\mathcal{W}}\exp (cJ_{n,U}^{1})$, and $N^{\prime }=%
\limfunc{span}(X,U)$. That (a) holds follows from the same argument as in
the proof of Lemma \ref{lem:vogel}. That (c) holds, follows from continuity
of $\hat{\omega}_{\mathcal{W}}$ on $\mathbb{R}^{n}$ (cf. the proof of Lemma %
\ref{lem:classic}), together with continuity of $\exp (cJ_{n,U}^{1})$ on the
complement of $\limfunc{span}(X,U)$ (cf. the proof of Lemma \ref{lem:vogel}%
). The first Part of (b) was established in the proof of Lemma \ref%
{lem:vogel}. The second Part of (b) follows from the corresponding
equivariance property of $\hat{\omega}_{\mathcal{W}}$, which was verified in
the proof of Lemma \ref{lem:classic}, together with the invariance property
in Equation (\ref{eqn:Jinv}) established in the proof of Lemma \ref%
{lem:vogel}. Part (d) follows from $\hat{\omega}_{\mathcal{W}}(y)\neq 0$ for 
$\lambda _{\mathbb{R}^{n}}$-almost every $y\in \mathbb{R}^{n}$ (cf. the
proof of Lemma \ref{lem:classic}) together with $\exp (cJ_{n,U}^{1}(y))>0$
for every $y\notin \limfunc{span}(X,U)$. The final claim follows from the
final statement in Lemma \ref{lem:auxA5} since (cf. the proof of Lemma \ref%
{lem:classic}) $\hat{\omega}_{\mathcal{W}}(y)\geq 0$ holds for every $y\in 
\mathbb{R}^{n}$ in case $\Pi _{\limfunc{span}(X)^{\bot }}\mathcal{W}\Pi _{%
\limfunc{span}(X)^{\bot }}$ is nonnegative definite.

2. The proof is very similar to the proof of the first part. It follows
along the same lines observing that the function defined via 
\begin{equation*}
y\mapsto \frac{\hat{u}^{\prime }(y)A^{\prime }A\hat{u}^{\prime }(y)}{\hat{u}%
^{\prime }(y)\hat{u}^{\prime }(y)}
\end{equation*}%
is well-defined and continuous on $\mathbb{R}^{n}\backslash \limfunc{span}%
(X) $, and is $G(\mathfrak{M})$-invariant. We skip the details. $%
\blacksquare $

\textbf{Proof of Corollary \ref{cor:vogelbunzel}:} Noting that 
\begin{equation*}
P_{0,I_{n}}(\check{\Omega}\text{ is nonnegative definite})=P_{0,I_{n}}(\hat{%
\omega}_{\mathcal{W}}\geq 0)
\end{equation*}%
in our present context, the first part follows upon combining Lemma \ref%
{lem:vogelbunzel} and Theorem \ref{size_bound_polynomial} (the statement
concerning the lower bound being 1 if $\Pi _{\limfunc{span}(X)^{\bot }}%
\mathcal{W}\Pi _{\limfunc{span}(X)^{\bot }}$ is nonnegative definite follows
from nonnegative definiteness of $\check{\Omega}_{\mathcal{W},U,c}^{\mathrm{%
BV},J}$, or of $\check{\Omega}_{\mathcal{W},c}^{\mathrm{BV}}$, respectively,
on the complement of $\lambda _{\mathbb{R}^{n}}$-null sets in that case).
For the last part of the corollary, we can apply a similar argument as the
one that was given to verify the analogous statement in Corollary \ref%
{cor:polylrv2}: Note that now $\{T\geq 0\}=\{\hat{\omega}_{\mathcal{W}}\geq
0\}\cup \{R\hat{\beta}=r\}\cup N$, where $N=\limfunc{span}(X,U)$ if $T$ is
based on $\check{\Omega}_{\mathcal{W},U,c}^{\mathrm{BV},J}$, and $N=\limfunc{%
span}(X)$ if $T$ is based on $\check{\Omega}_{\mathcal{W},c}^{\mathrm{BV}}$.
In both cases $N\cup \{R\hat{\beta}=r\}$ is a $\lambda _{\mathbb{R}^{n}}$%
-null set, and we see that (\ref{eqn:powerbound}) holds also in the
situation of the present lemma. The remainder of the proof is now analogous
to the argument given at the end of the proof of Corollary \ref{cor:polylrv2}%
. $\blacksquare $

\begin{lemma}
\label{lem:BVfeas} Let $a_{i}\in (0,\infty )$ for $i=0,\ldots ,m^{\prime }$ (%
$m^{\prime }\in \mathbb{N}$), $\bar{a}_{i}\in \mathbb{R}$ for $i=1,\ldots
,m^{\prime }$, $h_{i}\in \mathbb{R}$ for $i=0,\ldots ,m^{\prime \prime }$ ($%
m^{\prime \prime }\in \mathbb{N}$) with $h_{m^{\prime \prime }}\neq 0$, and $%
p_{i}\in \mathbb{R}$ for $i=0,\ldots ,m^{\prime \prime \prime }$ ($m^{\prime
\prime \prime }\in \mathbb{N}$) with $p_{m^{\prime \prime \prime }}\neq 0$.
Suppose further that $e_{n}(n)\notin \limfunc{span}(X)^{\bot }$. Then, $%
\tilde{N}=\limfunc{span}(X)$ (where $\tilde{N}$ has been defined in (\ref%
{eqn:excbv})), and the following holds:

\begin{enumerate}
\item If $U$ is an $n\times m$-dimensional matrix with $m\geq 1$ such that $%
(X,U)$ is of full column-rank $k+m<n$, then the estimators $\hat{\beta}$ and 
$\check{\Omega}_{U,a,\bar{a},h,p}^{\mathrm{BV},J}$ satisfy Assumption \ref%
{Ass_5} with $N=N_{\mathrm{BV},U}$, where $N_{\mathrm{BV},U}=\limfunc{span}%
(X,U)\cup \left\{ y\in \mathbb{R}^{n}\backslash \limfunc{span}(X,U):\hat{\rho%
}(y)\in \{\bar{a}_{1},\ldots ,\bar{a}_{m^{\prime }}\}\right\} $ in case $%
\hat{\rho}$ attains at least two different values on $\mathbb{R}%
^{n}\backslash \limfunc{span}(X)$, and $N_{\mathrm{BV},U}=\limfunc{span}%
(X,U) $ else. Furthermore, $\check{\Omega}_{U,a,\bar{a},h,p}^{\mathrm{BV},J}$
satisfies Assumption \ref{Ass_7}, and $\check{\Omega}_{U,a,\bar{a},h,p}^{%
\mathrm{BV},J}(y)$ is positive definite for every $y\in \mathbb{R}%
^{n}\backslash N_{\mathrm{BV},U}$ (in fact, for $y\in \mathbb{R}%
^{n}\backslash \limfunc{span}(X,U)$).

\item The estimators $\hat{\beta}$ and $\check{\Omega}_{a,\bar{a},h,p}^{%
\mathrm{BV}}$ satisfy Assumption \ref{Ass_5} with $N=N_{\mathrm{BV}}$, where 
$N_{\mathrm{BV}}=\limfunc{span}(X)\cup \left\{ y\in \mathbb{R}^{n}\backslash 
\limfunc{span}(X):\hat{\rho}(y)\in \{\bar{a}_{1},\ldots ,\bar{a}_{m^{\prime
}}\}\right\} $ in case $\hat{\rho}$ attains at least two different values on 
$\mathbb{R}^{n}\backslash \limfunc{span}(X)$, and $N_{\mathrm{BV}}=\limfunc{%
span}(X)$ else. Furthermore, $\check{\Omega}_{a,\bar{a},h,p}^{\mathrm{BV}}$
satisfies Assumption \ref{Ass_7}, and $\check{\Omega}_{a,\bar{a},h,p}^{%
\mathrm{BV}}(y)$ is positive definite for every $y\in \mathbb{R}%
^{n}\backslash N_{\mathrm{BV}}$ (in fact, for $y\in \mathbb{R}^{n}\backslash 
\limfunc{span}(X)$).
\end{enumerate}
\end{lemma}

\textbf{Proof:} The assumption $e_{n}(n)\notin \limfunc{span}(X)^{\bot }$
implies non-existence of a $y\in \mathbb{R}^{n}\backslash \limfunc{span}(X)$
so that $\sum_{i=1}^{n-1}\hat{u}_{i}^{2}(y)=0$, showing that $\hat{\rho}$ is
well-defined everywhere on $\mathbb{R}^{n}\backslash \limfunc{span}(X)$,
i.e., that $\tilde{N}=\limfunc{span}(X)$. We consider two cases: First,
assume that the design matrix $X$ is such that $\hat{\rho}=\rho $ holds
everywhere on $\mathbb{R}^{n}\backslash \limfunc{span}(X)$ for some fixed $%
\rho \in \mathbb{R}$. Then, the statements in 1. and 2., except for the
positive definiteness claims, follow from Lemma \ref{lem:vogelbunzel},
because $b_{\mathrm{BV}}(.,a,A)$ and $c_{\mathrm{BV}}(.,p)$ are then
constant equal to $b$ and $c$, say, respectively, on $\mathbb{R}%
^{n}\backslash \limfunc{span}(X)$ and thus $\check{\Omega}_{U,a,\bar{a}%
,h,p}^{\mathrm{BV},J}(y)=\check{\Omega}_{\mathcal{W},U,c}^{\mathrm{BV},J}(y)$
holds for every $y\notin \limfunc{span}(X,U)$, and $\check{\Omega}_{a,\bar{a}%
,h,p}^{\mathrm{BV}}(y)=\check{\Omega}_{\mathcal{W},c}^{\mathrm{BV}}(y)$
holds for every $y\notin \limfunc{span}(X)$ where the matrix $\mathcal{W}=(%
\mathcal{W}_{ij})=(\kappa _{D}(|i-j|/\max (bn,2)))$. Observe here that $%
\mathcal{W}$ is constant in $y$, is symmetric, and is positive definite. The
positive definiteness claims in 1. and 2. finally follow since $\hat{\omega}%
_{\mathcal{W}}(y)>0$ holds for $y\in \mathbb{R}^{n}\backslash \limfunc{span}%
(X)$ in view of positive definiteness of $\mathcal{W}$.

Next, we consider the case where $X$ is such that $\hat{\rho}$ attains at
least two different values on $\mathbb{R}^{n}\backslash \limfunc{span}(X)$.
We start with the statement in 1.: First of all, $N_{\mathrm{BV},U}$ is
easily seen to be $G(\mathfrak{M})$-invariant (because $\hat{\rho}:\mathbb{R}%
^{n}\backslash \limfunc{span}(X)\rightarrow \mathbb{R}$ is so). Second, we
can rewrite 
\begin{equation*}
N_{\mathrm{BV},U}=\bigcup_{i=1}^{m^{\prime }}\left\{ y\in \mathbb{R}%
^{n}:\sum_{i=2}^{n}\hat{u}_{i}(y)\hat{u}_{i-1}(y)-\bar{a}_{i}\sum_{i=1}^{n-1}%
\hat{u}_{i}^{2}(y)=0\right\} \cup \limfunc{span}(X,U).
\end{equation*}%
From that we see that $N_{\mathrm{BV},U}$ is a finite union of algebraic
sets, and hence an algebraic set. Thus, $N_{\mathrm{BV},U}$ is closed. Since
we also work under the hypothesis that $\hat{\rho}$ attains at least two
different values on $\mathbb{R}^{n}\backslash \limfunc{span}(X)$, we can
conclude that 
\begin{equation*}
\left\{ y\in \mathbb{R}^{n}:\sum_{i=2}^{n}\hat{u}_{i}(y)\hat{u}_{i-1}(y)-%
\bar{a}_{i}\sum_{i=1}^{n-1}\hat{u}_{i}^{2}(y)=0\right\} \neq \mathbb{R}^{n}
\end{equation*}%
holds for every $i=1,\ldots ,m^{\prime }$. It follows that the algebraic set
in the previous display is a $\lambda _{\mathbb{R}^{n}}$-null set for every $%
i=1,\ldots ,m^{\prime }$. Hence $N_{\mathrm{BV},U}$ is a closed $\lambda _{%
\mathbb{R}^{n}}$-null set as $\limfunc{span}(X,U)\neq \mathbb{R}^{n}$. To
prove the statements of 1., we now verify (a)-(d) in Lemma \ref{lem:auxA5}
for $H=I_{n}$, $\nu (.)=\hat{\omega}_{\mathcal{W}_{\mathrm{BV}}}(.)\exp (c_{%
\mathrm{BV}}(.,p)J_{n,U}^{1}(.))$, and $N^{\prime }=N_{\mathrm{BV},U}$. We
have already verified (a). Furthermore, note that $b_{\mathrm{BV}}(.,a,\bar{a%
})$ is continuous on $\mathbb{R}^{n}\backslash N_{\mathrm{BV},U}$. As a
consequence, $c_{\mathrm{BV}}(.,p)$, and $\mathcal{W}_{\mathrm{BV}}(.)$, and
thus $\hat{\omega}_{\mathcal{W}_{\mathrm{BV}}}$ are continuous on $\mathbb{R}%
^{n}\backslash N_{\mathrm{BV},U}$. We already know from the proof of Lemma %
\ref{lem:vogel} that $J_{n,U}^{1}$ is continuous on the complement of $%
\limfunc{span}(X,U)\subseteq N_{\mathrm{BV},U}$. It thus follows that $%
y\mapsto \hat{\omega}_{\mathcal{W}_{\mathrm{BV}}}(y)\exp (c_{\mathrm{BV}%
}(y,p)J_{n,U}^{1}(y))$ is continuous on $\mathbb{R}^{n}\backslash N_{\mathrm{%
BV},U}$. Hence, we have verified (c) in Lemma \ref{lem:auxA5}. To verify (b)
we recall from above that $N_{\mathrm{BV},U}$ is $G(\mathfrak{M})$%
-invariant. Furthermore, the required equivariance property in (b) holds as
a consequence of $G(\mathfrak{M})$-invariance of $\hat{\rho}$ and $%
J_{n,U}^{1}$ (cf. (\ref{eqn:Jinv})), and hence of $c_{\mathrm{BV}}(.,p)$ and 
$\mathcal{W}_{\mathrm{BV}}(.)$ on $\mathbb{R}^{n}\backslash N_{\mathrm{BV}%
,U} $, together with $\hat{u}(\delta y+X\eta )=\delta \hat{u}(y)$ for every $%
\delta \neq 0$, $y\in \mathbb{R}^{n}$ and $\eta \in \mathbb{R}^{k}$. That $%
\nu (y)=\hat{\omega}_{\mathcal{W}_{\mathrm{BV}}}(y)\exp (c_{\mathrm{BV}%
}(y,p)J_{n,U}^{1}(y))$ is even positive on $\mathbb{R}^{n}\backslash N_{%
\mathrm{BV},U}$ follows because $\mathcal{W}_{\mathrm{BV}}(y)$ is a positive
definite matrix for every $y\in \mathbb{R}^{n}\backslash \tilde{N}$ and $%
\tilde{N}=\limfunc{span}(X)$ holds. This implies (d) in Lemma \ref{lem:auxA5}%
, and also the sufficient condition for positive definiteness in the same
lemma. The statements in 2. for the case where $\hat{\rho}$ attains at least
two different values on $\mathbb{R}^{n}\backslash \limfunc{span}(X)$ are
almost identical, and we skip the details. $\blacksquare $

\textbf{Proof of Corollary \ref{cor:vogelbunzelfeasible}:} From Assumption %
\ref{as:poly} it follows that the last row of $X$ is not equal to zero,
i.e., $e_{n}(n)\notin \limfunc{span}(X)^{\bot }$ must hold. Hence, all
assumptions of Lemma \ref{lem:BVfeas} are satisfied. Combining this lemma
with Theorem \ref{size_bound_polynomial} proves the claims with $C_{\mathrm{%
BV}}(y,h)$ replaced by an arbitrary constant critical value $C$ (noting that
the lower bound obtained via Theorem \ref{size_bound_polynomial} equals $1$
due to nonnegative definiteness of $\check{\Omega}_{U,a,\bar{a},h,p}^{%
\mathrm{BV},J}$ ($\check{\Omega}_{a,\bar{a},h,p}^{\mathrm{BV}}$,
respectively) on the complement of the $\lambda _{\mathbb{R}^{n}}$-null set $%
N_{\mathrm{BV},U}$ ($N_{\mathrm{BV}}$, respectively)). But now we observe
that $y\mapsto C_{\mathrm{BV}}(y,h)$ is well-defined on $\mathbb{R}^{n}$
(recall the convention preceding Corollary \ref{cor:vogelbunzelfeasible}),
and by construction takes on only finitely many real numbers $C_{1}<\ldots
<C_{l}$, say. Hence, for every $f\in \mathfrak{F}$, every $\mu _{0}\in 
\mathfrak{M}_{0}$, every $\sigma ^{2}\in (0,\infty )$ we can conclude that 
\begin{equation*}
P_{\mu _{0},\sigma ^{2}\Sigma (f)}(\{y\in \mathbb{R}^{n}:T(y)\geq C_{\mathrm{%
BV}}(y,h)\})\geq P_{\mu _{0},\sigma ^{2}\Sigma (f)}(\{y\in \mathbb{R}%
^{n}:T(y)\geq C_{l}\}).
\end{equation*}%
Now apply what has been established before with $C=C_{l}$. This completes
the proof. $\blacksquare $

\textbf{Proof of Theorem \ref{size_one_cyclical}:} For any $i=1,2$ with $%
R_{\cdot i}\neq 0$, the $i$-th column of $E_{n,0}(\omega )$ does not belong
to $\mathfrak{M}_{0}^{lin}$. Hence $\func{span}(E_{n,0}(\omega ))\nsubseteq 
\mathfrak{M}_{0}^{lin}$, implying that $\rho (\omega )$ must be zero.
However, $\func{span}(E_{n,0}(\omega )\subseteq \func{span}(X)$ clearly
holds. All the other assumptions being obviously satisfied, Theorem \ref%
{size_one_AR2} completes the proof. $\blacksquare $

\textbf{Proof of Theorem \ref{size_bound_cyclical}:} We apply Theorem \ref%
{theo:AR2+}. It suffices to verify that $\gamma =\omega $ satisfies the
assumption $\limfunc{span}(E_{n,\rho (\gamma )}(\gamma ))\subseteq \mathrm{%
\limfunc{span}}(X)$ in that theorem. But this can be established exactly in
the same way as in the proof of Theorem \ref{size_one_cyclical}. $%
\blacksquare $

\bibliographystyle{ims}
\bibliography{refs}

\end{document}